\documentclass{gtart_h}  

\def\ifplaintex{\expandafter\ifx\csname documentclass\endcsname\relax}

\ifplaintex 
\else
\fi

\expandafter\ifx\csname beginpicture\endcsname\relax
\expandafter\ifx\csname documentclass\endcsname\relax
\input pictex \else
\input prepictex \input pictex \input postpictex \fi\fi

\def\gt{{\mathsurround=0pt\it $\cal G\mskip-2mu$eometry \&\ 
$\cal T\!\!$opology}}        

\def\gtp{{\mathsurround=0pt\it $\cal G\mskip-2mu$eometry \&\ 
$\cal T\!\!$opology $\cal P\!$ublications}}  

\def\lognumber#1{\def\thelognumber{#1}}
\def\volumenumber#1{\def\thevolumenumber{#1}}
\def\papernumber#1{\def\thepapernumber{#1}}
\def\volumeyear#1{\def\thevolumeyear{#1}}

\def\pagenumbers#1#2{\def\startpage{#1}\def\finishpage{#2}}
\def\published#1{\def\publishdate{#1}}
\def\proposed#1{\def\theproposer{#1}}
\def\seconded#1{\def\theseconders{#1}}
\def\received#1{\def\receiveddate{#1}}
\def\revised#1{\def\reviseddate{#1}}
\def\accepted#1{\def\accepteddate{#1}}
\def\asciititle#1{\def\theasciititle{#1}}
\def\covertitle#1{\def\thecovertitle{#1}}
\def\coverauthors#1{\def\thecoverauthors{#1}}
\def\asciiauthors#1{\def\theasciiauthors{#1}}
\def\asciiaddress#1{\def\theasciiaddress{#1}}
\def\asciiemail#1{\def\theasciiemail{#1}}

\long\def\asciiabstract#1{\long\def\theasciiabstract{#1}}
\def\asciikeywords#1{\def\theasciikeywords{#1}}

\def\shorttitle#1{\def\theshorttitle{#1}}

\let\\\par\let\thelognumber\relax
\let\thevolumenumber\relax\let\thepapernumber\relax
\let\thevolumeyear\relax\let\thesamplenumber\relax\let\startpage\relax
\let\finishpage\relax\let\publishdate\relax\let\receiveddate\relax
\let\reviseddate\relax\let\accepteddate\relax\let\theasciititle\relax
\let\thecovertitle\relax\let\theasciiauthors\relax\let\theasciiaddress\relax
\let\theasciiabstract\relax\let\theasciikeywords\relax
\let\theasciiemail\relax\let\theshortauthors\relax\let\theshorttitle\relax
\let\thecoverauthors\relax

\long\def\maketitlep{   

\count0=\startpage

\gt\hfill      
\beginpicture
\setcoordinatesystem units <0.33truein, 0.33truein> point at 2.2 0.9
\setplotsymbol ({$\cal G$})
\plotsymbolspacing=9truept
\circulararc 315 degrees from 0 1 center at 0 0
\setplotsymbol ({$\cal T$})
\circulararc 315 degrees from 1 -1 center at 1 0
\endpicture
\break
{\small\ifx\thesamplenumber\relax 
Volume \else Sample
\fi\thevolumenumber\ (\thevolumeyear)
\startpage--\finishpage\nl
Published: \publishdate}
\vglue 0.5truein plus 0.4fil minus 0.1truein

{\parskip=0pt\leftskip 0pt plus 1fil\def\\{\par\smallskip}{\ifplaintex\large
\else\Large\fi\bf\thetitle}\par\medskip}   

\vglue 0pt plus 0.1fil 

{\parskip=0pt\leftskip 0pt plus 1fil\def\\{\par}{\sc\theauthors}
\par\medskip}

\vglue 0pt plus 0.1fil 

{\small\parskip=0pt\let\newline\\
{\leftskip 0pt plus 1fil\def\\{\par}{\sl\theaddress}\par}
\expandafter\ifx\theemail\relax    
\relax\else\vglue 5pt plus 0.02fil minus 2pt\def\\{\stdspace{\rm 
and}\stdspace} 
\cl{Email:\stdspace\tt\theemail}\fi
\ifx\theurl\relax                  
\relax\else\vglue 5pt plus 0.02fil minus 2pt\def\\{\stdspace{\rm 
and}\stdspace}
\cl{URL:\stdspace\tt\theurl}\fi\par}

\vglue 7pt plus 0.3fil minus 3pt

{\bf Abstract}
\vglue 5pt plus 0.1fil minus 2pt

\theabstract

\vglue 7pt plus 0.3fil minus 3pt

{\bf AMS Classification numbers}\quad Primary:\quad \theprimaryclass

Secondary:\quad \thesecondaryclass

\vglue 5pt plus 0.3fil minus 2pt

{\bf Keywords:}\quad \thekeywords

\vglue 10pt plus 0.5fil minus 5pt

{\small  Proposed: \theproposer\hfill Received: \receiveddate\nl
Seconded: \theseconders\hfill 
\ifx\reviseddate\relax                         
Accepted: \accepteddate                        
\else
Revised: \reviseddate                          
\fi}
\eject
}       

\let\maketitlepage\maketitlep
\let\maketitle\maketitlepage

\font\phead=cmsl9 scaled 950
\font=cmsl9 scaled 1050
\font\pnum=cmbx10 scaled 913
\font\lnum=cmbx10 
\font\pfoot=cmsl9 scaled 950
\font=cmsl9 scaled 1050
\ifplaintex
\headline{\vbox to 0pt{\vskip -4.5mm\line{\small\phead\ifnum
\count0=\startpage ISSN 1364-0380 (on line)
1465-3060 (printed) \hfill {\pnum\folio}\else\ifodd\count0\def\\{ }%
\ifx\theshorttitle\relax\thetitle\else\theshorttitle\fi\hfill{\pnum\folio}
\else\def\\{ and }{\pnum\folio}\hfill\ifx\theshortauthors\relax\theauthors
\else\theshortauthors\fi\fi\fi}\vss}}
\footline{\vbox to 0pt{\vglue 0mm\line{\small\pfoot\ifnum\count0=\startpage
\copyright\ \gtp\hfill\else
\gt, Volume \thevolumenumber\ (\thevolumeyear)\hfill\fi}\vss
}}
\else
\makeatletter
\def\@oddhead{{\small\lhead\ifnum\count0=\startpage ISSN 1364-0380 (on line)
1465-3060 (printed) \hfill {\lnum\number\count0}\else\ifodd\count0
\def\\{ }\ifx\theshorttitle\relax \thetitle \else\theshorttitle\fi\hfill
{\lnum\number\count0}\else\def\\{ and }{\lnum\number\count0}
\hfill\ifx\theshortauthors\relax 
\theauthors\else\theshortauthors\fi\fi\fi}}\def\@evenhead{\@oddhead}
\def\@oddfoot{\small\lfoot\ifnum\count0=\startpage\copyright\ \gtp\hfill\else
\gt, Volume \thevolumenumber\ (\thevolumeyear)\hfill\fi}
\def\@evenfoot{\@oddfoot}
\makeatother
\fi

\newwrite\gtoutfile
\long\gdef\makeheadfile{  
{\def\\{, }\def\s{ }
\immediate\openout\gtoutfile head.xxx
\immediate\write\gtoutfile{Proxy-for: \ifx\theasciiauthors\relax
\theauthors\else\theasciiauthors\fi\s<\ifx\theasciiemail\relax\theemail\else\theasciiemail\fi>}
\immediate\write\gtoutfile{\noexpand\\}
\immediate\write\gtoutfile{Authors: \ifx\theasciiauthors\relax
\theauthors\else\theasciiauthors\fi}
{\def\\{ }\immediate\write\gtoutfile{Title: \ifx\theasciititle\relax
\thetitle\else\theasciititle\fi}}
\immediate\write\gtoutfile{Subj-class: GT or SG or MG etc}
\immediate\write\gtoutfile{MSC-class: \theprimaryclass\ifx\thesecondaryclass\relax\else, \thesecondaryclass\fi}
\immediate\write\gtoutfile{Journal-ref: Geom. Topol. \thevolumenumber
(\thevolumeyear) \startpage-\finishpage}
\immediate\write\gtoutfile{Comments: Published by Geometry and Topology at}
\immediate\write\gtoutfile{\s\s http://www.maths.warwick.ac.uk/gt/GTVol\thevolumenumber/paper\thepapernumber.abs.html}
\immediate\write\gtoutfile{\noexpand\\}
\immediate\write\gtoutfile{}
\ifx\theasciiabstract\relax
\immediate\write\gtoutfile{\theabstract}\else
\immediate\write\gtoutfile{\theasciiabstract}\fi
\immediate\write\gtoutfile{}
\immediate\write\gtoutfile{\noexpand\\}
\immediate\write\gtoutfile{}
\immediate\closeout\gtoutfile}}  

\def\maketitlepage{\maketitlep\makeheadfile}
\let\maketitle\maketitlepage

\lognumber{349}
\received{2 August 2003}
\volumenumber{9}\papernumber{15}\volumeyear{2005}
\pagenumbers{493}{569}   
\revised{5 April 2005}
\published{8 April 2005}
\accepted{5 April 2005}
\proposed{Robion Kirby}
\seconded{Walter Neumann, Shigeyuki Morita}

\usepackage{amssymb,amsmath,graphicx,psfrag}
\def\psfraga <#1,#2> #3#4{%
\psfrag {#3}{\smash{\rlap{\kern #1 \raise #2\hbox{#4}}}}}

\everymath{\displaystyle}
\allowdisplaybreaks
\def\S{Section }
\def\MR#1{\ignorespaces\qua\href{http://www.ams.org/mathscinet-getitem?mr=#1}{MR#1}}

\newtheorem{teo}{Theorem}[section]
\newtheorem{lem}[teo]{Lemma}
\newtheorem{cor}[teo]{Corollary}

\newtheorem{prop}[teo]{Proposition}
\newtheorem{defi}[teo]{Definition}
\newtheorem{ques}[teo]{Question}
\newtheorem{conj}[teo]{Conjecture}
\newtheorem{remark}[teo]{Remark}
\newtheorem{remarks}[teo]{Remarks}

\newcommand{\mr}{\mathbb{R}}
\newcommand{\mc}{\mathbb{C}}
\newcommand{\mz}{\mathbb{Z}}
\newcommand{\mh}{\mathbb{H}}

\newcommand{\mn}{\mathbb{N}}

\newcommand{\Aa}{{\mathcal A}}
\newcommand{\Bb}{{\mathcal B}}

\newcommand{\Dd}{{\mathcal D}}

\newcommand{\Ii}{{\mathcal I}}

\newcommand{\Ll}{{\mathcal L}}

\newcommand{\Pp}{{\mathcal P}}
\newcommand{\Rr}{{\mathcal R}}

\newcommand{\Tt}{{\mathcal T}}

\newcommand{\cG}{{\mathfrak c}}

\newcommand{\C}{{\mathbb C}}

\newcommand{\Z}{{\mathbb Z}}
\newcommand{\R}{{\mathbb R}}

\begin{document}

\title{Classical and quantum dilogarithmic invariants of flat $PSL(2,\C)$--bundles over 3--manifolds}
\asciititle{Classical and quantum dilogarithmic invariants of flat PSL(2,C)-bundles over 3-manifolds}
\covertitle{Classical and quantum dilogarithmic invariants of flat 
$PSL(2,{\noexpand\bf C})$--bundles over 3--manifolds}
\shorttitle{Dilogarithmic invariants of $PSL(2,\C)$--bundles over 3--manifolds}

\author{St\'ephane Baseilhac\\Riccardo Benedetti}
\asciiauthors{Stephane Baseilhac and Riccardo Benedetti}
\coverauthors{St\noexpand\'ephane Baseilhac\\Riccardo Benedetti}

\address{Universit\'e de Grenoble I, Institut Joseph Fourier, UMR CNRS
5582\\100 rue des Maths, B.P. 74, F-38402 Saint-Martin-d'H\`eres Cedex, FRANCE}
\secondaddress{Dipartimento di Matematica, Universit\`a di Pisa\\Via
F. Buonarroti, 2, I-56127 Pisa, ITALY}

\asciiaddress{Universite de Grenoble I, Institut Joseph Fourier, UMR CNRS
5582\\100 rue des Maths, B.P. 74, F-38402 Saint-Martin-d'Heres Cedex, 
FRANCE\\and\\Dipartimento di Matematica, Universita di Pisa\\Via
F. Buonarroti, 2, I-56127 Pisa, ITALY}

\gtemail{\mailto{baseilha@ujf-grenoble.fr}{\rm\qua 
and\qua}\mailto{benedett@dm.unipi.it}}
\asciiemail{baseilha@ujf-grenoble.fr, benedett@dm.unipi.it}

\begin{abstract}
We introduce a family of {\it matrix dilogarithms}, which
  are automorphisms of $\mc^N \otimes \mc^N$, $N$ being any odd
  positive integer, associated to hyperbolic ideal tetrahedra equipped
  with an additional decoration. The matrix dilogarithms satisfy
  fundamental {\it five-term identities} that correspond to decorated
  versions of the $2 \to 3$ move on $3$--dimensional triangulations.
  Together with the decoration, they arise from the solution we give
  of a {\it symmetrization problem} for a specific family of {\it
    basic} matrix dilogarithms, the classical ($N=1$) one being the
  Rogers dilogarithm, which only satisfy one special instance of
  five-term identity. We use the matrix dilogarithms to construct
  invariant state sums for closed oriented $3$--manifolds $W$ endowed
  with a flat principal $PSL(2,\mc)$--bundle $\rho$, and a fixed non
  empty link $L$ if $N>1$, and for (possibly ``marked'') cusped
  hyperbolic $3$--manifolds $M$. When $N=1$ the state sums recover
  known simplicial formulas for the volume and the Chern--Simons
  invariant. When $N>1$, the invariants for $M$ are new; those for
  triples $(W,L,\rho)$ coincide with the quantum hyperbolic invariants
  defined in \cite{BB2}, though our present approach clarifies
  substantially their nature. We analyse the structural coincidences
  versus discrepancies between the cases $N=1$ and $N>1$, and we
  formulate ``Volume Conjectures'', having geometric motivations,
  about the asymptotic behaviour of the invariants when $N\rightarrow
  \infty$.
\end{abstract}

\asciiabstract{%
We introduce a family of matrix dilogarithms, which are automorphisms
of C^N tensor C^N, N being any odd positive integer, associated to
hyperbolic ideal tetrahedra equipped with an additional
decoration. The matrix dilogarithms satisfy fundamental five-term
identities that correspond to decorated versions of the 2 --> 3 move
on 3-dimensional triangulations.  Together with the decoration, they
arise from the solution we give of a symmetrization problem for a
specific family of basic matrix dilogarithms, the classical (N=1) one
being the Rogers dilogarithm, which only satisfy one special instance
of five-term identity. We use the matrix dilogarithms to construct
invariant state sums for closed oriented 3-manifolds $W$ endowed with
a flat principal PSL(2,C)-bundle rho, and a fixed non empty link L if
N>1, and for (possibly "marked") cusped hyperbolic 3-manifolds M. When
N=1 the state sums recover known simplicial formulas for the volume
and the Chern-Simons invariant. When N>$, the invariants for M are
new; those for triples (W,L,rho) coincide with the quantum hyperbolic
invariants defined in [Topology 43 (2004) 1373-1423], though our
present approach clarifies substantially their nature. We analyse the
structural coincidences versus discrepancies between the cases N=1 and
N>1, and we formulate "Volume Conjectures", having geometric
motivations, about the asymptotic behaviour of the invariants when N
tends to infinity.}




\asciikeywords{Dilogarithms, state sum invariants,
quantum field theory, Cheeger-Chern-Simons invariants, scissors congruences,
hyperbolic 3-manifolds.}

{\small\maketitlepage}

\section{Introduction} \label{GTINTRO}

Since its beginning in the eighties, the theory of quantum invariants
of links and $3$--manifolds has rapidly grown up as a very active
domain of research with a large interaction between quite seemingly
independent branchs of mathematics and ideas from Quantum Field
Theories (QFT) in Physics. They are now organized in a well-structured
machinery based on the theory of representations of quantum groups
and, more generally, of linear monoidal categories, which is
recognized as a very powerful tool for producing `exact'
$3$--dimensional QFT (ie, functors from categories of manifold
cobordisms towards such linear categories). For reviews, see
eg \cite{O,T}. These exact theories also provide a new predictive
power and meaningful framework for the physics ideas they were
inspired from.

Nevertheless, in spite of its success and very aesthetic formalism,
$3$--dimensional `Quantum Topology' followed until recently a rather
divergent path with respect to more classical topological and
geometric themes, which were mostly developed during the last decades
into Thurston's geometrization program.  A conceptual breakthrough was
done by Kashaev with his Volume Conjecture \cite{K2}. He derived from
a family $\{K_N\}$, $N>1$ being any odd integer number, of conjectural
complex valued topological invariants of links in arbitrary closed
oriented $3$--manifolds, a well-defined family $\{\langle L \rangle_N\}$
of invariants of links $L$ in $S^3$ \cite{K3}, later identified by
Murakami--Murakami as the values of specific coloured Jones polynomials
at the roots of unity $\exp(2i\pi/N)$ \cite{MM}. He predicted also
that if $L$ is a hyperbolic link, then the asymptotic behaviour of
$\langle L \rangle_N$ when $N \rightarrow \infty$ recovers the volume
of the complement of $L$. The main motivation for this conjecture is
that the asymptotic behaviour of the elementary building blocks of
$K_N$ essentially involve classical dilogarithm functions, which are
known to be related to the computation of the volume of hyperbolic
polyhedra.

In our previous paper \cite{BB2} we constructed \emph{quantum
  hyperbolic invariants} (QHI) $H_N(W,L,\rho)$, well-defined possibly
up to a sign and an $N$th root of unity phase factor. Here $N > 1$ is
an odd integer, $L$ is a non empty link in a closed oriented
$3$--manifold $W$, and $\rho$ is a flat principal $PSL(2,\mc)$--bundle
over $W$.  These invariants eventually incorporate as a particular
case the Kashaev's conjectural ones, by using the trivial flat bundle.
A main ingredient of our construction was the use of so called {\it
  decorated $\Ii$--triangulations}, which are particular structured
families of oriented hyperbolic ideal tetrahedra with ordered
vertices, encoded by their triples of cross-ratio moduli, and equipped
with some additional decoration. Each QHI $H_N(W,L,\rho)$ can be
expressed as a state sum, ie, the total contraction of a pattern of
special automorphisms of $\mc^N \otimes \mc^N$ associated to the
tetrahedra of any such a decorated $\Ii$--triangulation. However, our
understanding of this remarkable family of tensors, in particular of
the nature of the decoration of the tetrahedra entering their
definition, was not satisfactory (see Remark \ref{compareQHI} for
further comments on this point).  As a consequence, also the nature of
the QHI remained somewhat obscure.  

The first aim of the present paper is to unfold and clarify the
structure of these tensors, called here `quantum' {\it matrix
  dilogarithms}. We formalize them and we show their fundamental
properties. They are explicitely given automorphisms of $\mc^N \otimes
\mc^N$, $N>1$ being any odd positive integer, associated to decorated
$\Ii$--tetrahedra. Their main structural property consists in
satisfying fundamental {\it five-term identities}: for every instance,
called a {\it transit}, of an $\Ii$--decorated version of the $2 \to 3$
bistellar (sometimes called Pachner, or Matveev--Piergallini) local
move on $3$--dimensional triangulations, the contractions of the two
patterns of associated matrix dilogarithms eventually lead to the same
tensor up to a determined phase ambiguity. The matrix dilogarithms, as
well as the additional decoration on the associated $\Ii$--tetrahedra,
arise from the solution of a {\it symmetrization problem} for a
specific family of {\it basic} matrix dilogarithms. These are derived
from the $6j$--symbols for the cyclic representation theory of a Borel
quantum subalgebra $\Bb_\zeta$ of $U_\zeta(sl(2,\mc))$, where
$\zeta=\exp(2i\pi/N)$. They satisfy only one particular instance of
five-term identity, the {\it Schaeffer's identity}, with some
geometric constraints on the involved cross-ratio moduli. The basic
matrix dilogarithms can be considered as natural non-commutative
analogues of the classical Rogers dilogarithm.  To stress this point,
our analysis of the symmetrization problem runs parallel to that for
the {\it classical} case ($N=1$), where we take the exponential of the
classical Rogers dilogarithm as basic dilogarithm. The symmetrization
problem makes the technical core of what we call the {\it semi-local}
part of the paper.

Later we face the {\it global} problems that arise in constructing
classical and quantum dilogarithmic invariants based on
globally decorated $\Ii$--triangulations of triples
$(W,L,\rho)$ or oriented non
compact complete hyperbolic $3$--manifolds $M$ of finite volume (for
short: {\it cusped} manifolds). 

For triples $(W,L,\rho)$, in the classical case the link is actually
immaterial; the dilogarithmic invariant only depends on the pair
$(W,\rho)$ and recovers the volume and the Chern--Simons invariant of
$\rho$. On the other hand, in the quantum case it is necessary to
incorporate a non empty (arbitrary) ``link fixing'' in the whole
construction, and the invariants are sensitive to the link. They
coincide with the QHI, but the present semi-local analysis
substantially clarifies their nature.

For cusped manifolds $M$, in the classical case the dilogarithmic
invariant recovers known simplicial formulas for the volume Vol$(M)$
and the Chern--Simons invariant CS$(M)$ \cite {NZ,N3,NY} (see also the
recent paper \cite{N4}). In the quantum case, the invariants $H_N(M)$
are new. Their construction is clean for the wide class of so called
``weakly-gentle'' cusped manifolds (see Definition \ref{W-T}); for
general cusped manifolds it is more tricky.  For weakly-gentle cusped
manifolds $M$, we recognize a strong structural coincidence between
the classical and the quantum invariants. Both are defined on the same
geometric objects, whereas to construct the quantum invariants for
general cusped manifolds we have to incorporate systems of arcs that
play the role of the link in the closed manifold case. (Note, however,
that it is reasonable to ask whether every cusped manifold is
weakly-gentle.) This leads us to formulate a version of the Volume
Conjecture for weakly-gentle cusped manifolds, relating the asymptotic
behaviour of $H_N(M)$ when $N\to \infty$ to ${\rm CS}(M) + i {\rm
  Vol}(M)$. Other forms of the Volume Conjecture (for example related
to Thurston's hyperbolic Dehn filling theorem) having substantial
geometric motivations are also proposed.

By using the matrix dilogarithms as fundamental ingredients, we have
developed in \cite{BB3} a family of exact finite dimensional {\it
  quantum hyperbolic field theories} (QHFT). The QHFT are
representations in the tensorial category of complex linear spaces of
a suitable 2+1 bordism category, based on arbitrary compact oriented
$3$--manifolds equipped with properly embedded tangles and with flat
principal $PSL(2,\C)$--bundles having arbitrary holonomy at the
meridians of the tangle components. The QHFT incorporate the present
dilogarithmic invariants as instances of partition functions.

There is a wide literature about the classical
Rogers dilogarithm and the computation of the volume and the
Chern--Simons invariant of $3$--manifolds equipped with flat
$PSL(2,\C)$--bundles. In particular, W Neumann's work
\cite{N3,N1,N2} on this subject has been a fundamental reference and a
source of inspiration for us.

\medskip

{\bf Acknowledgement}\qua We thank the referee for his remarks and
suggestions, that considerably improved the exposition of the paper.

\subsection{Description of the paper}

In Section \ref{GTSTATEMENT} we provide the complete statements of our
main results; in order to do it, we introduce the necessary
apparatus of notions and definitions. This is rather complicated
indeed, as it reflects the highly non trivial structure of the matrix
dilogarithms and the dilogarithmic invariants. This section is also
intended as a sort of self-contained account, without proofs, of the
content of the paper. For a deeper understanding, the reader is
addressed to the subsequent more technical sections.

In Sections \ref{GTBASICMD}, we introduce the {\it basic} matrix dilogarithms
$\Ll_N$ for every odd integer $N\geq 1$. The necessary quantum
algebraic background, in particular the derivation of $\Ll_N$, $N >
1$, from the representation theory of the Borel quantum algebra
$\Bb_\zeta$, shall be recalled in the Appendix. We formulate the {\it
  symmetrization problem}, which roughly asks to modify the basic
dilogarithms so as to make them `transit invariant', satisfying the
whole set of five-term relations. In Section \ref{CQDSYMROGERS} and
Section \ref{CQDSYMQUANTUM} we derive the essentially unique
solution of this problem, and this leads to the final matrix
dilogarithms $\mathcal{R}_N$, with their complicated additional
decoration. We show also (Lemma \ref{idemQHI}) that the
$\mathcal{R}_N$ coincide with the symmetrized quantum dilogarithms
used in \cite{BB2}, which implies that the quantum dilogarithmic invariants
of triples $(W,L,\rho)$ considered in Section \ref{GTINV} coincide with
the QHI. One aim of Section \ref{CQDSYMROGERS} is to provide, in
the simpler case $N=1$, a model for the contructions we need in the
quantum case.

In Section \ref{GTINV} we construct and discuss the classical and
quantum dilogarithmic invariants for triples $(W,L,\rho)$
and cusped manifolds $M$. 

In Section \ref{CQDSCISSORS} we construct further invariants called
{\it scissors congruence classes}. The terminology intentionally
refers to the background of the 3rd Hilbert problem (see \cite{DS,D2}
and \cite{N1}). The analysis of the relationship with the
dilogarithmic invariants is useful to settle out further discrepancies
between the classical and quantum cases and to formulate reasonable
intermediate questions towards the Volume Conjectures, that we recall
at the end of the paper.

\section{Statements of the main results}\label{GTSTATEMENT}

In this section we give the complete statements of the main results,
providing the necessary concepts and definitions. First we treat the
{\it semi-local} theory of matrix dilogarithms. Next we consider the construction of
invariant dilogarithmic state sums based on {\it globally} decorated
$\Ii$--triangulations of $3$--manifolds.

\subsection{Matrix dilogarithms and transit invariance}\label{LOCSTAT}

\subsubsection{Flat-charged $\Ii$--tetrahedra}

On the geometric/combinatorial side, the basic building blocks of
our constructions are the so called {\it flat/charged
$\Ii$--tetrahedra} that we are going to define.

An $\Ii$--{\it tetrahedron} $(\Delta,b,w)$ (see also \cite{BB2})
consists of

\medskip

(1)\qua An {\it oriented} tetrahedron $\Delta$, that we usually represent
    as positively embedded in $\R^3$ (oriented by its standard basis).

\medskip

(2)\qua A {\it branching} $b$ on $\Delta$, that is a choice of edge
 orientation associated to a total ordering $v_0,v_1,v_2,v_3$ of the
 vertices by the rule: each edge is oriented by the arrow emanating
 from the smallest endpoint. Denote by $E(\Delta)$ the set of
 $b$--oriented edges of $\Delta$, and by $e'$ the edge opposite to
 $e$. We put $e_0=[v_0,v_1]$, $e_1=[v_1,v_2]$ and
 $e_2=[v_0,v_2]=-[v_2,v_0]$. These are the edges of the face opposite
to the vertex $v_3$.

\medskip

(3)\qua A {\it modular triple}, $w=(w_0,w_1,w_2)=(w(e_0),w(e_1),w(e_2)) \in
(\C \setminus \{0,\ 1 \})^3$ such that
(indices mod($\mz/3\mz$)):
$$w_{j+1} = 1/(1-w_j)$$ hence
$$w_0w_1w_2 = -1 \ .$$
\noindent This gives a \emph{cross-ratio modulus}
$w(e)$ to each edge $e$ of $\Delta$, by imposing that $w(e)=w(e')$
for each edge $e$.

We say that $w$ is non degenerate if the imaginary parts of the
$w_j$ are not equal to zero; in such a case these imaginary parts
share the same sign $*_w = \pm 1$.

\paragraph{Complements on $\Ii$--tetrahedra}
The ordered triple of edges
$$(e_0=[v_0,v_1],e_2=[v_0,v_2],e_1'=[v_0,v_3])$$
departing from $v_0$ defines a \emph{$b$--orientation} of
$\Delta$. This orientation may or may not agree with the given
orientation of $\Delta$. In the first case we say that $b$ is of index
$*_b=1$, and it is of index $*_b=-1$ otherwise.

The $2$--faces of $\Delta$ can be named and ordered by their opposite
vertices.  For each $j=0,\ldots,3$ there are exactly $j$ $b$--oriented
edges incoming at the vertex $v_j$; hence there are only one source
and one sink of the branching. For any $2$--face $f$ of $\Delta$
the boundary of $f$ is not coherently oriented, only two edges of $f$
have a compatible {\it prevailing} orientation. In fact, each $2$--face
$f$ has two orientations; one is the boundary orientation induced by
the orientation of $\Delta$, via the convention ``last the ingoing
normal''; on the other hand, there is the {\it $b$--orientation}, that is the
orientation of $f$ which induces on $\partial f$ the prevailing
orientation among the three $b$--oriented edges. Remark that the
boundary and $b$--orientations coincide on exactly two $2$--faces
of $\Delta$.

Consider the half space model of the hyperbolic space $\mh^3$. We
orient it as an open set of $\R^3$. The natural boundary $\partial
\bar{\mh}^3=\mc\mathbb{P}^1 = \C \cup\{ \infty \}$ of $\mh^3$ is
oriented by its complex structure. We realize $PSL(2,\C)$ as the group
of orientation preserving (ie `direct') isometries of $\mh^3$, with
the corresponding conformal action on $\mc\mathbb{P}^1$.  Up to direct
isometry, an $\Ii$--tetrahedron $(\Delta,b,w)$ can be realized as an
hyperbolic ideal tetrahedron with 4 distinct $b$--ordered vertices
$u_0,u_1,u_2,u_3$ on $\partial \bar{\mh}^3$, in such a way that
$$w_0 = ( u_2- u_1)(u_3-u_0)/(u_2-u_0)(u_3- u_1).$$ These 4
points span a `flat' ($2$--dimensional) tetrahedron exactly when the modular
triple is degenerate (real). When it is non-degenerate, we get a positive
embedding of $\Delta$, with its own orientation, onto the corresponding
hyperbolic ideal tetrahedron in $\mh^3$ iff $*_b*_w =1$.

\paragraph{Flattenings and integral charges}\label{fcsec}
Given any $\Ii$--tetrahedron $(\Delta,b,w)$, we consider an
additional decoration made by two $\Z$--valued functions defined on the
edges of $\Delta$, called {\it flattening} and {\it integral charge}
respectively. These functions share the property that {\it opposite edges
take the same value}. Hence it is enough to specify them on the
edges $e_0,e_1,e_2$.

Before we do it, we fix once for ever our favourite standard branch
$\log$ of the logarithm, which has the imaginary part in $]-\pi,\pi]$. We
stress that this $\log$ is defined on $\C\setminus \{0\}$, although it
is not continuous at the negative real half-line.

Let $(\Delta,b,w)$ be an $\Ii$--tetrahedron, and $f=(f_0,f_1,f_2)$, 
with $f_i=f(e_i)\in \Z$. Set
 $${\rm l}_j = {\rm l}_j(b,w,f)=\log(w_j) + i\pi f_j$$
\noindent 
for $j=1$, $2$, $3$. We say that $(f_0,f_1,f_2)$ is a
\emph{flattening} of $(\Delta,b,w)$, and that $(\Delta,b,w,f)$ is a
\emph{flattened} $\Ii$--tetrahedron if
$${\rm l}_0 + {\rm l}_1 +{\rm l}_2 = 0.$$ We call ${\rm
l}_j$ a \emph{log-branch} of $(\Delta,b,w)$ for the edge $e_j$, and
set ${\rm l}=({\rm l}_0,{\rm l}_1,{\rm l}_2)$ for the total log-branch
associated to $f$.

An {\it integral charge} on a branched tetrahedron $(\Delta,b)$
is a function $c=(c_0,c_1,c_2)$, $c_i=c(e_i)\in \Z$, such that
$c_0+c_1+c_2=1$. We call the values of $c$ the {\it charges} of the
edges. A flattened $\Ii$--tetrahedron endowed with an integral charge is
said {\it flat/charged}.

\subsubsection{Matrix dilogarithms}\label{matdefsec}
Here we define the {\it matrix dilogarithms} that are associated to the
flat/charged $\Ii$--tetrahedra. First we describe how any function
 $$A\co  \C\setminus \{0,1\}\to {\rm Aut}(\C^N\otimes \C^N)$$ can be
interpreted as a function of $\Ii$--tetrahedra.  Later we will give the
explicit formulas for the matrix dilogarithms.

We equip $\C^N\otimes \C^N$ with the tensor product of the standard
basis of $\mc^N$, so that $A=A(x)\in {\rm Aut}(\C^N\otimes\C^N)$ is
given by its matrix elements $A^{\delta,\gamma}_{\beta,\alpha}$, where
$\alpha, \ldots, \delta \in \{0,\ldots,N-1\}$. We denote by
$\bar{A}=\bar{A}(x)$ the inverse of $A(x)$, with entries
$\bar{A}^{\beta,\alpha}_{\delta,\gamma}$. 

Take an $\Ii$--tetrahedron $(\Delta,b,w)$. At first, we use the
branching to select one cross-ratio modulus, say $x=w_0$. Then we use again the branching in
order to establish a one--one  correspondence
between the $2$--faces of $\Delta$ and the indices $\alpha, \ldots,
\delta$, and write
\begin{equation}\label{assocdefsec}
A(\Delta,b,w):= A(w_0)^{*_b}.
\end{equation}
The idea (see Section \ref{ftdefsec}) is that when $\Ii$--tetrahedra are glued
along faces, one should be able to form a new tensor by contracting
indices corresponding to paired faces. The slots for indices of the
resulting tensor are in one--one correspondence with the free faces of
the resulting complex. We define the correspondence as follows. Assume that
$*_b=1$. As usual, we name and order the $2$--faces by the opposite
vertices. So, the ordered faces $F_1$, $F_3$ are such that the boundary
and $b$--orientations coincide on them. Set the correspondence
$(F_1,F_3) \leftrightarrows (\alpha,\beta)$. Similarly, set
$(F_0,F_2) \leftrightarrows (\gamma,\delta)$, where $F_0$, $F_2$ are
the ordered faces on which the two orientations do not agree. We do the same
when $*_b=-1$, but in this case the two
orientations agree on $F_0$ and $F_2$. 

\begin{figure}[ht]
\begin{center}
\small
\psfraga <-2pt,-2pt> {a}{$\alpha$} 
\psfraga <-2pt,-2pt> {b}{$\beta$} 
\psfraga <-2pt,-2pt> {c}{$\gamma$} 
\psfraga <-2pt,-2pt> {d}{$\delta$} 
\psfraga <-2pt,-1pt> {x}{$x$} 
\psfraga <-2pt,-2pt> {Ax}{$A(x)$}
\includegraphics[width=7cm]{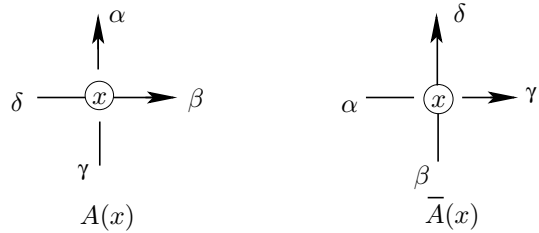}
\caption{\label{CQDtensor} Graphic tensors} 
\end{center}
\end{figure}

It is very convenient to adopt a
pictorial description of this correspondence between automorphisms and
$\Ii$--tetrahedra.  First, the tensors $A(x)$ and $\bar{A}(x)$ may be
given the graphical encoding shown in Figure \ref{CQDtensor}.

\begin{figure}[ht]
\begin{center}
\small
\psfraga <-1pt,0pt> {x}{$x$} 
\psfrag {*1}{$*_b=1$} 
\psfrag {*2}{$*_b=-1$} 
\includegraphics[width=8cm]{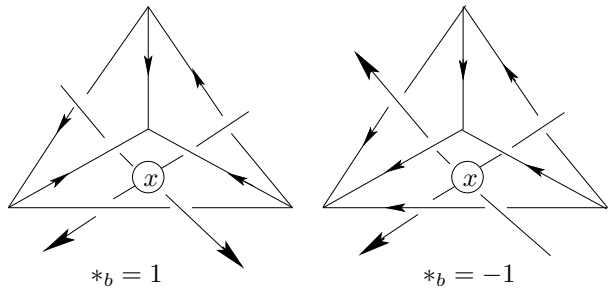}
\caption{\label{CQDidealtensor} $A(x)=A(\Delta,b,w)$, with $x=w_0$} 
\end{center}
\end{figure}

The two figures are normal crossings with an under/over crossing
specification and arc orientations. They are decorated with the
complex parameter $x$ and integers $\alpha, \ldots, \delta$; we have
omitted to draw the arrows on two of the arcs of each crossing,
because we stipulate that they are incoming at the central round
box. Each figure represents a matrix element; forgetting
$\alpha,\beta,\gamma$ and $\delta$, we represent the entire
automorphisms.  We stress that they are {\it planar} pictures,
realized in $\R^2 \cong \C$ with the canonical complex orientation
that is used to specify the index position. Finally we take
$\Ii$--tetrahedra with $w_0=x$ and $*_b = \pm 1$, and we realize the
above graphical encoding of the automorphisms as an {\it enriched}
version of the $1$--skeleton of the canonical cell decomposition of
int$(\Delta)$, which is dual to the natural triangulation of $\Delta$
(so that each arc of the graph is dual to a determined $2$--face of
$\Delta$). This is shown in Figure \ref{CQDidealtensor}. Note that the
embedding in $\Delta$ of this enriched $1$--skeleton is determined by
the branching, and, viceversa, the branching contains all the information in
order to reconstruct completely $(\Delta,b,w)$ itself (this is related
to the encoding of branched spines of $3$--manifolds via so called
normal o--graphs, see \cite{BP2}).

We can give now the explicit formulas for (the matrix
elements of) our matrix dilogarithms $\mathcal{R}_N(\Delta,b,w,f,c)$,
associated to flat/charged $\Ii$--tetrahedra, $N\geq 1$ being any
odd positive integer number.

  For $N=1$, we forget the integral charge $c$, so that
$\mathcal{R}_1$ is defined simply on flattened $\Ii$--tetrahedra. 
Namely, set
\begin{equation}\label{symRmatdil}
\mathcal{R}_1(\Delta,b,w,f) =
\exp\biggl(\frac{*_b}{i\pi}\biggl(-\frac{\pi^2}{6} - \frac{1}{2}
\int_0^{w_0} \biggl( \frac{{\rm l}_0(b,t,f)}{1-t} - \frac{{\rm
    l}_1(b,t,f)}{t} \biggr) \biggr) \biggr)\ dt
\end{equation}
where any ${\rm l}_j(b,t,f)$ is a log-branch 
as defined in Section \ref{fcsec}. This is just the exponential of a
multiple of the lift of the Rogers dilogarithm, discussed in Section \ref{unifRog}.

For $N=2m+1>1$ and every complex number $x$ set
$x^{1/N} = \exp(\log(x)/N)$, where $\log$ is the standard branch of the
logarithm which has the imaginary part in $]-\pi,\pi]$, as already fixed above
(by convention we put $0^{1/N} = 0$). 
Denote by $g$ the complex valued function, analytic
over the complex plane with cuts from the points $x =\zeta^k$ to
infinity ($k=1,\ldots,\ N-1$), defined by
$$g(x) := \prod_{j=1}^{N-1}(1 - x\zeta^{-j})^{j/N}$$
and set $h(x) := g(x)/g(1)$ (we have $g(1) =
\sqrt{N}\exp(-i\pi (N-1)(N-2)/12N)$). The function $g$ plays a main
role in the cyclic representation theory of a Borel quantum subalgebra of
$U_{\zeta}(sl(2,\mc))$ at $\zeta = \exp(2i\pi/N)$ (see the Appendix,
and in particular Theorem \ref{repRzeta-}).

For any $n \in \mn$ and $u'$, $v' \in \mc$ satisfying $(u')^N + (v')^N
= 1$, put
$$\omega(u',v'\vert n) = \prod_{j=1}^n \frac{v'}{1-u'\zeta^j}.$$
The functions $\omega$ are periodic in their integer argument, with
period $N$. Given a flat/charged $\Ii$--tetrahedron $(\Delta,b,w,f,c)$,
set
$$w_j'= \exp((1/N)(\log(w_j) + (f_j-*_bc_j)(N+1)\pi i)).$$
We define
\begin{equation}\label{symqmatdil}
\mathcal{R}_N(\Delta,b,w,f,c) =
\bigl((w_0')^{-c_1}(w_1')^{c_0}\bigr)^{\frac{N-1}{2}}\
(\Ll_N)^{*_b}(w_0',(w_1')^{-1})
\end{equation}
where (recall that $N=2m+1$)
$$\Ll_N(u',v')_{k,l}^{i,j}=h(u')\ \zeta^{kj+(m+1)k^2}\
\omega(u',v'\vert i-k) \ \delta(i + j - l)$$
and $\delta$ is the $N$--periodic Kronecker symbol, ie, $\delta(n) = 1$ if $n
\equiv 0$ mod($N$), and $\delta(n) = 0$ otherwise. Note that for every
$N\geq 1$, the exponent $*_b$ in (\ref{symRmatdil}) and
(\ref{symqmatdil}) is coherent
with that in (\ref{assocdefsec}). The
formula for $\Ll_N^{-1}$ is given in Proposition \ref{unitarity}.

\subsubsection{Transit configurations}\label{transitdefsec}
We define now the {\it transit configurations}, that is the suitable
$\Ii$--flat/charged versions of the $2\to 3$ bistellar (Pachner or
Matveev--Piergallini) local move on $3$--dimensional triangulations,
that will eventually support the fundamental five-term identities
between matrix dilogarithms.

It is useful to fix some general notation for triangulations of
(compact) $3$--dimensional polyhedra. A triangulation, say $T$, can be
considered as a finite family of {\it abstract} tetrahedra with a
fixed identification rule of some pairs of abstract $2$--faces, such
that, after the identification, each $2$--face is common to at most two
tetrahedra of $T$. We also assume that each abstract tetrahedron is
oriented, and that the face identifications reverse the orientation,
so that the resulting polyhedron is also oriented.  Denote by $E(T)$
the set of edges of $T$, by $E_{\Delta}(T)$ the whole set of edges of
the associated abstract tetrahedra, and by $\epsilon_T\co E_{\Delta}(T)
\longrightarrow E(T)$ the natural identification map.

\begin{figure}[ht]
\begin{center}
 \includegraphics[width=8cm]{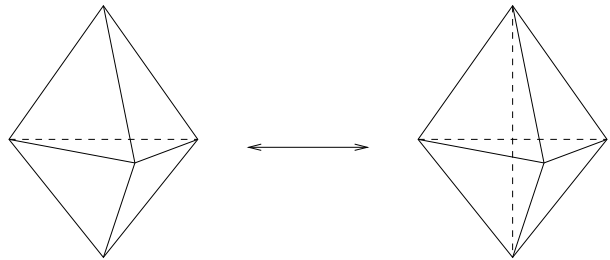}
\caption{\label{CQDfigmove1} The bare $2\to 3$ move between singular triangulations} 
\end{center}
\end{figure}

Consider the $2 \leftrightarrow 3$ move shown in Figure
\ref{CQDfigmove1}. We have two triangulations $T$ and $T'$ (by $2$ and
$3$ tetrahedra respectively) of a same oriented polyhedron, and each
tetrahedron inherits the induced orientation. Assume that each
tetrahedron of $T$ and $T'$ is $\Ii$--flat/charged.  We have to specify
the ``semi-local'' constraints satisfied by each ingredient of the
decoration: branchings, modular triples, flattenings and integral
charges. By ``semi-local'' we mean that the constraints hold between
the decorations of different triangulations of the same topologically
trivial support, where the decorations are given in a purely local way
on each tetrahedron. First of all we require that the local branchings
on $T$ and $T'$ fit well on common edges, and so define globally
branched triangulations $(T,b)$ and $(T',b')$.  

We start by defining the $\Ii$--{\it transits}, then we
will treat the flattening and integral charge transits.  A $2\to 3$
$\Ii$--transit $(T,b,w) \rightarrow (T',b',w')$ consists of a bare
$2\to 3$ move $T \to T'$ that extends to a branching move
$(T,b)\to (T',b')$, ie, the two branchings coincide on the `common'
edges of $T$ and $T'$. Moreover the modular triples have the following
behaviour. For each common edge $e\in E(T)\cap E(T')$ we have
\begin{equation}\label{ideqmodstatesec} \prod_{a\in \epsilon_T^{-1}(e)}w (a)^*=
 \prod_{a'\in \epsilon_{T'}^{-1}(e)}w' (a')^*
\end{equation}
\noindent 
where $*=\pm 1$ according to the $b$--orientation of the abstract
tetrahedron containing $a$ (respectively $a'$).

Note that (\ref{ideqmodstatesec}) implies
that the product of the $w'(a')^*$ around the ``new'' edge of $T'$ is
equal to $1$. So the inverse $3\to 2$ $\Ii$--transits are defined in the very
same way, providing that this last condition is verified on $T'$.

\begin{figure}[ht]
\begin{center}
\small
\psfrag {x}{$x$}
\psfrag {y}{$y$}
\psfrag {yx}{$\frac{y}{x}$}
\psfrag {yxy}{$\frac{y(1-x)}{x(1-y)}$}
\psfrag {1xy}{$\frac{(1-x)}{(1-y)}$}
\includegraphics[width=12cm]{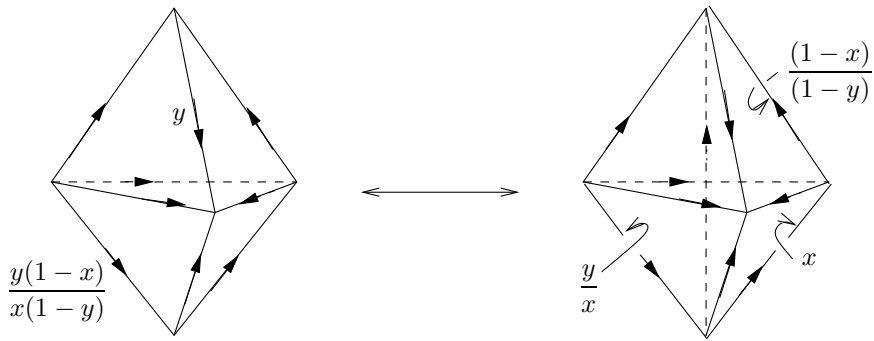}
\caption{\label{CQDidealt} A particular instance of $\Ii$--transit} 
\end{center}
\end{figure}

One particular instance of $\Ii$--transit is shown in Figure
\ref{CQDidealt}. Note that in this case all $*_b$ are equal
to $1$; $x,y$ etc. denotes the cross-ratio modulus $w_0$ of the
corresponding tetrahedron. Assume that all the modular triples
are non degenerate, and share the same sign $*_w=1$. Then we have an
oriented {\it convex} hyperbolic ideal polyhedron with $5$ vertices,
endowed with two different geometric triangulations by two (respectively
three) positively embedded non degenerate ideal tetrahedra. This
situation corresponds to a {\it scissors congruence} relation between
polyhedra in $\mh^3$.
The transit condition (\ref{ideqmodstatesec}), including the exponents
$*_b$, is the natural algebraic extension to situations including
arbitrarily oriented ideal tetrahedra, where the convexity is lost and
there are possible overlappings.

Next we define the notion of a $2\leftrightarrow 3$--transit for
flattened $\Ii$--tetrahedra. Consider a $2 \to 3$ $\Ii$--transit
$(T,b,w) \rightarrow (T',b',w')$ as above. The idea is just to take
formally the $\log$ of the relation (\ref{ideqmodstatesec}). Give a
flattening to each tetrahedron of the initial configuration, and
denote by ${\rm l}\co  E_{\Delta}(T) \rightarrow \mc$ the corresponding
log-branch function on $T$. A map $f'\co  E_{\Delta}(T') \rightarrow \mz$
defines a $2 \to 3$ {\it flattening transit} $(T,b,w,f) \rightarrow
(T',b',w',f')$ if for each common edge $e \in E(T) \cap E(T')$ we have
\begin{equation}\label{ideqfdefsec} 
\sum_{a\in \epsilon_T^{-1}(e)}*\ {\rm l}(a)= \sum_{a'\in
 \epsilon_{T'}^{-1}(e)}*\ {\rm l}'(a')
\end{equation}
\noindent 
where $*=\pm 1$ according to the $b$--orientation of the
tetrahedron that contains $a$ (respectively $a'$).  

It is easily seen that the flattening transits actually define
flattened $\Ii$--tetrah\-edra, and that the sum of the values of ${\rm
  l}'$ about the new edge of $T'$ is always equal to zero. So the
inverse $3 \to 2$ flattening transits are defined in the same way,
except that we also require that this last condition holds. Remark
that the flattenings of a flattening transit associated to a given
$\Ii$--transit $(T,b,w) \rightarrow (T',b',w')$ actually define a
flattening transit for every $\Ii$--transit $(T,b,u) \rightarrow
(T',b',u')$, if $w$ and $w'$ are non degenerate on the abstract
tetrahedra involved in the move and $u$ (respectively $u'$) is a
modular triple sufficiently close to $w$ (respectively $w'$).

It remains to define the transits for the integral charges.
These (like the single charge itself) do not depend on the modular
triples, and even not on the signs $*_b$. A $2\to 3$ branched
move $(T,b,c)\to (T',b',c')$ between charged tetrahedra defines an {\it integral
charge transit} if for each common edge $e \in E(T) \cap E(T')$ we have
\begin{equation}\label{ideqcdefsec} \sum_{a\in \epsilon_T^{-1}(e)}c(a)=
 \sum_{a'\in \epsilon_{T'}^{-1}(e)}c'(a').
\end{equation}
\noindent 
This implies that the sum of the charges around the new edge of $T'$
is always equal to $2$. So we require that this last property is
satisfied when we
define the inverse $3\to 2$ charge transits.  

The $2\to 3$ {\it flat/charged $\Ii$--transits} are defined by
assembling the above definitions.

\subsubsection{Five term relations}\label{ftdefsec}

Here we describe the {\it contraction} of patterns of automorphisms of
$\C^N\otimes \C^N$ associated to patterns of $\Ii$--tetrahedra.

Let $Q$ be any oriented triangulated $3$--dimensional compact
polyhedron. For simplicity, we assume that $Q$ is connected. As
already said, a triangulation $T$ of $Q$ can be considered as a finite
family of abstract tetrahedra $\Delta^i$, with orientation reversing
identifications of some pairs of abstract $2$--faces.  Assume that $T$
is equipped with a {\it global} branching $b$. This means that $b$ is
a system of orientations of the edges of $T$ that restricts to a
branching $b^i$ on each $\Delta^i$ (hence, the face identifications
are compatible with these local branchings). Assume moreover that each
$\Delta^i$ is given a structure of $\Ii$--tetrahedron $(\Delta^i, b^i,
w^i)$, and that we have a function $A\co  \C\setminus \{0,1\}\to {\rm
  Aut}(\C^N\otimes \C^N)$. The correspondence
$A(\Delta^i,b^i,w^i)=A(w^i_0)^{*_{b^i}}$ in (\ref{assocdefsec}) gives us a
pattern of automorphisms of $\C^N\otimes \C^N$.

A {\it state} of $(T,b,w)$ is a function which associate to every
$2$--simplex $t$ of the $2$--skeleton of $T$ an integer $s(t)\in
\{0,\dots,N-1 \}$. So, every state determines a matrix entry for each
$A(\Delta^i,b^i,w^i)$. As two tetrahedra $\Delta^k$, $\Delta^l$ induce
opposite orientations on a common face $t$, the index $s(t)$ is
``down'' for one of $A(\Delta^k,b^k,w^k)$ or $A(\Delta^l,b^l,w^l)$, while it
is ``up'' for the other. By applying the Einstein's rule of ``summing
on repeated indices'', we get the contraction, or {\it trace}, of
this pattern of tensors. We denote this trace by
\begin{equation}\label{ideq3defsec}
\prod_{\Delta \subset T} A(\Delta,b,w).
\end{equation}
The type of the resulting tensor depends on the free $2$--faces, and
their boundary and $b$--orientations. This trace construction can be
very effectively figured out (in the style of spin networks), if we
look at the enriched interior $1$--skeleton of the cell decomposition
dual to the triangulation. For example, in Figure \ref{CQDSchaeffer}
we show the graphical representation (following Figure
\ref{CQDidealtensor}) of the contractions of tensors corresponding to
the two patterns of $\Ii$--tetrahedra involved in Figure
\ref{CQDidealt}.
\begin{figure}[ht]
\begin{center}\small
\psfraga <-2pt, 0pt> {x}{$x$}
\psfraga <-2pt, 0pt> {x1}{$x_1$}
\psfraga <-2pt, 0pt> {x2}{$x_2$}
\psfraga <-2pt, 0pt> {x3}{$x_3$}
\psfraga <-2pt, 0pt> {y}{$y$}
\includegraphics[width=8cm]{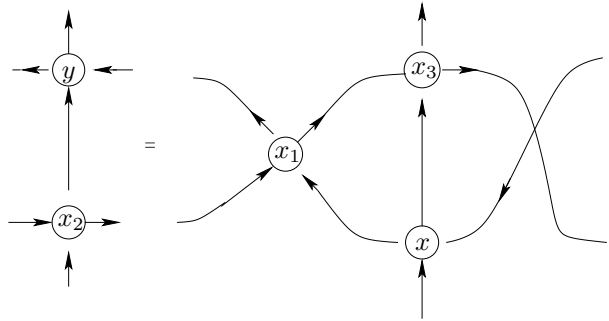}
\caption{\label{CQDSchaeffer} Matrix Schaeffer's identity ($x_1=y/x$, $x_2=y(1-x)/x(1-y)$, and $x_3=(1-x)/(1-y)$)} 
\end{center}
\end{figure}

Now we can state the main results about the semi-local structure
of the matrix dilogarithms. First, remark that if we change the
branching of a flat/charged $\Ii$--tetrahedron $(\Delta,b,w,f,c)$ by a
permutation $p \in S_4$ of its vertices, we get another flat/charged
$\Ii$--tetrahedron $(\Delta,b',w',f',c')= p(\Delta,b,w,f,c)$, where for
each edge $e$ of $\Delta$ we have $w'(e)=w(e)^{\epsilon(p)}$,
$f'(e)=\epsilon(p)f(e)$ and $c'(e)=c(e)$, $\epsilon(p)$ being the
signature of $p$.

There are two main statements, strictly related one to each other. The
first describes the good behaviour of the matrix dilogarithms with
respect to the above action of $S_4$ on flat/charged $\Ii$--tetrahedra
(roughly speaking, it says that the matrix dilogarithms are {\it
  symmetric}); note that when the branching changes, a different
member of the modular triple is selected in the definition of the
matrix dilogarithms. The second statement concerns the fundamental five term
identities. We give an unified statement for all odd $N\geq 1$,
however, it is understood that for $N=1$ we can forget the integral
charges, as they do not enter the definition of the classical
dilogarithm $\Rr_1$. Remark that, as any transit is, in particular, a branching
transit, the traces of the two patterns of associated matrix
dilogarithms are tensors of the same type.
   
\begin{teo}\label{mainlocal} Up to a possible sign or $N$th root of
  unity phase factor, the following properties hold true:

{\rm(1)}\qua {\bf Symmetry}\qua For any permutation $p \in S_4$,
$\mathcal{R}_N(p(\Delta,b,w,f,c))$ is conjugated to
$\mathcal{R}_N(\Delta,b,w,f,c)$, via matrices that depend only on $p$ 
and the branching $b$.

{\rm(2)}\qua {\bf Transit invariance}\qua
For any $2 \rightarrow 3$ flat/charged $\Ii$--transit 
$(T,b,w,f,c)$ $\rightarrow (T',b',w',f',c')$, the traces of the two
patterns of associated matrix dilogarithms lead to the same tensor. In
formula:
$$\prod_{\Delta \subset T} \mathcal{R}_N(\Delta,b,w,f,c) \equiv_N \pm
\prod_{\Delta' \subset T'} \mathcal{R}_N(\Delta',b',w',f',c')$$
where $\equiv_N$ means equality up to multiplication by $N$th roots of
unity.
\end{teo}
We have given here a qualitative formulation of (1); for a
more definite statement, including explicit formulas of the
conjugation matrices, see Corollary \ref{6jsymbis}.

\subsection{Dilogarithmic invariant state sums of globally flat/charged 
$\Ii$--triangulations}\label{GLOBALSTAT}

In this paper, we consider global applications of the matrix
dilogarithms either to compact closed oriented $3$--manifolds $W$
equipped with a flat principal $PSL(2,\C)$--bundle $\rho$, or to cusped
hyperbolic $3$--manifolds $M$, equipped with the
holonomy $\rho$ of the hyperbolic structure (see \cite{BB3} for a
wider range of applications).

\subsubsection{Globally flat charged $\Ii$--triangulations of
  $(W,\rho)$ or $M$} \label{ITdefsec}

The first step is to look for {\it global} 
$\Ii$--{triangulations} of such equipped manifolds.
These are possibly singular
(see Section \ref{GTINV} for more information about
triangulations) globally branched triangulations $(T,b)$ of $W$ or
$M$, such that each tetrahedron $\Delta^i$ of $T$ is equipped with a
modular triple $w^i$, making it an $\Ii$--tetrahedron
$(\Delta^i,b^i,w^i)$. Moreover, we require that at every edge $e$ of
$T$, we have the following {\it edge compatibility condition}:
\begin{equation}\label{ideq1defsec}
\prod_{a \in \epsilon_T^{-1}(e)} w^j(a)^{*_{b^j}} =1
\end{equation}
where $*_{b^j}=\pm 1$ according to the $b^j$--orientation of the
tetrahedron $\Delta^j$ that contains $a$. Note that
(\ref{ideq1defsec}) is the relation satisfied by the cross-ratio
moduli at the edge produced by a $2\to 3$ $\Ii$--transit.  So the edge
compatibility condition is natural to have a class of triangulations
which is {\it stable} for the $2\to 3$ transits. On the other hand, it
is necessary in order to construct hyperbolic $3$--manifolds by gluing
hyperbolic ideal tetrahedra.

In the case of pairs $(W,\rho)$ these $\Ii$--triangulations always
exist, and can be obtained via the {\it idealization} of {\it
$\Dd$--triangulations}, two notions introduced in
\cite{BB2}. We recall them briefly. Let $(\Delta,b,z)$ be a branched
tetrahedron endowed with a $PSL(2,\C)$--valued $1$--cocycle $z$. We
write $z_j = z(e_j)$ and $z'_j = z(e'_j)$. For instance, the
cocycle relation on the $2$--face opposite to $v_3$ reads
$z_0z_1z_2^{-1} = 1$. We say that $(\Delta,b,z)$ is {\it idealizable}
if
$$ u_0=0,\ u_1= z_0(0),\ u_2= z_0z_1(0),\ u_3= z_0z_1z'_0(0) $$ are 4
distinct points in $\C \subset \mc\mathbb{P}^1= \partial
\bar{\mh}^3$. These 4 points span a (possibly flat) hyperbolic ideal
tetrahedron with ordered vertices.

Let $(W,\rho)$ be as above. We consider the pair $(W,\rho)$ up to
orientation preserving homeomorphisms of $W$ and flat bundle
isomorphisms of $\rho$. Equivalently, $\rho$ is identified with a
conjugacy class of representations of the fundamental group of $W$ in
$PSL(2,\C)$.

A $\Dd$--{\it triangulation} of $(W,\rho)$ consists of a
triple $\Tt=(T,b,z)$ where: $T$ is a triangulation of $W$; $b$ is a
global branching of $T$; $z$ is a $PSL(2,\C)$--valued $1$--cocycle on $(T,b)$
representing $\rho$ and such that $(T,b,z)$ is {\it idealizable}, ie,
all its abstract tetrahedra $(\Delta^i,b^i,z^i)$ are idealizable.

If $(\Delta,b,z)$ is idealizable, for all $j=0,1,2$ one can associate
to $e_j$ the cross-ratio modulus $w_j\in
\C\setminus \{0,1\}$ of the hyperbolic ideal tetrahedron spanned by
$(u_0,u_1,u_2,u_3)$.
We call the $\Ii$--tetrahedron $(\Delta,b,w)$ with
$w=(w_0,w_1,w_2)$ the {\it idealization} of $(\Delta,b,z)$. 

For any $\Dd$--triangulation $\Tt=(T,b,z)$ of $(W,\rho)$, its {\it
idealization} $\Tt_{\Ii}=(T,b,w)$ is given by the family $\{(\Delta^i,
b^i, w^i)\}$ of idealizations of the $(\Delta^i,b^i,z^i)$. It is a
fact (see \cite{BB2}, and also Section \ref{GTINV} for more details)
that the idealization of any $\Dd$--triangulation is an
$\Ii$--triangulation, that is it verifies the edge compatibility
condition (\ref{ideq1defsec}). Moreover, every pair $(W,\rho)$ admits
$\Dd$--triangulations.

The situation is more subtle for cusped manifolds $M$.
It is well-known that every such a manifold $M$ admits {\it quasi geometric}
geodesic triangulations by immersed hyperbolic ideal tetrahedra of non-negative
volume (see Remark \ref{rem-tame} (2)). The volume of $M$ is just given by the sum of the volumes of these
tetrahedra. Hence, a quasi geometric
geodesic triangulation of $M$ possibly contains flat tetrahedra
of null volume, but there are strictly positive ones. 
With the usual notation, such a triangulation
gives rise to a pair $(T,w)$, where each modular triple has non-negative
imaginary part, and is only cyclically ordered.

\begin{defi}\label{tame} {\rm A cusped manifold $M$ is said to be {\it
      gentle} if it admits a quasi geometric geodesic triangulation
    $(T,w)$ such that $T$ admits a global branching $b$. In such a
    case, set $w'=w^{*_b}$. Then $(T,b,w')$ is said to be a {\it quasi
      geometric $\Ii$--triangulation} of $M$. For each non degenerate
    $\Ii$--tetrahedron of such an $\Ii$--triangulation, we have
    $*_b*_{w'} = 1$.}
\end{defi}
\begin{remarks}\label{rem-tame}
{\rm (1)\qua To be ``gentle'' is a somewhat demanding
assumption. Nevertheless, many cusped manifolds are gentle. The simplest
example is the complement of the figure-eight knot in the three-sphere. For the
sake of simplicity, in the present section we state the results under
this assumption. However, in Section \ref{GTINV} we show that
the same conclusions hold under the much milder assumption to be
``weakly-gentle'' (see Definition \ref{W-T}). In fact, it is reasonable to ask
whether every cusped manifold is weakly-gentle. If not,
the construction of the {\it quantum} invariants for general cusped
manifolds is more tricky (see Definition \ref{N-W-T}).

(2)\qua We recall a basic procedure to construct quasi geometric geodesic
triangulations of a given cusped manifold $M$. We start with the {\it
  Epstein--Penner canonical cell decomposition} of $M$ \cite {EP}.
This is obtained by identifying pairs of boundary faces of a finite
number of convex ideal hyperbolic polyhedra $\{ G_j \}$, each having a
finite number of faces. Fix a total ordering of the vertices of each
$G_j$ and use it, as usual, to triangulate $G_j$ without adding new
vertices. If the orderings match on the paired faces, we eventually
get a geodesic triangulation of $M$ by strictly positive ideal
tetrahedra, which naturally inherits a global branching from the total
vertex orderings on the $G_j$. If the orderings do not agree on some
pair of identified faces, we have to introduce some degenerate
tetrahedra to get a (quasi geometric) triangulation. This
triangulation does not inherit a global branching from the
construction, but it might nonetheless support some branching.}
\end{remarks}

From now on, in the present section, we consider either pairs
$(W,\rho)$, or gentle cusped manifolds $M$ equipped with
$\Ii$--triangulations as just described.

The notion of {\it global} flattening on an $\Ii$--triangulation of
$(W,\rho)$ or $M$ is obtained by imposing that the associated
log-branches formally satisfy, at each edge $e$ of $T$, the log of the
edge compatibility condition  (\ref{ideq1defsec}). More precisely
\begin{equation}\label{ideq2defsec}
\sum_{a\in \epsilon_{T}^{-1}(e)} *\ {\rm l}(a)=0.
\end{equation}
\noindent Again, this is the natural constraint to get a class of
triangulations which is stable with
respect to the flattening transits. 

Arguing in the same way for the integral charge transits, one would
require that the sum of the charges around every edge of $T$ is equal
to $2$. But a simple `Gauss--Bonnet' argument on each triangulated
sphere making the link of a vertex of a triangulation $T$ of
$(W,\rho)$ shows that such tentative global integral charges {\it do
not exist} (for the triangles of such a link triangulation would
inherit charges $c$ such that the $c\pi$ should behave like the
angles of a {\it flat} triangulation of the $2$--sphere). A way to
overcome this problem is to fix an arbitrary non empty link $L$ in $W$
(considered up to ambient isotopy) and to incorporate this {\it link
fixing} in all the constructions. This eventually leads to the
following notion of $\Dd$--triangulation for a triple $(W,L,\rho)$. A
{\it distinguished triangulation} of $(W,L)$ is a pair $(T,H)$ such
that $T$ is a triangulation of $W$ and $H$ is a {\it Hamiltonian}
subcomplex of the $1$--skeleton of $T$ which realizes the link $L$
(Hamiltonian means that $H$ contains all the vertices of $T$). A
$\Dd$--{\it triangulation} $\Tt=(T,H,b,z)$ for a triple $(W,L,\rho)$
consists of a $\Dd$--triangulation $(T,b,z)$ for $(W,\rho)$ such that
$(T,H)$ is a distinguished triangulation of $(W,L)$. An
$\Ii$--triangulation for $(W,L,\rho)$ is the idealization of a
$\Dd$--triangulation of $(W,L,\rho)$. Finally we can state the notion of global
integral charge:

Let $X$ be either a triple $(W,L,\rho)$ or a gentle cusped
manifold $M$, and $\Tt_{\Ii}$ be an $\Ii$--triangulation of
$X$. A {\it global integral charge} on $\Tt_{\Ii}$ is a
collection of integral charges on the tetrahedra of $\Tt_{\Ii}$ such
that the sum of the charges around every edge of $T$ not belonging to
$H$ is equal to $2$, while the sum of the charges around every edge in
$H$ is equal to $0$ ($H=\emptyset$ when $X=M$).

\subsubsection{Invariant state sums}\label{invdefsec}

Let $(\Tt_{\Ii},f,c)$ be a globally flat/charged
$\Ii$--triangulation of $X$. We can associate to each tetrahedron the
corresponding matrix dilogarithm
$\mathcal{R}_N(\Delta^i,b^i,w^i,f^i,c^i)$, and take the trace as in
(\ref{ideq3defsec}), that we denote $\mathcal{R}_N(\Tt_{\Ii},f,c)$. As
there are no free $2$--faces, we get a scalar. We can give a more
familiar state sum description of this scalar. Recall that a state of
$T$ is function defined on the $2$--simplices of the $2$--skeleton of $T$,
with values in $\{0,\dots,N-1\}$. Any such a state $\alpha$
determines a matrix element
$\mathcal{R}_N(\Delta^i,b^i,w^i,f^i,c^i)_\alpha$ for each matrix
dilogarithm $\mathcal{R}_N(\Delta^i,b^i,w^i,f^i,c^i)$. Set
$$\mathcal{R}_N(\Tt_{\Ii},f,c)_\alpha = 
\prod_i \mathcal{R}_N(\Delta^i,b^i,w^i,f^i,c^i)_\alpha.$$
\noindent Then 
\begin{equation}\label{statesum}
\mathcal{R}_N(\Tt_{\Ii},f,c) = 
\sum_\alpha \mathcal{R}_N(\Tt_{\Ii},f,c)_\alpha.
\end{equation} 
Finally, we can state the main global results about the classical and
quantum dilogarithmic invariants. As for Theorem \ref{mainlocal}, we
give unified statements for all odd $N\geq 1$, but for $N=1$ we can
forget the integral charge and work directly with flattened
$\Ii$--triangulations of $(W,\rho)$ or $M$. On the other hand, in the
quantum case the link $L$ is encoded by the global integral charge,
which is {\it entirely responsible for the link contribution to the
state sums}. Remark also that both global flattenings and integral
charges induce a cohomology class in $H^1(X;\Z/2\Z)$, which is transit
invariant. The invariants depend on the choice of these
classes. Here we prefer to normalize the choice, by requiring that
these classes are trivial. The corresponding flat/charged
$\Ii$--triangulations are said to be (cohomologically) {\it
  normalized}.

\begin{teo}\label{mainglobal} Let  $X$ be either a triple $(W,L,\rho)$
  or a gentle cusped manifold $M$. We have:

{\rm(1)}\qua $X$ admits normalized globally flat/charged $\Ii$--triangulations
    $(\Tt_{\Ii},f,c)$.

{\rm(2)}\qua Let $v$ be the number of vertices of $T$ ($v=0$ for $M$). For every odd integer $N\geq 1$, the value of the
state sum $N^{-v}\mathcal{R}_N(\Tt_{\Ii},f,c)$ does not
depend on the choice of the normalized flat/charged $\Ii$--triangulation of
$X$, possibly up to a sign and multiplication by $N$th roots of
unity. Hence, up to the same ambiguity, it defines a {\rm
  dilogarithmic invariant} $H_N(X)$.
\end{teo}

\noindent A direct consequence of the proof of Theorem \ref{qtransit}
is that for $N \equiv 1$ mod($4$), $N > 1$, the invariants $H_N(X)$
have no sign ambiguity.  The existence of global flattenings and
charges in (1) is based on previous results of Neumann about the
combinatorics of $3$--dimensional triangulations.  For proving (2), we
consider the $\Ii$--decorated versions of few other local moves on
$3$--dimensional triangulations, besides the $2 \leftrightarrow 3$ one,
and the corresponding matrix dilogarithm identities.  Then we show
that arbitrary flat/charged $\Ii$--triangulations of $X$ can be
connected via a finite sequence of such transits together with $2
\leftrightarrow 3$ transits. The reader is addressed to Section
\ref{GTINV} for more information on these invariants.
 
\section {Basic matrix dilogarithms and 
the symmetrization problem}\label{GTBASICMD}
We use the interpretation of automorphisms $A(x)\in \C^N\otimes \C^N$
as functions of $\Ii$--tetrahedra, where $x\in \C \setminus \{0,1\}$,
and the notions of transits and five term identities introduced in Section \ref{LOCSTAT}.

\begin{defi}\label{CQDMD}{\rm A {\it basic matrix dilogarithm of rank
      $N$} is a map $\Ll\co  \C \setminus \{ 0,\ 1 \} \to {\rm
      Aut}(\C^N\otimes \C^N)$ which satisfies the five-term identity
    shown in Figure \ref{CQDSchaeffer}, providing that all the modular
    triples are non degenerate and have imaginary parts of the same
    sign. We call this particular five-term identity with these
    constraints on the cross-ratio moduli the {\it matrix
      Schaeffer's identity}.}
\end{defi}
Recall that Figure \ref{CQDSchaeffer} corresponds to the $\Ii$--transit
of Figure \ref{CQDidealt}, where all the tetrahedra have the same index
$*_b=1$. Note that the Schaeffer's identity holds exactly, with no
phase ambiguity.

\medskip

{\bf The family $\{ \Ll_N \}$}\qua We introduce here the
explicit family $\{ \Ll_N \}$ of {\it basic} matrix dilogarithms of rank $N$
used in this paper. Recall that $N$ is an odd positive integer.

\medskip

\noindent {\bf The classical dilogarithm $\Ll_1$}\qua Definition
\ref{CQDMD} is modeled on the
fundamental functional identity satisfied by the classical Rogers
dilogarithm. As usual, denote by $\log$ the standard branch of the logarithm,
with imaginary part in $]-\pi,\pi]$.  The \emph{Rogers dilogarithm} is the
function over $\C$, complex analytic over $\mathfrak{D} = \mc
\setminus \{(-\infty;0) \cup (1;+\infty) \}$, defined by
\begin{equation}\label{Rdilog}
 {\rm L}(x) = -\frac{\pi^2}{6} -\frac{1}{2} \int_0^x \biggl(
\frac{\log(t)}{1-t} + \frac{\log(1-t)}{t} \biggr) \ dt
\end{equation}
where we integrate first along the path $[0;1/2]$ on 
the real axis and then 
along any path in $\mathfrak{D}$ from $1/2$ to $x$. 
Here we add $-\pi^2/6$ so that ${\rm L}(1)=0$. 
When $\vert x - 1/2 \vert <1/2$ we may also write L as 
$$ {\rm L}(x) = -\frac{\pi^2}{6} + \frac{1}{2} \log(x)\log(1-x) + 
\sum_{n=1}^{\infty} \frac{x^n}{n^2}.$$
The sum in the right-hand side is the power series 
expansion in the open unit disk $\vert x \vert <1$ of the 
\emph{Euler dilogarithm} ${\rm Li}_2$, defined by 
$${\rm Li}_2(x) =  - \int_0^x \frac{\log(1-t)}{t} \ dt$$
and complex analytic over $\mc \setminus (1;+\infty)$. For a
detailed study of the dilogarithm functions and their relatives, see
\cite{L} or the review \cite{Z}.  The function L is related to the
{\it Bloch--Wigner dilogarithm}
\begin{equation} \label{BWvol}
{\rm D}_2(x) = {\rm Im} \bigl( {\rm Li}_2(x) \bigr) 
+ \arg(1-x)\log\vert x\vert
\end{equation}
\noindent which is obtained by adding to ${\rm Im}( {\rm Li}_2(x))$ a
correction term that compensates its jump along the branch cut
$(1;+\infty)$.  The function ${\rm D}_2(x)$ is a real analytic
continuation of ${\rm Im}( {\rm Li}_2(x))$ on $\mc \setminus \{0,1\}$,
and it is continuous (but not differentiable) at $0$ and $1$.  It
gives the volume of $\Ii$--tetrahedra by the formula
$${\rm Vol}(\Delta,b,w) = *_b{\rm D}_2(w_0)$$
\noindent and we have the 6--fold symmetry relations
\begin{equation}\label{BWsym}
{\rm D}_2(w_0) = {\rm D}_2(w_1)= {\rm D}_2(w_2)=
- {\rm D}_2(w_0^{-1}) = -{\rm D}_2(w_1^{-1}) = 
-{\rm D}_2(w_2^{-1}).
\end{equation}
Moreover, if we apply the formula (\ref{BWvol}) to the $\Ii$--transit
of Figure \ref{CQDidealt} we get the five-term functional relation
\begin{equation}\label{BWfivet}
 {\rm D}_2(y) + {\rm D}_2(\frac{1-x^{-1}}{1-y^{-1}}) = {\rm D}_2(x)
 +{\rm D}_2(y/x) + {\rm D}_2(\frac{1-x}{1-y})
\end{equation} 
\noindent when $x \ne y$.  All the other five term relations obtained
by changing the branching in Figure \ref{CQDidealt} also hold true, due
to (\ref{BWsym}).

\noindent One would like to think of the Rogers 
dilogarithm L as the natural complex analytic analogue of ${\rm D}_2(x)$. 
But L verifies similar 
five-term relations only by putting strong restrictions on the variables. 
Namely, the analog of 
(\ref{BWfivet}) is the classical Schaeffer's identity 
\begin{equation}\label{Rfivet}
{\rm L}(x) - {\rm L}(y) +{\rm L}(y/x) - {\rm L}(\frac{1-x^{-1}}{1-y^{-1}}) + 
{\rm L}(\frac{1-x}{1-y})=0
\end{equation}
which for real $x$, $y$ holds only when $0 < y < x < 1$. This identity
characterizes the Rogers dilogarithm: if $f(0;1) \to \mr$ is a 3 times
differentiable function satisfying (\ref{Rfivet}) for all $0 <y <x
<1$, then $f(x) = k{\rm L}(x)$ for a suitable constant $k$ (see eg
\cite{D1}, Appendix). By analytic continuation, the relation
(\ref{Rfivet}) holds true for complex parameters $x$, $y$, providing
that the imaginary part of $y$ is non-zero and $x$ lies inside the
triangle formed by $0$, $1$ and $y$. This is
equivalent to all variables having imaginary parts with the same
sign, as in Definition \ref{CQDMD}. We set
$$ \Ll_1(x) = \exp((1/\pi i){\rm L}(x)).$$
Clearly $\Ll_1$ is a basic matrix dilogarithm of rank $1$. We take the
exponential in order to unify the treatment of the classical and
quantum $(N>1)$ cases.  \medskip

\noindent {\bf The quantum dilogarithms $\Ll_N$}\qua Let
$N=2m+1>1$. Recall the notation introduced in Subsection 
\ref{matdefsec}. We put $u \in \mc \setminus \{0,1\}$, $v = 1-u$, and
define 
\begin{equation}\label{LNBdil}
\Ll_N(u)_{k,l}^{i,j}=\Ll_N(u^{\frac{1}{N}},v^{\frac{1}{N}})_{k,l}^{i,j} =
h(u^{\frac{1}{N}})\ \zeta^{kj+(m+1)k^2}\
\omega(u^{\frac{1}{N}},v^{\frac{1}{N}}\vert i-k) \ \delta(i + j - l).
\end{equation}
\noindent 
Up to a different parametrization, the function $\Ll_N(u',v')$ is the
Faddeev--Kashaev's matrix of $6j$--symbols for the cyclic representation
theory of a Borel quantum subalgebra $\mathcal{B}_{\zeta}$ of
$U_{\zeta}(sl(2,\mc))$, where $\zeta=\exp(2i\pi/N)$ (see Remarks
\ref{compareQHI} and \ref{rempar}). We prove in Section
\ref{CQDSYMQUANTUM} that $\Ll_N(u)$ is actually a basic matrix
dilogarithm of rank $N$, as in Definition \ref{CQDMD}. We can state
now the following problem.

\medskip

{\bf Symmetrization problem for $\Ll_N$}\qua {\sl For every $N$,
find a suitable {\rm symmetrized} version $\Rr_N$ of $\Ll_N$ which satisfies
all the instances of five-term identities, for all transit configurations,
and without any constraint on the modular triples.}
\medskip

It turns out that the solution of this problem is strictly related to
the study of a suitable {\it uniformization} of $\Ll_N$ and to the
behaviour of $\Ll_N$ with respect to the {\it tetrahedral symmetries}.
The flattenings and integral charges arise naturally from this
solution.

\section {The symmetrization problem for $\Ll_1$}\label{CQDSYMROGERS}

\subsection {Uniformization}\label{unifRog}
We use the ``uniformization mod($\pi^2\mz$)'' R of L due to W Neumann
\cite{N1, N2} (see also the recent \cite{N4}).

Let us recall its definition. 
Let $\widehat{\mc} = \widehat{\mc}_{00} \cup \widehat{\mc}_{01} \cup
\widehat{\mc}_{10} \cup \widehat{\mc}_{11}$, where
$\widehat{\mc}_{\varepsilon \varepsilon'}$ ($\varepsilon,
\varepsilon'=0,1$) is the Riemann surface of the function defined on
$\mathfrak{D}=\mc \setminus \{(-\infty;0) \cup (1;+\infty) \}$ by
$$x \mapsto (\log(x) +\varepsilon i\pi ,\log((1-x)^{-1})+
\varepsilon'i\pi ).$$
Thus $\widehat{\mc}$ is the ramified abelian covering of
$\mc \setminus \{0,1\}$ obtained from $\mathfrak{D} \times \mz^2$ by
the identifications
$$\begin{array}{l}
\{(-\infty;0)+i0\} \times \{p\} \times \{q\} \sim \{(-\infty;0) -i0\} 
\times \{p+2\} \times \{q\}\\
\{(1;+\infty)+i0\} \times \{p\} \times \{q\} \sim \{(1;+\infty)-i0\} 
\times \{p\} \times \{q+2\} .
\end{array}$$
Here $(-\infty;0)\pm i0$ comes from the upper/lower fold of $\mathfrak{D}$
with respect to $(-\infty;0)$, and similarly for $(1;+\infty)\pm
i0$. The function 
\begin{equation}\label{logb1}
{\rm l}(x;p,q) = (\log(x) + pi\pi,\log((1-x)^{-1}) +qi\pi)
\end{equation}
is well-defined and analytic on $\widehat{\mc}$. Consider the following lift on
$\widehat{\mc}$ of the Rogers dilogarithm L, defined in (\ref{Rdilog}):
\begin{equation}\label{Ndef}
{\rm R}(x;p,q) = {\rm L}(x) + \frac{i\pi}{2} (p\log(1-x) +q\log(x)).
\end{equation}
It is known that:

\begin{lem}\label{R} The formula (\ref{Ndef}) defines 
an analytic map ${\rm R}\co  \widehat{\mc}\to \C/\pi^2\Z$.
\end{lem}

The idea of interpreting $x$ as a modulus of a hyperbolic ideal
tetrahedron, and $p,\ q$ as additional decorations, comes from
\cite{N1}. We implement this idea in the set up of $\Ii$--tetrahedra formalized in Subsection \ref{LOCSTAT}. Given an
$\Ii$--tetrahedron $(\Delta,b,w)$, let us consider a $\Z$--valued
function $f$ of the edges of $\Delta$ such that, for every edge,
$f(e)=f(e')$. As for $w=(w_0,w_1,w_2)$, we write $f=(f_0,f_1,f_2)$
with the ordering given by the branching $b$. Then we set
$${\rm R}(\Delta,b,w,f) = {\rm R}(w_0;f_0,f_1).$$

\subsection{Tetrahedral symmetries}

Let $(\Delta,b,w,f)$ be as in Section \ref{unifRog}. Here we analyze
under which condition on $f$ the function ${\rm R}(\Delta,b,w,f)$
respects the tetrahedral symmetries. Note that if $f$ is a flattening
of $(\Delta,b,w)$ and $w$ is non degenerate, then $f$ is a flattening
of $(\Delta,b,u)$ for every modular triple $u$ sufficiently close to
$w$.

By acting with a permutation $p \in S_4$ on the vertices of $\Delta$,
one passes from $b$ to a new branching $b'$. This gives
$(\Delta,b',w',f')$, with $w'(e)=w(e)^{\epsilon(p)}$ and
$f'(e)=\epsilon(p)f(e)$ for any edge $e$ of $\Delta$, where
$\epsilon(p)$ is the signature of $p$. Beware that all these data are
renamed according to the new ordering of the vertices given by $b'$.
We have:

\begin{lem}\label{Rsym3} 
For any non-degenerate enriched $\Ii$--tetrahedron $(\Delta,b,w,f)$,
the identities
$${\rm R}(\Delta,b',u',f') = \epsilon(p)\ {\rm R}(\Delta,b,u,f)\quad
{\rm mod}(\pi^2/6) \mz$$ hold true for every permutation $p$ and for
every modular triple $u$ sufficiently close to $w$ if and only if $f$
is a flattening of $(\Delta,b,w)$. These identities are also satisfied
if $w$ is degenerate and we replace $u$ by $w$.
\end{lem}

{\bf Proof}\qua  The basic remark (already made in \cite{N3})
is that the Rogers dilogarithm L has symmetries only up to some
elementary functions. Indeed, by differentiating both sides of each
identity we see that
\begin{eqnarray}\label{Rsym}
{\rm L}\bigl((1-x)^{-1}\bigr) = {\rm L}(x) - \varepsilon (i\pi/2)
\log(1-x) + (\pi^2/6) \nonumber\hspace{0.5cm} \\ {\rm L}(1-x^{-1}) =
{\rm L}(x) - \varepsilon (i\pi/2) \log(x) - (\pi^2/6) \nonumber
\hspace{1.42cm}\\ {\rm L}(x^{-1}) = - {\rm L}(x) + \varepsilon
(i\pi/2) \log(x) \hspace{3.18cm} \\ {\rm L}(1-x) = - {\rm L}(x) -
(\pi^2/6) \nonumber \hspace{4.05cm} \\ {\rm L}\bigl(x/(x-1)\bigr) = -
{\rm L}(x) + \varepsilon (i\pi/2) \log(1-x) - (\pi^2/3)
\hspace{0.21cm} \nonumber
 \end{eqnarray}
when ${\rm Im}(x) \ne 0$, with $\varepsilon = 1$ if ${\rm
Im}(x) > 0$ and $\varepsilon = -1$ if ${\rm Im}(x) < 0$.  A
straightforward computation shows that these relations imply:
\begin{eqnarray}\label{Rsym2}
{\rm R}\bigl((1-x)^{-1};p,q\bigr) = {\rm R}(x;-\varepsilon-p-q,p) 
+ (\pi^2/6) + (p\pi^2/2) \nonumber \hspace{0.02cm} \\
{\rm R}(1-x^{-1};p,q) =  {\rm R}(x;q,-\varepsilon-p-q) -  (\pi^2/6)  
- (q\pi^2/2) \nonumber \hspace{0.36cm} \\
{\rm R}(x^{-1};p,q) = - {\rm R}(x;-p,p+q-\varepsilon) 
- (p\pi^2/2) \hspace{2.1cm} \\ 
{\rm R}(1-x;p,q)  = - {\rm R}(x;-q,-p) -  (\pi^2/6) 
\nonumber \hspace{2.98cm} \\ 
{\rm R}\bigl(x/x-1;p,q\bigr)  =   - {\rm R}(x;p+q-\varepsilon,-q)  -  
(\pi^2/3) + (q\pi^2/2)   \hspace{0.04cm} \nonumber
\end{eqnarray}
\noindent 
under the same assumption. Lemma \ref{R} implies that these relations
are still valid up to $\pi^2$ when $x \in \mathbb{R} \setminus \{ 0,1
\}$. We get the result by renaming the variables according
to the branching. For instance, in the first equality, setting
$(x;-\varepsilon-p-q,p)=(u_0;f_0,f_1)$ we have
$((1-x)^{-1};p,q\bigr)=(u_1;f_1,f_2)$, which is obtained from
$(u_0;f_0,f_1)$ after the
permutation $(012)$.\endproof

\subsection {Complete five term relations}\label{cftrel}
\noindent

Recall the notion of $2\to 3$ flattening transit from Subsection
\ref{ftdefsec}.

\begin{lem}\label{S_Rtransit} 
Let $(T,b,w,f) \rightarrow (T',b',w',f')$ be a $2\to 3$ flattening
transit, such that $(T,b,w) \rightarrow (T',b',w')$ is the
$\Ii$--transit configuration of Figure \ref{CQDidealt}, without any constraint
on the moduli $w$ and $w'$. Then we have
\begin{equation}\label{Rtidsec4}
\sum_{\Delta \subset T} {\rm R}(\Delta,b,w,f) = 
\sum_{\Delta' \subset T'} {\rm R}(\Delta',b',w',f')\quad
 {\rm mod}(\pi^2\mz)\ .
\end{equation} 
 \end{lem} 
 {\bf Proof}\qua This lemma is equivalent to Proposition 2.5 of
 \cite{N2} (see also \cite{N4}). It is based on a clever analytic
 continuation argument that we reproduce for the sake of completeness,
 and because it will be reconsidered in the proof of Theorem
 \ref{qtransit} (quantum case). Denote by $(\Delta^i,b^i,w^i,f^i)$ the flattened
 $\Ii$--tetrahedron opposite to the $i$-th vertex (for the ordering
 induced by $b$). The moduli give us a point
$$(w_0^0,w_0^1,w_0^2,w_0^3,w_0^4)=(x,y,y/x,y(1-x)/x(1-y),(1-x)/(1-y))
\in (\mc \setminus \{0,1\})^5.$$
\noindent Let $\mathfrak{G} \subset (\mc \setminus \{0,1\})^5$ be the
set of such points. Consider the map
$$F\co  \widehat{\mc}^5=\prod_{i=0}^{i=4} \ \{ (w_0^i;f_0^i,f_1^i) \}
\longrightarrow (\mc\setminus \{0,1\})^5 \ $$
\noindent defined by forgetting the $f_j^i$.  Note that the
log-branch functions
$${\rm l}_j^i\co  \widehat{\mc}=\{ (w_0^i;f_0^i,f_1^i) \} \rightarrow \mc $$
\noindent are all analytic, by (\ref{logb1}) and the fact that ${\rm
  l}_2^i=-{\rm l}_0^i-{\rm l}_1^i$ on each flattened
$\Ii$--tetrahedron. Moreover, the relations (\ref{ideqfdefsec}) are
linear identities between the ${\rm l}_j^i$, with $*=1$ for each
summand. Hence they define an \emph{analytic subset}
$\widehat{\mathfrak{G}}$ of $F^{-1}(\mathfrak{G})$.

\begin{figure}[ht]
\begin{center}
{\tiny
\psfrag {0}{$0$}
\psfrag {1}{$1$}
\psfrag {x}{$x$}
\psfrag {y}{$y$}
\psfrag {11x}{$1-1/x$}
\psfrag {11y}{$1-1/y$}
\psfrag {1x}{$1/(1-x)$}
\psfrag {1y}{$1/(1-y)$}
\includegraphics[width=7cm]{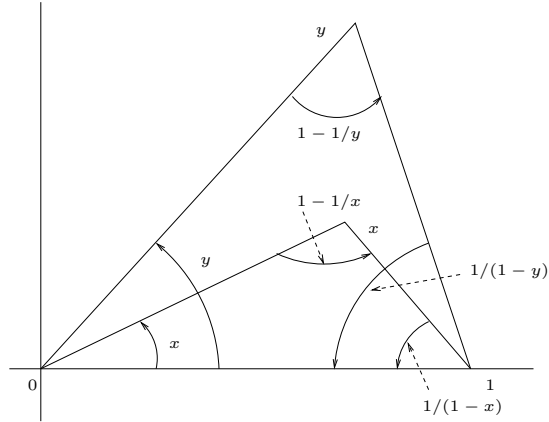}}
\caption{\label{Impositive} Position of $x$ with respect to $y$ in
$\mathfrak{G}^+$, and the associated moduli}
\end{center}
\end{figure}

Denote by $\mathfrak{G}^+ \subset \mathfrak{G}$ the space
where the $w_0^i$ have positive imaginary parts (what
follows could be done with the subset where the $w_0^i$ have
negative imaginary parts). From Figure \ref{Impositive} and the above
description in terms of $x$ and $y$, we see that the points of
$\mathfrak{G}^+$ are characterized by the property that $x$ lies
inside the triangle formed by $0$, $1$ and $y$ with ${\rm Im}(y) > 0$,
so that $\mathfrak{G}^+$ is connected and contractible. Moreover, if
we let ${\rm Im}(x)$ and ${\rm Im}(y)$ go towards $0$ with $0 < {\rm
Re}(y) < {\rm Re}(x) <1$, we come to the subset of $\mathfrak{G}$
where $0 < y < x < 1$ with real $x$ and $y$. We know that
the Schaeffer's identity
$${\rm L}(x) - {\rm L}(y) +{\rm L}(y/x) - {\rm
L}(\frac{1-x^{-1}}{1-y^{-1}}) + {\rm L}(\frac{1-x}{1-y})=0$$
\noindent holds on this subset. Since it is contained in the frontier of
$\mathfrak{G}^+$, and the left-hand side of the Schaeffer's identity
is analytic on $\mathfrak{G}^+$, we deduce by analytic continuation
that the latter holds true on the whole of $\mathfrak{G}^+$.

Next we describe $\widehat{\mathfrak{G}} \cap
F^{-1}(\mathfrak{G}^+)$. In $\mathfrak{G}^+$ the imaginary parts of
the $w_0^i$ are positive, so this is also the case for all the
other moduli of the $\Ii$--transit configuration of
Figure \ref{CQDidealt}. Hence for any edge $e \in E(T) \cap E(T')$ we get
\begin{equation}\label{split}
\sum_{a\in \epsilon_T^{-1}(e)}\log(w(a))= \sum_{a'\in
\epsilon_{T'}^{-1}(e)}\log(w'(a')).
\end{equation}
This implies that the relations (\ref{ideqfdefsec}) are valid over
$F^{-1}(\mathfrak{G}^+)$ if and only if the flattening functions $f\co 
E_{\Delta}(T) \rightarrow \mz$ and $f'\co  E_{\Delta}(T') \rightarrow
\mz$ verify
\begin{equation}\label{sumtransitforf}
\sum_{a\in \epsilon_T^{-1}(e)}f(a)= \sum_{a'\in
\epsilon_{T'}^{-1}(e)}f'(a').
\end{equation}
Let us write (\ref{Rtidsec4}) over $\widehat{\mathfrak{G}} \cap
F^{-1}(\mathfrak{G}^+)$. The dilogarithmic terms of each side are
respectively ${\rm L}(y) + {\rm L}(1-x^{-1}/1-y^{-1})$ and ${\rm L}(x)
+ {\rm L}(y/x) + {\rm L}(1-x/1-y)$, which are equal due to the
Schaeffer's identity. A straightforward computation using
(\ref{split}) shows that the logarithmic terms at each side are equal
if and only if we have :
\begin{equation}\label{sys}\begin{array}{c}
f_1^0 - f_0^2 -  f_1^2 + f_0^3  + f_1^3 = 0 \\ 
-f_1^1 + f_1^2 - f_1^3 = 0 \\
f_0^0 -f_1^3 + f_1^4 = 0 \\
-f_0^1 + f_0^3 + f_1^3 -f_0^4 - f_1^4 = 0 \\
f_0^2 - f_0^3 + f_0^4 = 0.
\end{array}\end{equation}
\noindent Solving this system and using $f_0^i + f_1^i + f_2^i = -1$
(the $w_j^i$ have positive imaginary parts), we find some of the
relations (\ref{sumtransitforf}). Hence the identity (\ref{Rtidsec4})
is true over $\widehat{\mathfrak{G}} \cap F^{-1}(\mathfrak{G}^+)$.
Since $\widehat{\mathfrak{G}}$ is an analytic subset of
$\widehat{\mc}^5$, we deduce from Lemma \ref{R} that (\ref{Rtidsec4})
is also true on the whole of $\widehat{\mathfrak{G}}$ up to $\pi^2$.
\endproof 

We can state now the solution of the symmetrization problem for $\Ll_1$:

\begin{teo}\label{Rtransit} 
Let $(T,b,w,f) \rightarrow (T',b',w',f')$ be any $2\to 3$ flattening
transit. Then we have
\begin{equation}\label{Rtidsec6}
\sum_{\Delta \subset T} *\ {\rm R}(\Delta,b,w,f) = 
\sum_{\Delta' \subset T'} *\ {\rm R}(\Delta',b',w',f')\quad
 {\rm mod}(\pi^2\mz)
\end{equation} 
\noindent where $*=\pm 1$ according to the $b$--orientation of $\Delta$
 (resp.\ $\Delta'$). \end{teo} 

{\bf Proof}\qua By Lemma \ref{S_Rtransit}, the theorem holds
true for the special transit of Figure \ref{CQDidealt}. Any other $2
\leftrightarrow 3$ transit is obtained from this one by changing the
branching. Correspondingly, let us apply Lemma \ref{Rsym3} to
(\ref{Rtidsec4}). We find {\it local} defects, one for each
tetrahedron, which are integer multiples of $\pi^2/6$. We claim that
these defects {\it globally} compensate. As before,
denote by $\Delta^i$ the tetrahedron opposite to the $i$-th vertex in
Figure \ref{CQDidealt}. Any change of branching is obtained as a
composition of the transpositions $(01)$, $(12)$, $(23)$ and $(34)$ of
the vertices. The following table describes for each $\Delta^i$ the
defect induced by these transpositions:
$$
\def\strutt{\vrule width 0pt height 12pt}
\def\struttt{\vrule width 0pt height 8pt}
\begin{tabular}{||c||c|c||c|c|c||} \hline &
  \rule{0cm}{0.4cm}{$\Delta^1$}& \strutt$\Delta^3$ & $\Delta^0$ & $\Delta^2$
  & $\Delta^4$ \\ \hline

\rule{0cm}{0.8cm} \raisebox{0.2cm}{$(01)$} & \raisebox{0.2cm}{$0$} &
\raisebox{0.2cm}{$\frac{f_0^3\pi^2}{2}$} & \raisebox{0.2cm}{$0$} &
\raisebox{0.2cm}{$\frac{f_0^2\pi^2}{2}$} &
\raisebox{0.2cm}{$\frac{f_0^{\struttt4}\pi^2}{2}$}\\ \hline

\rule{0cm}{0.8cm} \raisebox{0.2cm}{$(12)$} & \raisebox{0.2cm}{$0$} &
\raisebox{0.2cm}{$-\frac{\pi^2}{3} -\frac{f_1^3\pi^2}{2}$} &
\raisebox{0.2cm}{$\frac{f_0^0\pi^2}{2}$} & \raisebox{0.2cm}{$0$} &
\raisebox{0.2cm}{$-\frac{\pi^2}{3} -\frac{f_1^{\struttt4}\pi^2}{2}$} \\ \hline

\rule{0cm}{0.8cm} \raisebox{0.2cm}{$(23)$} &
\raisebox{0.2cm}{$-\frac{\pi^2}{3} -\frac{f_1^1\pi^2}{2}$} &
\raisebox{0.2cm}{$0$} & \raisebox{0.2cm}{$-\frac{\pi^2}{3}
-\frac{f_1^0\pi^2}{2}$} & \raisebox{0.2cm}{$0$} &
\raisebox{0.2cm}{$\frac{f_0^{\struttt4}\pi^2}{2}$} \\ \hline

\rule{0cm}{0.8cm} \raisebox{0.2cm}{$(34)$} &
\raisebox{0.2cm}{$\frac{f_0^1\pi^2}{2}$} & \raisebox{0.2cm}{$0$} &
\raisebox{0.2cm}{$\frac{f_0^0\pi^2}{2}$} &
\raisebox{0.2cm}{$\frac{f_0^{\struttt2}\pi^2}{2}$} & \raisebox{0.2cm}{$0$} \\
\hline
\end{tabular}$$
Note that the reduction mod($2$) of the relations
(\ref{sumtransitforf}) are always satisfied over
$\widehat{\mathfrak{G}}$. So this table shows that for any change of
the branching in Figure \ref{CQDidealt} the symmetry defects at both sides
of (\ref{Rtidsec6}) are the same up to $\pi^2$.\endproof 

\begin{remarks}\label{commut} {\rm Dealing with the classical
    ``commutative'' dilogarithm one can prove Theorem \ref{Rtransit}
    without using the tetrahedral symmetries (see \cite{N4}, up to
    some differences in the set up). On the other hand, the path we
    have followed displays the interesting ``local defects vs global
    compensations'' phenomenon. This path is strictly analogous to
    what we shall do in the quantum case. As already remarked in
    \cite{N1}, the proof of Lemma \ref{S_Rtransit} shows that the
    flattening transits realize the most general relations between
    enriched $\Ii$--tetrahedra for which the identities
    (\ref{Rtidsec6}) are universaly true, that is independently of the
    specific values of $w$ and $w'$.}
\end{remarks}

\noindent Finally, as in (\ref{symRmatdil}) we set
$$\Rr_1(\Delta,b,w,f)=\exp((*_b/i\pi){\rm R}(\Delta,b,w,f))$$
The map $\Rr_1$ gives us the symmetrized matrix dilogarithm of rank $1$.
Clearly it satisfies the conclusion of Theorem \ref{mainlocal}.
 
\section {The symmetrization problem for $\Ll_N$, $N > 1$}\label{CQDSYMQUANTUM}
The Appendix collects some quantum algebraic facts used in the present
section. In Section \ref{LNbdilog}, we describe the lifted matrix
Schaeffer's identity for the matrix $\widehat{\Ll}_N$ obtained in
Subsection \ref{repR} of the Appendix. In Section \ref{symqrel} we
compute the tetrahedral symmetries of $\widehat{\Ll}_N$ and we
prove Theorem \ref{mainlocal} for $N>1$.

As before, let $N=2m+1>1$ be any odd positive integer, and denote by
log the standard branch of the logarithm, which has the imaginary part
in $]-\pi,\pi]$. For any complex number $x \ne 0$ write
$x^{1/N}=\exp((1/N)\log(x))$. We denote $\zeta=\exp(2i\pi/N)$. Remark
that $\zeta^{m+1} = -\exp(i\pi/N)$, so that $\zeta^{N(m+1)} = 1$. For any $u \in \mc \setminus \{0,1\}$ and $p \in
\mz$ define
\begin{equation}\label{rootN}
u_p'= u_0'\ \zeta^{(m+1)p} = \exp((1/N)(\log(u) + p(N+1)\pi i)).
\end{equation}
We can lift $\Ll_N$, given in (\ref{LNBdil}), over the
Riemann surface $\widehat{\mc}$ of Section \ref{CQDSYMROGERS} by
setting
\begin{equation}\label{newL}
\widehat{\Ll}_N(u;p,q)_{i,j}^{k,l} = \Ll_N(u_p',v_{-q}')_{i,j}^{k,l}
= \frac{g(u_p')}{g(1)}\ \zeta^{il+(m+1)i^2}\
\omega(u_p',v_{-q}'\vert k-i) \ \delta(k + l - j)
\end{equation}
for any $(u;p,q) \in \{ \mc \setminus \{(-\infty,0) \cup
(1,+\infty)\}\} \times \mz^2$, where $v=1-u$. The matrix
$\widehat{\Ll}_N(u;p,q)$ is invertible, with inverse given in
Proposition \ref{unitarity}. Remark that if $u \in
(-\infty,0)$ then $\widehat{\Ll}_N(u+i0;p,q) =
\widehat{\Ll}_N(u-i0;p+2,q)$ because
$$(1/N)(\log(u+i0) + p(N+1)\pi i) = (1/N)(\log(u-i0) + (p+2)(N+1)\pi
i) -2\pi i.$$ On another hand, if $u \in (1,+\infty)$ then $u_p'$ lies
on the ray $\{t\zeta^{(m+1)p}\}$, $t > 1$. As this is a branch cut of
the function $g$ in (\ref{newL}), we have 
$$\widehat{\Ll}_N(u+i0;p,q) =
\zeta^{-(m+1)p}\widehat{\Ll}_N(u-i0;p,q+2).$$ Hence, denoting by $U_N$
the multiplicative group of $N$th roots of unity, we see that the
matrix valued map $\hat{\Ll}_N\co  \widehat{\mc} \rightarrow {\rm
  M}_{N^2}(\mc/U_N)$ is complex analytic (compare with Lemma \ref{R}).
Recall that $\equiv_N$ denotes the equality up to multiplication by
$N$th roots of unity.

\subsection{Lifted basic five term relation}\label{LNbdilog}

We say that $(\Delta,b,w,a)$ is an {\it enriched} $\Ii$--tetrahedron if
$(\Delta,b,w)$ is an $\Ii$--tetrahedron and $a$ is a $\mz$--valued
function on the edges of $\Delta$ such that $a(e)=a(e')$ for every
pair of opposite edges $e$ and $e'$. We identify $a$ with
$(a_0,a_1,a_2)$, where $a_j=a(e_j)$ and the ordering of the edges is
induced by the branching $b$. Similarly to (\ref{rootN}), given $a$ we
define $N$th roots of the moduli by
\begin{equation}\label{Nroot}
w_j'= w_{a_j}'=\exp((1/N)(\log(w_j) + a_j(N+1)\pi i)).
\end{equation}
We call $w'\co  E(\Delta) \rightarrow \mc \setminus \{0,1\}$
the {\it $N$th--branch} of $w$ for $a$ (for short: $N$th--branch map),
and its values are the {\it $N$th--root moduli}. We write
\begin{equation}\label{prodconst}
\tau = -w_0'w_1'w_2'.
\end{equation}
\begin{defi}\label{defenricht}{\rm Consider a $2 \to 3$ $\Ii$--transit
    $(T^0,b^0,w^0) \rightarrow (T^1,b^1,w^1)$ whose underlying
    branching transit is as in Figure \ref{CQDidealt}. Suppose that we
    have a map $a^0$ that enriches the tetrahedra of $T^0$ involved in
    the move. A map $a^1$ that enriches those in $T^1$
    defines a {\it $N$th--branch transit} if for each common edge $e
    \in E(T^0) \cap E(T^1)$ we have}
\begin{equation}\label{ideqmodprime} 
\prod_{\tilde{e}^0\in \epsilon_{T^0}^{-1}(e)}(w^0)'(\tilde{e}^0)=\prod_{\tilde{e}^1\in
\epsilon_{T^1}^{-1}(e)}(w^1)'(\tilde{e}^1)
\end{equation}
\noindent {\rm where the identification map
  $\epsilon_{T^i}\co E_{\Delta}(T^i) \rightarrow E(T^i)$ is as in Subsection \ref{transitdefsec}.}
\end{defi}
It is easily seen that (\ref{ideqmodprime}) implies that the
$N$th--roots of unity $\tau$ in (\ref{prodconst})
are the same for all tetrahedra. Also, the
product of the $N$th--root moduli about the new edge of $T^1$ is equal to
$\tau^2$.

For any enriched $\Ii$--tetrahedron
$(\Delta,b,w,a)$ we define
\begin{equation}\label{dilogtet}
\widehat{\Ll}_N(\Delta,b,w,a) = 
(\widehat{\Ll}_N)^{*_b}(w_0;a_0,a_1) = (\Ll_N)^{*_b}(w_0',(w_1')^{-1}).
\end{equation}

The notion of five term identity (in particular the matrix
Schaeffer's one) naturally lifts to enriched $\Ii$--tetrahedra and
$N$th--branch transits. In the rest of this section we prove:

\begin{teo} \label{basNtr} The matrix Schaeffer's identity
corresponding to any $N$th--branch $\Ii$--transit holds true for the
tensors $\widehat{\Ll}_N(\Delta,b,w,a)$, with furthermore no
restriction on the cross-ratio moduli.
\end{teo}

Two remarks are in order. First, when the moduli of the tetrahedra
involved in the move satisfy the conditions of Definition \ref{CQDMD},
the $N$th--root moduli given by the standard log (ie, with $a\equiv
0$) make a $2 \leftrightarrow 3$ $N$th--branch transit. So this theorem
implies that the matrix $\Ll_N$, as defined in (\ref{LNBdil}), {\it is
  a basic matrix dilogarithm of rank $N$}. Second, here we impose a
specific branching transit because we have not yet analyzed the
symmetries of $\widehat{\Ll}_N(\Delta,b,w,a)$. We shall relax this
assumption in Section \ref{symqrel}.

\medskip

\begin{figure}[ht]\small
\psfrag{a}{$(w_0^3,a_0^3)$}
\psfrag{b}{$(w_0^1,a_0^1)$}
\psfrag{c}{$(w_0^2,a_0^2)$}
\psfrag{d}{$(w_0^0,a_0^0)$}
\psfrag{e}{$(w_0^4,a_0^4)$}
\begin{center}
 \includegraphics[width=10cm]{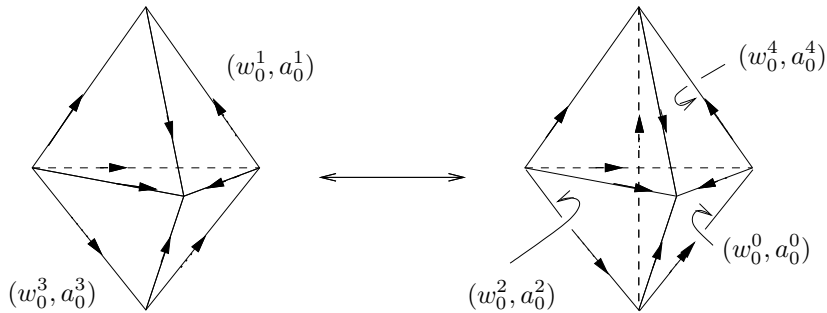}
 \caption{\label{CQDidealtnew} The enriched $\Ii$--transit supporting
   the matrix Schaeffer's identity}
\end{center}
\end{figure}

{\bf Proof of Theorem \ref{basNtr}}\qua Denote by $(\Delta^i,b^i,w^i,a^i)$
the enriched $\Ii$--tetrahedron opposite to the $i$-th vertex in
Figure \ref{CQDidealtnew}. Using Figure \ref{CQDidealtensor} we see that the
associated (Schaeffer's) five-term identity reads
$$\begin{array}{l} \widehat{\Ll}_N(\Delta^1,b^1,w^1,a^1)_{23}\ \
\widehat{\Ll}_N(\Delta^3,b^3,w^3,a^3)_{12} = \hspace{4cm}\\ \\
\hspace{1cm} \widehat{\Ll}_N(\Delta^4,b^4,w^4,a^4)_{12}\ \
\widehat{\Ll}_N(\Delta^2,b^2,w^2,a^2)_{13}\ \
\widehat{\Ll}_N(\Delta^0,b^0,w^0,a^0)_{23}.
\end{array}$$
Both sides are operators acting on $\mc^N \otimes \mc^N \otimes
\mc^N$. The indices show the tensor factor on which the
$\widehat{\Ll}_N$ act, for instance $Y_1^{-1}Z_2^{-1}Y_2 = Y^{-1}
\otimes Z^{-1}Y \otimes {\rm id}_{\mc^N}$, and so on. 

The splitting
formula induced by Theorem \ref{repRzeta-} gives
\begin{equation}\label{formmat2a}
\Ll_N(w_0',(w_1')^{-1}) = \Upsilon \ \cdot \ \Psi(-Y^{-1} \otimes Z^{-1}Y)
\end{equation}
where
$$\Upsilon = \frac{1}{N}\ \sum_{i,j =0}^{N-1} \zeta^{ij}\ Z^{-i} \otimes Y^j$$
and
\begin{equation}\label{formmat2b}
\Psi(-Y^{-1} \otimes Z^{-1}Y)= \frac{g(w_0')}{g(1)}\ \sum_{t=0}^{N-1}
\prod_{s=1}^{t} \frac{(w_0')^{-1}(w_1')^{-1}}{1-(w_0')^{-1}\zeta^{-s}}
\ (-Y^{-1} \otimes Z^{-1}Y)^t\hspace{0.5cm}
\end{equation}
is obtained by reversing the computation after formula
(\ref{formmat}). A remarkable fact is that $\Psi$ is a solution (in
fact, the unique up to multiplication by scalars) of the functional
relation
\begin{equation}\label{factor2}
\Psi(\zeta^{-1}A) =\Psi(A)\ \left( \frac{1-(w_0')^{-1}
(w_1')^{-1}A}{(w_0')^{-1}}\ \right) = \Psi(A)\ \left(w_0' -
(w_1')^{-1}A \right)
\end{equation}  
where $A=-Y^{-1} \otimes Z^{-1}Y$. By (\ref{formmat2a}) we have to prove:
\begin{equation} \label{splitform}
\Upsilon_{23}\Psi_{23}^1\ \Upsilon_{12}\Psi_{12}^3 =
\Upsilon_{12}\Psi_{12}^4\ \Upsilon_{13}\Psi_{13}^2\
\Upsilon_{23}\Psi_{23}^0
\end{equation}
\noindent where the $\Psi^i$ are given by (\ref{formmat2b}) for each
enriched tetrahedron $(\Delta^i,b^i,w^i,a^i)$, and we omit their
matrix arguments for simplicity. The first step is to split this
relation into the pentagon relation (\ref{pentagone}) for $\Upsilon$,
and a five term identity for $\Psi$ that we shall consequently
prove. Write
$$U = -Y_1^{-1}Z_2^{-1}Y_2 = -\zeta^{-1} (XZ)_1^{-1}X_2\quad ,\quad V
= -Y_2^{-1}Z_3^{-1}Y_3$$
\noindent where the matrices $X$, $Y$ and $Z$ are defined in the proof
of Theorem \ref{repRzeta-}. By commuting the variables we easily verify
that
$$\begin{array}{l} \Psi_{13}^2(-\zeta^{-1} (XZ)_1^{-1}X_3)\
\Upsilon_{23} = \Upsilon_{23}\ \Psi_{13}^2(-\zeta^{-1}
(XZ)_1^{-1}Z_2^{-1}X_3) \hspace{4cm}\\ \\
\hspace{4.15cm} = \Upsilon_{23}\ \Psi_{13}^2\left(-(-\zeta^{-1}
(XZ)_1^{-1}X_2)(-\zeta^{-1} (XZ)_2^{-1}X_3)\right) \\ \\
\hspace {4.15cm} = \Upsilon_{23}\ \Psi_{13}^2(-UV)
\end{array}$$
$$\begin{array}{lll} \Psi_{12}^4(-\zeta^{-1} (XZ)_1^{-1}X_2)\
\Upsilon_{13} & = & \Upsilon_{13}\ \Psi_{12}^4(-\zeta^{-1}
(XZ)_1^{-1}X_2(XZ)_3)\hspace {3cm}\\ \\ & = & \Upsilon_{13}\
\Psi_{12}^4(U(XZ)_3)\\ \\ \Psi_{12}^4(-\zeta^{-1} (XZ)_1^{-1}X_2)\
\Upsilon_{23} & = & \Upsilon_{23}\ \Psi_{12}^4(-\zeta^{-1}
(XZ)_1^{-1}X_2(XZ)_3^{-1}) \\ \\ & = & \Upsilon_{23}\
\Psi_{12}^4(U(XZ)_3^{-1}).
\end{array}$$
\noindent Then the right-hand side of (\ref{splitform}) is equal to
(for simplicity, we only indicate some of the matrix arguments)
$$\begin{array}{l} \Upsilon_{12}\Psi_{12}^4\ \Upsilon_{13}\Psi_{13}^2\
\Upsilon_{23}\Psi_{23}^0(V)= \Upsilon_{12}\Psi_{12}^4\
\Upsilon_{13}\Upsilon_{23}\ \Psi_{13}^2(-UV)\ \Psi_{23}^0(V)
\hspace{3cm}\\ \\
\hspace{4.35cm} = \Upsilon_{12}\Upsilon_{13}\
\Psi_{12}^4(U(XZ)_3)\Upsilon_{23}\ \Psi_{13}^2(-UV)\ \Psi_{23}^0(V) \\
\\
\hspace{4.35cm} =\Upsilon_{12}\Upsilon_{13}\Upsilon_{23}\
\Psi_{12}^4(U)\ \Psi_{13}^2(-UV)\ \Psi_{23}^0(V).
\end{array}$$
\noindent and the left-hand side immediately gives
$$\Upsilon_{23}\Psi_{23}^1\ \Upsilon_{12}\Psi_{12}^3 =
\Upsilon_{23}\Upsilon_{12}\ \Psi_{23}^1\Psi_{12}^3.$$
\noindent As $\Upsilon$ is a linear representation of the canonical
element $S_\zeta$ of the algebra $Q(\mathcal{B}_{\zeta}^0)$ (see
Section \ref{algprel}), it is a solution of the pentagon relation
(\ref{pentagone}). So we are left to show that
\begin{equation} \label{FadKas}
\Psi_{23}^1(V)\Psi_{12}^3(U) =
\Psi_{12}^4(U)\Psi_{13}^2(-UV)\Psi_{23}^0(V)\quad .
\end{equation}
\noindent We first prove that $\Psi_{12}^4(U)^{-1}\ \Psi_{23}^1(V)\
\Psi_{12}^3(U)\ \Psi_{23}^0(V)^{-1}$ commutes with $UV$. For that it
is enough to observe that $UV = \zeta VU$ and to use (\ref{factor2}). Namely, 
$$\begin{array}{l} \Psi_{12}^4(U)^{-1}\ \Psi_{23}^1(V)\
 \Psi_{12}^3(U)\ \Psi_{23}^0(V)^{-1}\ (UV) \hspace{6cm}\\ \\
\hspace{2cm}= \Psi_{12}^4(U)^{-1}\ \Psi_{23}^1(V)\ 
\Psi_{12}^3(U)\ (UV)\ \Psi_{23}^0(\zeta^{-1}V)^{-1} \\ \\
\hspace{2cm}= \Psi_{12}^4(U)^{-1}\ \Psi_{23}^1(V)\ (UV) \
\Psi_{12}^3(\zeta U)\ \Psi_{23}^0(\zeta^{-1}V)^{-1}\\ \\
\hspace{2cm}= \Psi_{12}^4(U)^{-1}\ (UV)\ \Psi_{23}^1(\zeta^{-1}V)\
\Psi_{12}^3(\zeta U)\ \Psi_{23}^0(\zeta^{-1}V)^{-1}\\ \\
\hspace{2cm}= (UV)\ \Psi_{12}^4(\zeta U)^{-1} \
\Psi_{23}^1(\zeta^{-1}V)\ \Psi_{12}^3(\zeta U)\
\Psi_{23}^0(\zeta^{-1}V)^{-1}.
\end{array}$$
\noindent Moreover, (\ref{factor2}) allows us to turn the last four
terms into
$$\begin{array}{l} \Psi_{12}^4(U)^{-1}\ \left((w_0')^4 - \zeta
((w_1')^4)^{-1}\ U\right) \ \Psi_{23}^1(\zeta^{-1}V)\
\Psi_{12}^3(\zeta U)\ \Psi_{23}^0(\zeta^{-1}V)^{-1} \hspace{4cm}\\ \\
= \Psi_{12}^4(U)^{-1}\ \left( \Psi_{23}^1(\zeta^{-1}V)\ (w_0')^4 -
\zeta ((w_1')^4)^{-1}\ \Psi_{23}^1(V)\ U \right) \\
 \hspace{8cm}\times \Psi_{12}^3(\zeta U)\
 \Psi_{23}^0(\zeta^{-1}V)^{-1}\\ \\ = \Psi_{12}^4(U)^{-1}\
 \Psi_{23}^1(V)\ \left( ((w_0')^1 - ((w_1')^1)^{-1}\ V) (w_0')^4 -
 \zeta ((w_1')^4)^{-1}\ U \right) \\
 \hspace{8cm}\times \Psi_{12}^3(\zeta U)\ \Psi_{23}^0(\zeta^{-1}V)^{-1}
\end{array}$$
$$\begin{array}{l} = \Psi_{12}^4(U)^{-1}\ \Psi_{23}^1(V)\ \bigl (
 \Psi_{12}^3(\zeta U)\ \bigl((w_0')^1(w_0')^4 - \zeta ((w_1')^4)^{-1}\
 U \bigr) - \\
 \hspace{5cm} ((w_1')^1)^{-1}(w_0')^4\ \Psi_{12}^3(U)\ V \bigr)\
\Psi_{23}^0(\zeta^{-1}V)^{-1} \\ \\ = \Psi_{12}^4(U)^{-1}\
\Psi_{23}^1(V)\ \Psi_{12}^3(\zeta U)\ \bigl((w_0')^1(w_0')^4 - \zeta
((w_1')^4)^{-1}\ U - \\
 \hspace{3cm} ((w_1')^1)^{-1}(w_0')^4\bigl((w_0')^3 - \zeta ((w_1')^3)^{-1}\ 
U\bigr)\ V \bigr)\ \Psi_{23}^0(\zeta^{-1}V)^{-1}.
\end{array}$$
Since there is an $N$th-branch transit, we have (see Figure \ref{CQDidealtnew}):\begin{equation}\label{relmod}
(w_1')^0(w_0')^4 = (w_1')^1 \quad ,\quad (w_2')^3(w_0')^1 = (w_2')^4 \quad , \quad (w_0')^0(w_1')^4 = (w_1')^3.
\end{equation}
\noindent Together with the relations (\ref{prodconst}), the last two imply $(w_0')^1(w_0')^4 = (w_0')^3(w_0')^0$. Then, in the last expression, the term between parenthesis reads
$$\bigl((w_0')^3 - \zeta ((w_1')^3)^{-1}\ U\bigr)\ \bigl((w_0')^0 - ((w_1')^0)^{-1}\ V\bigr).$$
\noindent By applying (\ref{factor2}) two more times we eventually find
$$\begin{array}{l}
\Psi_{12}^4(U)^{-1}\ \Psi_{23}^1(V)\ \Psi_{12}^3(U)\ \Psi_{23}^0(V)^{-1}\ (UV)= \hspace{6cm}\\
\hspace{5cm}(UV)\ \Psi_{12}^4(U)^{-1}\ \Psi_{23}^1(V)\ \Psi_{12}^3(U)\ \Psi_{23}^0(V)^{-1}.
\end{array}$$
\noindent This shows that $P(-UV) = \Psi_{12}^4(U)^{-1}\ \Psi_{23}^1(V)\ \Psi_{12}^3(U)\ \Psi_{23}^0(V)^{-1}$ is a linear functional of $-UV$ (there must be a $-$ sign in front of $UV$, because $U^N=V^N=-(UV)^N=-{\rm Id}_{\mc^N}$). To conclude the proof, it is enough to show that $P(-UV)$ satisfies
\begin{equation} \label{eqnouv}
P(-\zeta^{-1}UV) = P(-UV) \ \bigl((w_0')^2 - ((w_1')^2)^{-1}\  (-UV)\bigr).
\end{equation}
\noindent Indeed, this equation defines $P(-UV)$ as well as $
\Psi_{13}^2(-UV)$ up to a scalar, and by Lemma \ref{unitarity} we know that for each $i$ the matrix $\widehat{\Ll}_N(\Delta^i,b^i,w^i,a^i)$ has determinant $1$. Consider the change of variable $V \rightarrow \zeta^{-1}V$ in $P(-UV)$. We have
$$\begin{array}{l}
\Psi_{12}^4(U)^{-1}\ \Psi_{23}^1( \zeta^{-1}V)\ \Psi_{12}^3(U)\ \Psi_{23}^0(\zeta^{-1}V)^{-1}  \\ \\
 \hspace{0.3cm} =  \Psi_{12}^4(U)^{-1}\ \Psi_{23}^1(V)\ \bigl((w_0')^1 - ((w_1')^1)^{-1}\ V\bigr)\ \Psi_{12}^3(U)\ \Psi_{23}^0(\zeta^{-1}V)^{-1}\\ \\
 \hspace{0.3cm} =  \Psi_{12}^4(U)^{-1}\ \Psi_{23}^1(V)\ \bigl( (w_0')^1\ \Psi_{12}^3(U) - ((w_1')^1)^{-1} \Psi_{12}^3(\zeta^{-1}U)\ V \bigr)\ \Psi_{23}^0(\zeta^{-1}V)^{-1}\\ \\
\hspace{0.3cm} = \Psi_{12}^4(U)^{-1}\ \Psi_{23}^1(V)\ \Psi_{12}^3(U)\ \bigl(  (w_0')^1 - ((w_1')^1)^{-1}\bigl((w_0')^3 - ((w_1')^3)^{-1}\ U\bigr)V \bigr)\\
\hspace{10cm}\times \ \Psi_{23}^0(\zeta^{-1}V)^{-1} \\ \\
\hspace{0.3cm}= \Psi_{12}^4(U)^{-1}\ \Psi_{23}^1(V)\ \Psi_{12}^3(U)\ \bigl( \Psi_{23}^0(\zeta^{-1}V)^{-1}\bigl((w_0')^1- ((w_1')^1)^{-1}(w_0')^3\bigr) \\
\hspace{6cm} + ((w_1')^1)^{-1} ((w_1')^3)^{-1}\ \Psi_{23}^0(V)^{-1}\ UV \bigr).
\end{array}$$
\noindent Again, since we have an $N$th-branch transit the following relations hold true:
$$(w_0')^3=(w_0')^2(w_0')^4\quad ,\quad (w_0')^1=(w_0')^0(w_0')^2\quad , \quad (w_1')^2=(w_1')^1(w_1')^3. $$
\noindent Together with the first relation in (\ref{relmod}), the first above gives $((w_1')^1)^{-1}(w_0')^3=(w_0')^2((w_1')^0)^{-1}$. So the term between parenthesis is equal to
$$\Psi_{23}^0(\zeta^{-1}V)^{-1}\bigl((w_0')^0- ((w_1')^0)^{-1}\ V\bigr)(w_0')^2 + ((w_1')^2)^{-1}\ \Psi_{23}^0(V)^{-1}\ UV$$
and we find
$$\begin{array}{l}
\Psi_{12}^4(U)^{-1}\ \Psi_{23}^1( \zeta^{-1}V)\ \Psi_{12}^3(U)\ \Psi_{23}^0(\zeta^{-1}V)^{-1}\hspace{4cm} \\ \\
\hspace{1.5cm}  = \Psi_{12}^4(U)^{-1}\ \Psi_{23}^1(V)\ \Psi_{12}^3(U)\ \Psi_{23}^0(V)^{-1} \left((w_0')^2 - ((w_1')^2)^{-1} (-UV)\right).
\end{array}$$
\noindent This proves (\ref{eqnouv}), whence the theorem. \endproof

\subsection{Tetrahedral symmetries}\label{symqrel}

As in Section \ref{LNbdilog}, we define an enriched $\Ii$--tetrahedron
$(\Delta,b,w,a)$ from a flat/charged one $(\Delta,b,w,f,c)$ by putting
$a=f-*_bc$. In that case, (\ref{Nroot}) gives $\tau=\zeta^{-*_b(m+1)}$ in
(\ref{prodconst}).

Recall that the symmetry group on four elements numbered from $0$ to
$3$ is generated by the transpositions $(01)$, $(12)$ and $(23)$. We
saw in Section \ref{CQDSYMROGERS} that if we change the branching of a
flattened $\Ii$--tetrahedron $(\Delta,b,w,f)$ by a permutation $(ij)$ of its
vertices, we get another flattened $\Ii$--tetrahedron $(ij)(\Delta,b,w,f)$. This
behaviour naturally extends to flat/charged $\Ii$--tetrahedra.

Let $S$ and $T$ be the $N \times N$ invertible square matrices with
entries
$$T_{i,j} = \zeta^{i^2(m+1)}\delta(i+ j)\quad ,\quad S_{i,j} =
N^{-1/2}\zeta^{ij}.$$ We have $S^4 = {\rm id}_{\mc^N}$ and $S^2 =
\zeta'(ST)^3$ for some $N$th--root of unity $\zeta'$. Hence the
matrices $S$ and $T$ define a projective $N$--dimensional
representation $\Theta$ of $SL(2,\mz)$. The following proposition
describes the tetrahedral symmetries of
$\widehat{\Ll}_N(\Delta,b,w,f,c)$ in terms of $\Theta$.
\begin{prop} \label{6jsym} 
Let $(\Delta,b,w,f,c)$ be a flat/charged $\Ii$--tetrahedron with
$*_b=+1$. We have
$$\begin{array}{l} \widehat{\Ll}_N\bigl((01)(\Delta,b,w,f,c)\bigr)
\equiv_N (w_0')^{\frac{1-N}{2}}\ T_1^{-1}\
\widehat{\Ll}_N(\Delta,b,w,f,c) \ T_1\\
\widehat{\Ll}_N\bigl((12)(\Delta,b,w,f,c)\bigr) \equiv_N
(w_1')^{\frac{N-1}{2}}\ S_1^{-1}\ \widehat{\Ll}_N(\Delta,b,w,f,c)\ T_2
\\ \widehat{\Ll}_N\bigl((23)(\Delta,b,w,f,c)\bigr) \equiv_N
(w_0')^{\frac{1-N}{2}}\ S_2^{-1}\ \widehat{\Ll}_N(\Delta,b,w,f,c)\ S_2
\end{array}$$
\noindent 
where $w_i'=(w_i)_{f_i-*_bc_i}'$, $T_1=T \otimes 1$ etc., and we write
$\equiv_N$ for the equality up to multiplication by $N$th roots of
unity. Moreover, for any enriched non degenerate $\Ii$--tetrahedron
$(\Delta,b,w,a)$, these identities hold true for $(\Delta,b,u,a)$ with
$u$ sufficiently close to $w$ if and only if $a=f-*_bc$ is a flat/charge for $(\Delta,b,w)$.
\end{prop}

To prove this result, we use some formulas described
in the Appendix. The ambiguity up to $N$th roots of unity is due to
Lemma \ref{factor1} (i) and (iii) and a similar
identity. In Figure \ref{tressages} we show these symmetry relations (up
to scalars) by using the graphic tensors of Subsection \ref{matdefsec};
there we put $a=(f,c)$, $\bar{T}=T^{-1}$ and $\bar{S}=S^{-1}$. Note
that the matrices $T$ and $S$ and their inverses act as duality
morphisms. We need the following inversion formula:

\begin{figure}[ht]
\begin{center}
{\tiny
\psfrag {01}{$(01)(w,a)$}
\psfrag {12}{$(12)(w,a)$}
\psfrag {23}{$(23)(w,a)$}
\psfraga <-2pt, 0pt> {w,a}{$w,s$}
\psfraga <-1pt, 1pt> {T}{$T$}
\psfraga <-2pt, 1pt> {S}{$S$}
\includegraphics[width=10cm]{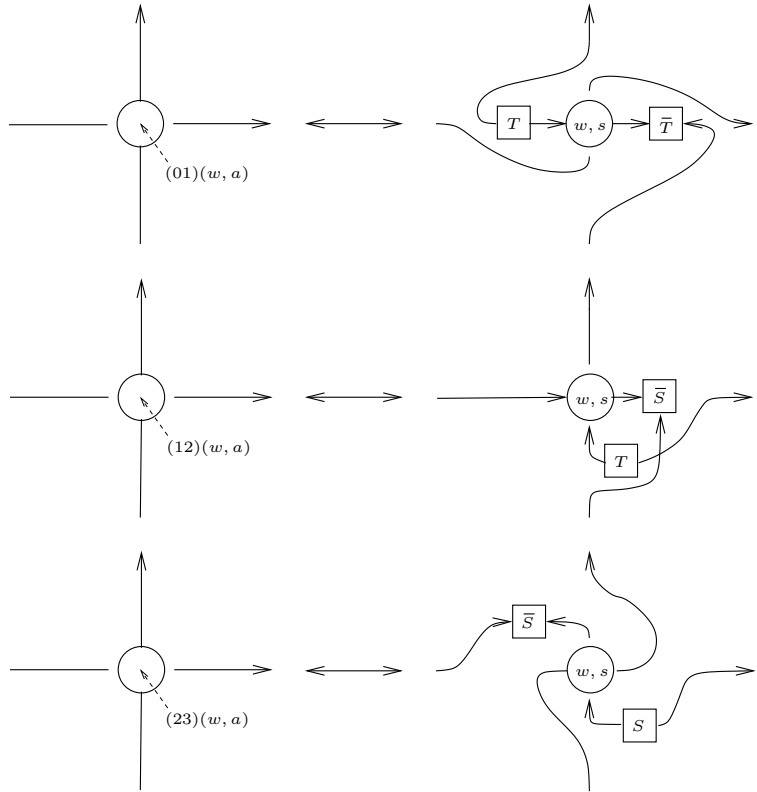}}
\caption{\label{tressages} The symmetry relations of $\widehat{\Ll}_N$
for flat/charged $\Ii$--tetrahedra}
\end{center}
\end{figure}

\begin{lem}\label{inversion} 
Suppose that $w_0'w_1'w_2'=-\tau$. Then we have
$$\prod_{j=1}^n\frac{(w_1')^{-1}}{1-w_0'\zeta^j} \ \cdot
\prod_{j=1}^{N-n}\frac{w_2'}{1-(w_0')^{-1}\zeta^{j-1}} =
\zeta^{-(m+1)(N-n)(N-n-1)}\ \tau^{-n}.$$
\end{lem}

{\bf Proof}\qua  This consists of the following straightforward
computation:
$$\begin{array}{l} \prod_{j=1}^n\frac{(w_1')^{-1}}{1-w_0'\zeta^j} \
\cdot \prod_{j=1}^{N-n}\frac{w_2'}{1-(w_0')^{-1}\zeta^{j-1}} =
\frac{w_2(w_1')^{-n}(w_2')^{-n}}{\prod_{j=1}^n(1-w_0'\zeta^j)
\prod_{j=1}^{N-n}(1-(w_0')^{-1}\zeta^{j-1})}
\\ \\
\hspace{2cm} =
\frac{w_2(w_1')^{-n}(w_2')^{-n}(-1)^{N-n}}{\zeta^{(m+1)(N-n)(N-n-1)}
(w_0')^{n-N}\prod_{j=1}^n(1-w_0'\zeta^j)\prod_{j=1}^{N-n}(1-w_0'\zeta^{1-j})}
\\ \\
\hspace{2cm} = - \zeta^{-(m+1)(N-n)(N-n-1)}\ \tau^{-n}\ w_2w_0\
(1-w_0)^{-1} \\
\hspace{2cm} = \zeta^{-(m+1)(N-n)(N-n-1)}\ \tau^{-n}.\hfill\qed
\end{array}$$

\medskip

{\bf Proof of Proposition \ref{6jsym}}\qua  By (\ref{dilogtet}) and
Proposition \ref{unitarity} we have
$$\begin{array}{l}
\widehat{\Ll}_N\bigl((01)(\Delta,b,w,f,c)\bigr)_{k,j}^{i,l} =
(\widehat{\Ll}_N)^{-1}((w_0')^{-1},w_2')_{k,j}^{i,l} \hspace{4.5cm}\\
 \hspace{1cm} = [(w_0')^{-1}]\ \frac{g(1)}{g((w_0')^{-1})}\ \zeta^{-ij
-(m+1)i^2}\ \delta(k+j-l)\
\prod_{s=1}^{k-i}\frac{1-(w_0')^{-1}\zeta^{s-1}}{w_2'} \\
\hspace{1cm}\equiv_N (w_0')^{\frac{1-N}{2}} \ \frac{g(w_0')}{g(1)}\
\zeta^{-ij -(m+1)i^2}\ \delta(k+j-l)\ \zeta^{(m+1)(i-k)^2}\
\prod_{s=1}^{i-k}\ \frac{(w_1')^{-1}}{1-w_0'\zeta^{s}} \\
\hspace{1cm}\equiv_N (w_0')^{\frac{1-N}{2}} \ \frac{g(w_0')}{g(1)}\
\zeta^{-il +(m+1)k^2}\ \delta(k+j-l)\ \prod_{s=1}^{i-k}\
\frac{(w_1')^{-1}}{1-w_0'\zeta^{s}}
\end{array}$$
where we use Lemma \ref{factor1} (i) and Lemma \ref{inversion} in the 
second equality ($\tau=\zeta^{-(m+1)}$). Now, this may be written as
$$\begin{array}{l} (w_0')^{\frac{1-N}{2}} \
\sum_{i',k'=0}^{N-1}\zeta^{-(m+1)(i')^2} \delta(i+i')\
\zeta^{(m+1)(k')^2}\delta(k+k') \hspace{5cm} \\
\hspace{4cm} \times \ \frac{g(w_0')}{g(1)}\
\zeta^{i'l+(m+1)(i')^2}\ \delta(j-k'-l)\ \prod_{s=1}^{k'-i'}\
\frac{(w_1')^{-1}}{1-w_0'\zeta^{s}}\\
\hspace{1cm} = (w_0')^{\frac{1-N}{2}} \ \sum_{i',k'=0}^{N-1}
(T^{-1})_i^{i'} \ \widehat{\Ll}_N(w_0',(w_1')^{-1})_{i',j}^{k',l}\ \
T_{k'}^k
\end{array}$$
which proves the first relation. The third comes from a very similar
computation:
$$\begin{array}{l}
\widehat{\Ll}_N\bigl((23)(\Delta,b,w,f,c)\bigr)_{i,l}^{k,j} =
(\widehat{\Ll}_N)^{-1}((w_0')^{-1},w_2')_{i,l}^{k,j} \hspace{4.5cm}\\
 \equiv_N (w_0')^{\frac{1-N}{2}} \
 \frac{g(w_0')}{g(1)}\ \zeta^{-kj +(m+1)i^2}\ \delta(i+l-j)\
 \prod_{s=1}^{k-i}\ \frac{(w_1')^{-1}}{1-w_0'\zeta^{s}}\\
  \equiv_N (w_0')^{\frac{1-N}{2}} \
 \frac{g(w_0')}{g(1)}\ \bigl(N^{-1} \sum_{l'=0}^{N-1}
 \zeta^{l'(i+l-j)} \bigr)\ \zeta^{-kj +(m+1)i^2}\ \delta(i+l-j)
 \prod_{s=1}^{k-i}\ \frac{(w_1')^{-1}}{1-w_0'\zeta^{s}}
\end{array}$$
$$\begin{array}{l}
 \hspace{0.5cm} \equiv_N (w_0')^{\frac{1-N}{2}} \ \sum_{j',l'=0}^{N-1}
 (N^{-1/2}\zeta^{ll'})\ (N^{-1/2}\zeta^{-jj'})\ \\
\hspace{4cm}\times \ \frac{g(w_0')}{g(1)}\ \zeta^{il' +(m+1)i^2}\
\delta(j'-k-l')\ \prod_{s=1}^{k-i}\
\frac{(w_1')^{-1}}{1-w_0'\zeta^{s}}\\
 \hspace{0.5cm} \equiv_N (w_0')^{\frac{1-N}{2}} \ \sum_{j',l'=0}^{N-1}
 (S^{-1})_{j}^{j'}\ \widehat{\Ll}_N(w_0',(w_1')^{-1})_{i,j'}^{k,l'}\
 S_{l'}^{l}.
\end{array}$$
where we use again Lemma \ref{factor1} (i) and Lemma
\ref{inversion} in the first equality. The second symmetry relation is
more sophisticated. Consider the right-hand side. It gives
$$\begin{array}{l} \sum_{i',l'=0}^{N-1}(S^{-1})_{i}^{i'}\
 \widehat{\Ll}_N(w_0',(w_1')^{-1})_{i',j}^{k,l'}\ T_{l'}^{l}
 \hfill\\ = \sum_{i',l'=0}^{N-1}
 \bigl(N^{-1/2}\zeta^{-ii'}\bigr)\
 \bigl(\zeta^{(m+1)(l')^2}\delta(l+l')\bigr)\hfill\\ 
 \hspace {4cm} \times \ \frac{g(w_0')}{g(1)}\ \zeta^{i'l' +(m+1)(i')^2}\delta(k+l'-j)\
 \prod_{s=1}^{k-i'}\frac{(w_1')^{-1}}{1-w_0'\zeta^{s}}\\ =
 N^{-1/2}\zeta^{(m+1)l^2}\ \frac{g(w_0')}{g(1)}\ \delta(j+l-k)\
 \sum_{i'=0}^{N-1} \zeta^{-i'l+(m+1)(i')^2-ii'}\
 \prod_{s=1}^{k-i'}\ \frac{(w_1')^{-1}}{1-w_0'\zeta^{s}}\\ =
 N^{-1/2}\zeta^{(m+1)(l^2-k^2)}\ \frac{g(w_0')}{g(1)}\
 \delta(j+l-k)\ \sum_{i'=0}^{N-1} \zeta^{i'(k-i-l)}\
 \prod_{s=1}^{i'-k}\ \frac{1-(w_0')^{-1}\zeta^{s-1}}{w_2'}\\ =
 \frac{g(w_0')}{g(1)}\ N^{-1/2}\zeta^{(m+1)(l^2-k^2)+k(k-i-l)}\
 \delta(j+l-k)\\
\hspace{4cm}\times \ \sum_{i'=0}^{N-1}
\bigl(\zeta^{k-i-l}(w_2')^{-1}\bigr)^{i'-k} \ \prod_{s=1}^{i'-k}\
\bigl(1-(w_0')^{-1}\zeta^{s-1}\bigr)\\ =
\frac{g(w_0')}{g(1)}\ N^{-1/2}\zeta^{(m+1)(l^2-k^2)+k(k-i-l)}\
 \delta(j+l-k)\ f(0,(w_0')^{-1}\zeta^{-1}\vert
(w_2')^{-1}\zeta^{j-i})
\end{array}$$
where the function $f(x,y\vert z)$ is defined in the proof of
Proposition \ref{unitarity}. As described there, we have
$$\frac{f(0,(w_0')^{-1}\zeta^{-1}\vert
(w_2')^{-1}\zeta^{j-i})}{f(0,(w_0')^{-1}\zeta^{-1}\vert (w_2')^{-1})}
= \prod_{s=1}^{j-i} \
\frac{1-(w_2')^{-1}\zeta^{s-1}}{w_1'\zeta^{-(m+1)}\zeta^s}$$
where we note that
$-(w_0')^{-1}(w_2')^{-1}\zeta^{-1}=w_1'\zeta^{-(m+1)}$. Simplifying
the powers of $\zeta$ with the help of the Kronecker symbol
$\delta(j+l-k)$, we see immediately that the right-hand side reads
$$N^{-1/2}\ \frac{g(w_0')}{g(1)}\ f(0,(w_0')^{-1}\zeta^{-1}\vert
(w_2')^{-1}) \ \zeta^{-il-(m+1)i^2}\ \delta(j+l-k)\
\prod_{s=1}^{j-i} \ \frac{1-(w_2')^{-1}\zeta^{s-1}}{w_1'}$$
\noindent A very similar computation to Lemma \ref{factor1} (iii) shows that
$$f(0,y \vert z) \equiv_N \frac{(-yz)^{\frac{N-1}{2}}\
g(1)}{g(y^{-1}/\zeta)g(z/\zeta)}$$
\noindent for $z^N=1/(1-y^N)$. This gives
$$\begin{array}{lll} f(0,(w_0')^{-1}\zeta^{-1}\vert (w_2')^{-1}) &
\equiv_N & \frac{(-(w_0')^{-1}(w_2')^{-1})^{\frac{N-1}{2}}\
g(1)}{g(w_0')\ g((w_2')^{-1}/\zeta)} \hspace{2cm} \\ \\ & \equiv_N &
\frac{(w_1')^{\frac{N-1}{2}}\ g(1)}{g(w_0')g((w_2')^{-1})}\ \
\frac{w_1'}{{1-(w_2')^{-1}}}
\end{array}$$
\noindent 
where we use Lemma \ref{gdilocyc} in the last equality. Hence, noting
that $\vert g(1) \vert =N^{1/2}$, we find
$$\begin{array}{lll} N^{-1/2}\ \frac{g(w_0')}{g(1)}\
f(0,(w_0')^{-1}\zeta^{-1} \vert (w_2')^{-1}) & \equiv_N &
N^{-1}(w_1')^{\frac{N-1}{2}}\frac{w_1'}{1-(w_2')^{-1}}\
\frac{g(1)}{g((w_2')^{-1})}\\ \\ & \equiv_N & (w_1')^{\frac{1-N}{2}}\
\frac{g(1)\ [(w_2')^{-1}]}{g((w_2')^{-1})}
\end{array}$$
\noindent with $[x]=N^{-1} (1-x^N)/(1-x)$ as in Proposition \ref{unitarity}. So
$$\begin{array}{l} \sum_{j',k'=0}^{N-1}(S^{-1})_{i}^{i'}\
\widehat{\Ll}_N(w_0',(w_1')^{-1})_{i',j}^{k,l'}\ T_{l'}^{l}
\hspace{6cm} \\
\hspace{1cm}\equiv_N (w_1')^{\frac{1-N}{2}} \ \frac{g(1)\ [(w_2')^{-1}]}
{g((w_2')^{-1})}\ \zeta^{-il-(m+1)i^2}\ \delta(j+l-k)\ \prod_{s=1}^{j-i} 
\ \frac{1-(w_2')^{-1}\zeta^{s-1}}{w_1'}\\ \\
\hspace{1cm} = (w_1')^{\frac{1-N}{2}}\
(\widehat{\Ll}_N)^{-1}((w_2')^{-1}, w_1')_{j,l}^{i,k}\quad .
\end{array}$$
\noindent This proves the second symmetry relation. The last claim
follows from the fact that we need $\tau=\zeta^{-(m+1)}$ in all the
above computations.\endproof

Let us recall the matrix $\Rr_N$ from Subsection \ref{matdefsec}:
\begin{defi}\label{symRNdef}
{\rm For each odd integer $N>1$, the {\it symmetrized matrix dilogarithm} of flat/charged
$\Ii$--tetrahedra is defined as}
\begin{eqnarray}
\Rr_N(\Delta,b,w,f,c) & = &
\bigl((w_0')^{-c_1}(w_1')^{c_0}\bigr)^{\frac{N-1}{2}}\
(\widehat{\Ll}_N)^{*_b}(w_0;f_0-*_bc_0,f_1-*_bc_1) 
 \nonumber \\
 & = & 
\bigl((w_0')^{-c_1}(w_1')^{c_0}\bigr)^{\frac{N-1}{2}}\ (\Ll_N)^{*_b}
(w_0',(w_1')^{-1})\nonumber
\end{eqnarray}
\noindent {\it where $w_i'=(w_i)_{f_i-*_bc_i}'$.}
\end{defi}
Note that the log
of the scalar $(w_0')^{-c_1}(w_1')^{c_0}$ is of the same form as the
function we add to the classical Rogers dilogarithm to define the
uniformized dilogarithm ${\rm R}(x;p,q)$ in (\ref{Ndef}). We have the
following result, which is the precise form of (1) in Theorem \ref{mainlocal}: 

\begin{cor}\label{6jsymbis} 
For any flat/charged $\Ii$--tetrahedron $(\Delta,b,w,f,c)$ with $*_b=+1$,
the following symmetry relations hold true:
$$\begin{array}{l} \Rr_N\bigl((01)(\Delta,b,w,f,c)\bigr)
\equiv_N \pm\ T_1^{-1}\ \Rr_N(\Delta,b,w,f,c) \ T_1\\
\Rr_N\bigl((12)(\Delta,b,w,f,c)\bigr) \equiv_N \pm\
S_1^{-1}\ \Rr_N(\Delta,b,w,f,c)\ T_2 \\
\Rr_N\bigl((23)(\Delta,b,w,f,c)\bigr) \equiv_N \pm\
S_2^{-1}\ \Rr_N(\Delta,b,w,f,c)\ S_2 \quad .
\end{array}$$
where $\equiv_N$ means equality up to multiplication by $N$th roots of unity.
\end{cor}

{\bf Proof}\qua For the first equality, on the
left hand side the scalar factor in Definition \ref{symRNdef} reads
$$\left((w_0')^{-1})^{-c_2}((w_2')^{-1})^{c_0}\right)^{\frac{N-1}{2}} =\left((w_0')^{-c_1+1}((w_0'w_2')^{-c_0})\right)^{\frac{N-1}{2}}$$
because $c_0 + c_1 +c_2 = 1$. As $w_0'w_1'w_2'=-\zeta^{-(m+1)}$, this
is equal to $(w_0')^{\frac{N-1}{2}}$ times
$((w_0')^{-c_1}((w_1')^{c_0}))^{\frac{N-1}{2}}$ up to sign and
multiplication by $N$th roots of unity. By Proposition \ref{6jsym},
this coincides with the scalar factor on the right hand side. The same
computation works for the other transpositions.\endproof

\subsection{Complete five term relations}\label{compqrel}

Recall from Subsection \ref{transitdefsec} the notion of
{\it flat/charged} $\Ii$--transit. It is immediate that a flat/charged
$\Ii$--transit associated to the branching configuration of Figure
\ref{CQDidealt} defines an $N$th--branch transit (see Definition \ref{defenricht}) for the induced enriched $\Ii$--tetrahedra. Finally
we can prove the statement (2) of Theorem \ref{mainlocal}:

\begin{teo} \label{qtransit} For every odd integer $N>1$, the tensors
  $\Rr_N$ satisfy the five-term identities associated to arbitrary
  flat/charged $\Ii$--transits with no constraints on the cross ratio
  moduli, possibly up to sign and multiplication by $N$th roots of
  unity. Moreover, the flat/charged $\Ii$--transits realize the minimal
  relations between enriched $\Ii$--tetrahedra for this property to hold.
\end{teo}

By Proposition \ref{6jsym}, the last claim shows that having the {\it whole
  set} of five-term identities is equivalent to having the symmetry
relations (and similarly when $N=1$ after replacing the flat/charge
with flattenings). On the other hand, the analyticity of
$\Rr_N(\Delta,b,w,a)$ for all branchings $b$ implies only $a=f+\lambda
c$ for an arbitrary $\lambda \in \mc$, where $f$ is a flattening and
$c$ an integral charge. This, in turn, is enough for the matrix Schaeffer's
identity.

\medskip

{\bf Proof}\qua  Consider first a transit whose branching is as
in Figure \ref{CQDidealt}. By Theorem \ref{basNtr} and Definition \ref{symRNdef},
in that case the statement is true if we remove the powers
of $(w_0')^{-c_1}(w_1')^{c_0}$ from the $\Rr_N$. We claim
that the products $P_0$ and $P_1$ of these scalars at both sides of
the transit are also equal, up to a sign if $N \equiv 3$ mod($4$). 

\noindent Denote by $\Delta^j$
the tetrahedron opposite to the $j$-th vertex (for the ordering of the
vertices induced by the branching); do the same for the log-branches
and integral charges. Recall the notation ${\rm l}_i^j$ for the
log-branch of the $i$-th edge of $\Delta^j$.  Then $P_0$ reads
$$\exp\biggl(\frac{N-1}{2}\bigl(-c_1^1{\rm l}_0^1+
c_0^1{\rm l}_1^1-c_1^3{\rm l}_0^3 +c_0^3{\rm l}_1^3\bigr)\biggr)$$
\noindent 
times $(-1)^{\frac{N-1}{2}(f_0^1c_1^1+f_1^1c_0^1+f_0^3c_1^3+f_1^3c_0^3)}$, 
and $P_1$ is
$$\exp\biggl(\frac{N-1}{2}\bigl(-c_1^0{\rm l}_0^0+c_0^0{\rm
l}_1^0-c_1^2{\rm l}_0^2+c_0^2{\rm l}_1^2-c_1^4{\rm l}_0^4+c_0^4{\rm
l}_1^4\bigr) \biggr)$$
\noindent 
times
$(-1)^{\frac{N-1}{2}(f_0^0c_1^0+f_1^0c_0^0+f_0^2c_1^2+
f_1^2c_0^2+f_0^4c_1^4+f_1^4c_0^4)}$. Consider
the sums in the above exponentials. They are formally the same as
minus the {\it sums of logarithmic terms} at both sides of
(\ref{Rtidsec4}) (just replace each log of a modulus by the
corresponding log-branch, and similarly each flattening by a
charge). The proof of Proposition \ref{Rtransit} tells us that when the
log of the moduli satisfy (\ref{split}), these
sums are equal when the relations (\ref{sumtransitforf}) are
true. Since the relations (\ref{ideqfdefsec}) are formally
identical to (\ref{split}), and the transit of integral charges is
defined by the identities (\ref{sumtransitforf}), we deduce that
$P_0=\pm P_1$. So the statement is true for a flat/charged
$\Ii$--transit with a branching as in Figure \ref{CQDidealt}.

We get the result for all the other instances of $2 \leftrightarrow 3$
flat/charged $\Ii$--transits by using Lemma \ref{6jsymbis}. Indeed,
note that the matrices $S$ and $T$ always act on the input spaces of
the tensors $\Rr_N$, whereas their inverses always act on the output
spaces. As the composition of the $\Rr_N$ is defined by the {\it
  oriented} graphs described in Subsection \ref{matdefsec}, the matrix
action compensates on the common faces of two tetrahedra. So the five-term relation corresponding to any $2 \leftrightarrow 3$ flat/charged
$\Ii$--transit is conjugated to the one for Figure \ref{CQDidealt}.

By the proof of Theorem \ref{basNtr}, we know that the lifted matrix
Schaeffer's identity for $\widehat{\Ll}_N$ holds true if and only if
the relations (\ref{ideqmodprime}) are valid. This implies that: for
each enriched $\Ii$--tetrahedra $(\Delta,b,w,a)$ we have
$w_0'w_1'w_2'=-\tau$ for a fixed root of unity $\tau$, each $a$ must
be of the form $a=f+\lambda c$, $\lambda \in \mz$, for some flattening
$f$ and integral charge $c$, and $f$ and $c$ transit as usual. We
claim that the whole set of five-term identities hold simultaneously
if and only if $\lambda = -*_b$. For instance, when $\tau=1$ (that is
$\lambda = 0$, the case of flattened $\Ii$--tetrahedra), we see easily
that the first relation in Proposition \ref{6jsym} becomes
$$\widehat{\Ll}_N\bigl((01)(\Delta,b,w,f)\bigr) \equiv_N M_1\
(w_0')^{\frac{1-N}{2}}\ T_1^{-1}\ \widehat{\Ll}_N(\Delta,b,w,f) \ T_1\
M_1^{-1}$$ where $M$ is the $N \times N$--matrix with entries $M_j^i =
\zeta^{(m+1)i}\ \delta(i-j)$. The second relation becomes even more
complicated. The same phenomenon happens whenever we do not put
$a=f-*_b c$. In general there is no global compensation of these
further matrix actions when considering symmetries on the lifted
matrix Schaeffer's identity for flattened $\Ii$--tetrahedra. So we are
forced to consider flat/charged $\Ii$--tetrahedra to get the versions
of the lifted matrix Schaeffer's identity for all instances of
enriched $\Ii$--transits. \endproof

As promised in the Introduction, the next lemma shows that for $N >1$
the tensors $\Rr_N$ coincide, up to an $N$th root of unity, with the
symmetrized quantum dilogarithms used in \cite{BB2} (see Definition
3.2 in that paper, where the $N$th root moduli $\underline{w}_i'$
below are given by $N$th roots of certain ratios of cocycle
parameters, as explained in Remark \ref{compareQHI}).

\begin{lem}\label{idemQHI} 
For any flat/charged $\Ii$--tetrahedron we have
$$\begin{array}{l}
\Rr_N(\Delta,b,w,f,c)_{k,l}^{i,j} = \bigl((w_0')^{-c_1}(w_1')^{c_0}
\bigr)^{\frac{N-1}{2}} \ \zeta^{(m+1)c_1(i-k)-(m+1)^2f_1c_0}\hspace{4cm} \\
\hspace{3.5cm} \times \ \frac{g(\underline{w}_0')}{g(1)}\
\zeta^{kj+(m+1)k^2}\ \delta(i+j-l)\ \prod_{s=1}^{i-k-(m+1)c_0}\
\frac{(\underline{w}_1')^{-1}}{1-\underline{w}_0'\zeta^s}
\end{array}$$
where $\underline{w}_n'= (w_n)_{f_n}'$ and $w_n'= (w_n)_{f_n-*_bc_n}'$ for $n=0,1,2$.
\end{lem}

{\bf Proof}\qua Consider the right hand side. We have
$$\zeta^{(m+1)c_1(i-k)}\ \prod_{s=1}^{i-k-(m+1)c_0}\
\frac{(\underline{w}_1')^{-1}}{1-\underline{w}_0'\zeta^s} =
\prod_{s=1}^{i-k}\ \frac{(w_1')^{-1}}{1-w_0'\zeta^s}\
\prod_{s=1}^{-(m+1)c_0}\
\frac{(\underline{w}_1')^{-1}}{1-\underline{w}_0'\zeta^s}$$
\noindent and Lemma \ref{gdilocyc} implies
$$\frac{g(\underline{w}_0')}{g(1)}\ \prod_{s=1}^{-(m+1)c_0}\
\frac{(\underline{w}_1')^{-1}}{1-\underline{w}_0'\zeta^s} =
\frac{g(w_0')}{g(1)}\ \zeta^{(m+1)^2f_1c_0}.$$ Gathering these
formulas we find the result. \endproof

\begin{remark}\label{compareQHI}{\rm 
    In \cite{BB2} we started
    with the Faddeev-Kashaev's matrix of $6j$--symbols for the cyclic
    representation theory of a Borel quantum subalgebra
    $\mathcal{B}_{\zeta}$ of $U_{\zeta}(sl(2,\mc))$, where
    $\zeta=\exp(2i\pi/N)$. We associated this matrix to a tetrahedron
    equipped with a $1$--cocycle taking values in the Borel subgroup
    $B$ of upper triangular matrices of $SL(2,\mc)$. This was possible
    due to a very natural parametrization of the cyclic irreducible
    representations of $\mathcal{B}_{\zeta}$ (see Remark
    \ref{rempar}). Only {\it a posteriori} we noticed that the
    relevant parameters were certain ratios of the matrix entries
    of the cocycle values, that corresponded to the
    cross-ratio moduli of determined hyperbolic ideal
    tetrahedra. We
    used the {\it idealization} procedure to transfer in the set up of
    $\Ii$--triangulations the results and the computations
    previously obtained in terms of cocycle parameters.

    \noindent However, the fundamental objects were just the quantum
    basic matrix dilogarithms $\Ll_N$. Studying them directly is far
    from a mere rephrasing of the results in \cite{BB2}. As a
    by-product, here we show that the symmetrization is intrisically
    related to the algebra $\mathcal{B}_\zeta$. Also, in \cite{BB2} a
    choice of flattening was {\it hidden} in a preliminary choice at
    hand of the $N$th roots of the cocycle parameters; the integral
    charge was intended as a way to contruct {\it another} tensor from
    $\Ll_N$, rather than a way to evaluate the {\it same} matrix
    valued function on different {\it variable} branches of the
    Riemann surface $\widehat{\mc}$, globally organized by the
    combination of flattenings and charges.}
\end{remark}

\section{Classical and quantum dilogarithmic invariants}\label{GTINV}

In this section we define the state sum invariants based
on the matrix dilogarithms $\Rr_N$. As this construction is a
generalization and a refinement of the QHI's one, we shall often refer
to \cite{BB2}. For the sake of clarity, we begin with some general
facts about $3$--manifold triangulations.

\medskip

{\bf A few generalities on 3--manifold triangulations}\qua
Consider a tetrahedron $\Delta$ with its usual triangulation with $4$
vertices, and let $\Gamma_\Delta$ be the interior of the $2$--skeleton
of the dual cell decomposition. A {\it simple} polyhedron $P$ is a
$2$--dimensional compact polyhedron such that each point of $P$ has a
neighbourhood which can be embedded into an open subset of
$\Gamma_\Delta$.  A simple polyhedron $P$ has a natural stratification
given by its singularities; $P$ is {\it standard} if all the strata of
this stratification are open cells; depending on the dimension, we
call them {\it vertices, edges} and {\it regions}.

Every compact $3$--manifold $Y$ (which for simplicity we assume
connected) with non-empty boundary has {\it standard spines}
\cite{Cas}, that is standard polyhedra $P$ together with an embedding
in Int($Y$) such that $Y$ is a regular neighbourhood of $P$.
Moreover, $Y$ can be reconstructed from any of its standard spines.
Since we shall always work with combinatorial data encoded by
triangulations/spines, which define the corresponding manifold only up
to PL-homeomorphisms, we shall
systematically forget the underlying embeddings.

A {\it singular triangulation} of a 3--dimensional polyhedron $Q$ is a
triangulation for which self-adjacencies and multiple adjacencies of
$3$--simplices along $2$--faces are allowed. This shall be simply called
a {\it triangulation} of $Q$. For any $Y$ as above, let us denote by
$Q(Y)$ the space obtained by collapsing each connected component of
$\partial Y$ to a point. A (topological) {\it ideal triangulation} of
$Y$ is a triangulation $T$ of $Q(Y)$ such that the vertices of $T$ are
precisely the points of $Q(Y)$ corresponding to the components of
$\partial Y$. By removing small open neigbourhoods of the vertices of
$Q(Y)$, any ideal triangulation leads to a cell decomposition of $Y$
by {\it truncated tetrahedra}, which induces a (singular)
triangulation on $\partial Y$. If $Y$ is oriented, the tetrahedra of
any triangulation inherit the induced orientation. From now on we will
consider only oriented manifolds.

For any ideal triangulation $T$ of $Y$, the $2$--skeleton of the {\it
dual} cell-decomposition of $Q(Y)$ is a standard spine $P(T)$ of
$Y$. This procedure can be reversed, so that we can associate to each
standard spine $P$ of $Y$ an ideal triangulation $T(P)$ of $Y$ such
that $P(T(P))=P$. Hence standard spines and ideal triangulations are
dual equivalent viewpoints which we will freely intermingle.

Consider now a compact closed $3$--manifold $W$. For any
$r_0\geq 1$, let $Y$ be the manifold
obtained by removing $r_0$ disjoint open balls from $W$. By definition
$Q(Y)=W$ and any ideal triangulation of $Y$ is a singular
triangulation of $W$ with $r_0$ vertices. It is easily seen
that all triangulations of $W$ come in this way from ideal
triangulations.

\begin{figure}[ht]
\begin{center}
 \includegraphics[width=6cm]{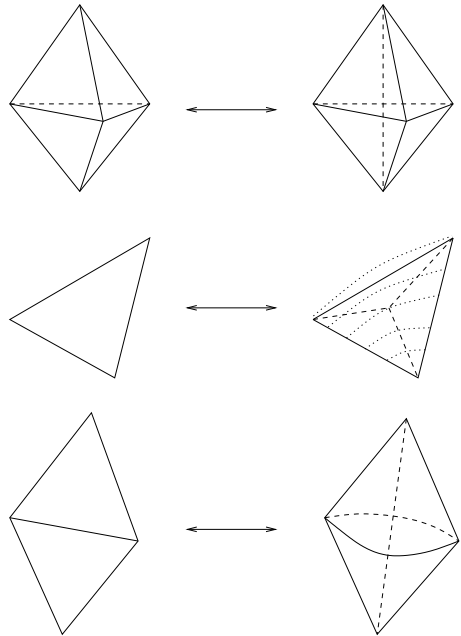}
\caption{\label{CQDfigmove1bis} The moves between singular triangulations} 
\end{center}
\end{figure}

\begin{figure}[ht]
\begin{center}
 \includegraphics[width=6cm]{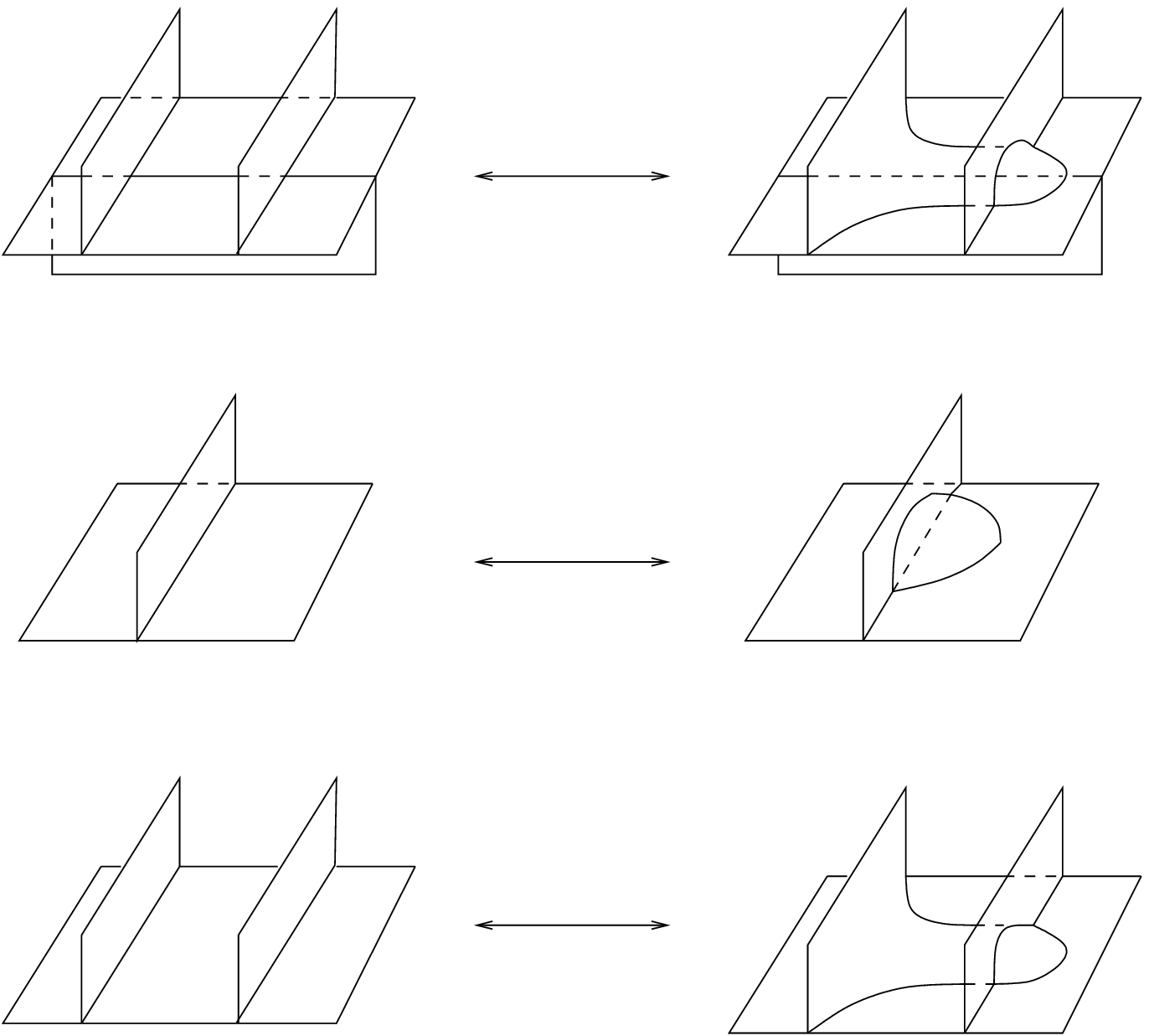}
\caption{\label{CQDfigmove2} The moves on standard spines} 
\end{center}
\end{figure}

In Figure \ref{CQDfigmove1bis} and Figure \ref{CQDfigmove2} we recall
some elementary moves on the triangulations and simple spines of a
polyhedron $Q(Y)$. They are called the (respectively dual) $2\to 3$
move, bubble move, and $0 \to 2$ move. The bubble move consists in
replacing a $2$--simplex by the cone on a $2$--sphere triangulated by
two $2$--simplices. It is a fact (see \cite{Mat,
  Pi}) that standard spines of the same compact oriented $3$--manifold
$Y$ with boundary and with at least two vertices (which, of course, is
a painless requirement) may always be connected by means of a finite
sequence of the dual $2 \to 3$ move and its inverse. The dual result
holds for ideal triangulations.  In order to handle triangulations of
closed $3$--manifolds we need a move which allows us to vary the
number of vertices, like the bubble move.

We say that a triangulation $T$ of a compact closed
$3$--manifold $W$ is {\it quasi regular} if all the edges have distinct
vertices. Every $W$ admits quasi-regular triangulations. 

\medskip

{\bf Global branching}\qua A global branching on an ideal
triangulation $T$ of $Y$ (see Section \ref{ftdefsec}) can be defined
in terms of the dual orientations of the regions of the standard spine
$P(T)$, which are dual to the edges of $T$. In fact, a global
branching gives the spine $P(T)$ a structure of oriented {\it branched
  surface} embedded in $Y$ (this also justifies the term ``branching''
- this point of view is widely developed in \cite{BP2}). Every $Y$ has
triangulations supporting a branching \cite[Theorem 3.4.9]{BP2}.

We introduced in Subsection \ref{ITdefsec} explicitely described
$\Ii$--{\it triangulations} for compact closed pairs $(W,\rho)$ or
gentle cusped manifolds $M$. This notion extends
straightforwardly to (topological) ideal triangulations of any compact
oriented $Y$.

It is easy to see that the $S_4$--action on $\Ii$--tetrahedra defined
before Theorem \ref{mainlocal} extends to $\Ii$--triangulations: a
global change of the branching turns any $\Ii$--triangulation $(T,b,w)$
into another $\Ii$--triangulation $(T,b',w')$.

\medskip

{\bf Pseudo-developing maps}\qua We describe here an important
geometric object associated to any $\Ii$--triangulation, that clarifies
the role of the edge compatibility condition (\ref{ideq1defsec}). Given an
$\Ii$--triangulation $(T,b,w)$ of $Y$, lift $T$ to a cellulation
$\widetilde{T}$ of the universal covering $\widetilde{Y}$, and fix a
base point $\tilde{x}_0$ in the $0$--skeleton of $\widetilde{T}$;
denote by $x_0$ the projection of $\tilde{x}_0$ onto $Y$. For each
tetrahedron in $\widetilde{T}$ that contains $\tilde{x}_0$, use the
moduli of its projection in $T$ to define an hyperbolic ideal
tetrahedron, by respecting the gluings in $\widetilde{T}$. Doing
similarly with the vertices adjacent to $\tilde{x}_0$ and so on, we
construct an image in $\bar{\mh}^3$ of a complete lift of $T$ in
$\widetilde{T}$, having one tetrahedron in each $\pi_1(Y,x_0)$--orbit.
The edge compatibily implies that for any two paths of tetrahedra in
$\widetilde{T}$ having a same starting point, we get the same end
point.  This construction extends to a piecewise-linear map $D\co 
\widetilde{Y} \rightarrow \bar{\mh}^3$, equivariant with respect to
the action of $\pi_1(Y,x_0)$ and $PSL(2,\mc)$. So we eventually find:
a representation $\tilde{\rho}\co  \pi_1(Y,x_0) \rightarrow PSL(2,\mc)$
satisfying $D(\gamma(x))=\tilde{\rho}(\gamma) D(x)$ for each $\gamma
\in PSL(2,\mc)$; a piecewise-straight continuous section of the flat
bundle $\widetilde{Y} \times_{\tilde{\rho}} \bar{\mh}^3 \rightarrow
Y$, with structural group $PSL(2,\mc)$ and total space the quotient of
$\widetilde{Y} \times \bar{\mh}^3$ by the diagonal action of
$\pi_1(Y,x_0)$ and $\tilde{\rho}$. The map $D$ behaves formally as a
developing map for a $(PSL(2,\mc),\mh^3)$--structure on $Y$ (see eg
\cite[Ch. B]{BP1} for this notion). By the arbitrary choices we made,
only the conjugacy class $\rho$ of $\tilde{\rho}$ is well defined, and
$D$ is only defined up to the left action of $PSL(2,\mc)$. This class
is preserved by $\Ii$--transits.

\medskip

{\bf Flat/charged $\Ii$--transits}\qua In Subsection
\ref{transitdefsec} we defined the notion of a $2\leftrightarrow 3$
flat/charged $\Ii$--transit. We need to discuss this notion for the
other triangulation local moves. We stipulate that whenever they apply
to a global triangulation, everything remains unchanged outside the
support of the given move. For clarity we repeat the definitions also
for the $2\to 3$ move.  

An {\it $\Ii$--transit} $(T,b,w) \rightarrow (T',b',w')$ of
$\Ii$--triangulations of the same $3$--manifold $Y$ consists of a bare
triangulation move $T \to T'$ that extends to a branching move
$(T,b)\to (T',b')$, ie, the two branchings coincide on the `common'
edges of $T$ and $T'$; moreover the modular triples have the following
behaviour.

For a $2\to 3$ move we require that for each common edge $e\in
   E(T)\cap E(T')$ we have
\begin{equation}\label{ideqmod} \prod_{a\in \epsilon_T^{-1}(e)}w (a)^*=
 \prod_{a'\in \epsilon_{T'}^{-1}(e)}w' (a')^*
\end{equation}
where $*=\pm 1$ according to the $b$--orientation of the abstract
tetrahedron containing $a$ (respectively $a'$).

For the $0 \to 2$ and bubble move we require that for each edge $e\in
E(T')$ we have
\begin{equation}
\label{ideq2mod} \prod_{a'\in \epsilon_{T'}^{-1}(e)}w' (a')^*=1.\end{equation}
\noindent The $\Ii$--transits for negative $3 \to 2$ moves are defined
in the same way, and for negative $2 \to 0$ and bubble moves $w'$ is
defined by simply forgetting the moduli of the two disappearing
tetrahedra.

Consider a $2 \to 3$ $\Ii$--transit $(T,b,w) \rightarrow (T',b',w')$ as
in Figure \ref{CQDidealt}. Give a flattening to each tetrahedron of the
initial configuration, and denote by ${\rm l}\co  E_{\Delta}(T)
\rightarrow \mc$ the corresponding log-branch function on $T$.  A map
$f'\co  E_{\Delta}(T')$ $\rightarrow \mz$ defines a $2 \to 3$ {\it
flattening transit} $(T,b,w,f) \rightarrow (T',b',w',f')$ if for each
common edge $e \in E(T) \cap E(T')$ we have 
\begin{equation}\label{ideq} \sum_{a\in \epsilon_T^{-1}(e)}*\ {\rm l}(a)=
 \sum_{a'\in \epsilon_{T'}^{-1}(e)}*\ {\rm l}'(a')
\end{equation}
\noindent {\rm where $*=\pm 1$ according to the $b$--orientation of the
tetrahedron that contains $a$ (respectively $a'$).

A map $f'\co  E_{\Delta}(T')
\rightarrow \mz$ defines a $0 \to 2$ (respectively bubble) \emph{flattening
transit} if for each edge $e \in E(T')$ we have}
\begin{equation}\label{ideq2}
\sum_{a'\in \epsilon_{T'}^{-1}(e)} *\ {\rm l}'(a')=0.
\end{equation}
For negative $2 \to 0$ and bubble moves the
flattening transits are defined by simply forgetting the flattenings
of the two disappearing tetrahedra.

Note that the relations (\ref{ideq2}) mean that the two new tetrahedra
have \emph{the same} log-branches, for their $b$--orientations are
always opposite (and similarly for (\ref{ideq2mod})).

We saw in Subsection \ref{ITdefsec} that we need to fix an (arbitrary)
non empty link $L$ in $W$ in order to remove an obstruction to the
existence of global integral charges on triangulations of $W$. So, we have
to refine all the transits in order to make the set of distinguished
triangulations $(T,H)$ for $(W,L)$ stable with respect to them. First, a negative
$3\to 2$ or $2 \to 0$ move is admissible if and only if the
disappearing edge does not belong to the link. Moreover (see
\cite[Section 4.1.1]{BB2} for details):

We say that a bubble move $(T,H) \rightarrow (T',H')$ on a $2$--face
$t$ of $T$ is {\it distinguished} if $t$ contains an edge $e$ of $H$,
and $H'$ is defined by replacing $e$ in $H$ with the two new edges of
$T'$ making with $e$ the boundary of a $2$--face. We have a {\it bubble
  charge transit} $(T,H,c) \rightarrow (T',H',c')$ if the sum of
charges: stays equal about the two edges of $t$ distinct from $e$;
goes from $0$ to $2$ about $e$; is equal to $0$ about the two new
edges of $H'$; is equal to $2$ for the remaining edge of the two new
tetrahedra of $T'$.

Such bubble transits preserve the Hamiltonian property of
the link realization. The flat/charged $\Ii$--transits are just
obtained by assembling the above definitions.

\medskip

{\bf $\Dd$--triangulations of pairs $(W,\rho)$}\qua This
notion, as well the {\it idealization} procedure, was introduced in
Subsection \ref{ITdefsec}. Recall that quasi-regular triangulations of
$W$ support idealizable $\Dd$--triangulations for any pair $(W,\rho)$
(in fact {\it generic} 1-cocycles on a quasi-regular triangulation are
idealizable). Moreover, there is the following natural notion of:

\medskip

{\bf $\Dd$--transits}\qua Let $(T_0,b_0)\to (T_1,b_1)$ be a
transit of branched triangulations of $W$ and $z_k \in
Z^1(T_k;PSL(2,\C))$, $k=0,1$. We have a {\it cocycle transit}
$(T_0,z_0) \leftrightarrow (T_1,z_1)$ if $z_0$ and $z_1$ agree on the
common edges of $T_0$ and $T_1$. This makes an {\it idealizable
cocycle transit} if both $z_0$ and $z_1$ are idealizable $1$--cocycles,
and in this case we say that $(T_0,b_0,z_0) \leftrightarrow
(T_1,b_1,z_1)$ is a $\Dd$--transit. It is easy to see that $z_0$ and
$z_1$ as above represent the same flat bundle $\rho$. We have:

\begin{prop}\label{DdomI} {\rm (\cite{BB2}, Proposition 2.16)}
  Consider a fixed pair $(W,\rho)$, and denote by $\Ii$ the
  idealization map $\Tt \to \Tt_{\Ii}$ on the $\Dd$--triangulations of
  $(W,\rho)$. For any $\Dd$--transit $\mathfrak{d}$ there exists an
  $\Ii$--transit $\mathfrak{i}$ (respectively for any $\mathfrak{i}$ there
  exists $\mathfrak{d}$) such that $\mathfrak{i} \circ \Ii = \Ii \circ
  \mathfrak{d}$.
\end{prop}

Similarly, the natural behaviour of $\Dd$--triangulations
with respect to branching changes dominates, via the idealization, the
one of $\Ii$--triangulations.

Globally flat/charged $\Ii$--triangulations for either triples
$(W,L,\rho)$ or gentle cusped manifolds $M$ were defined in
\ref{ITdefsec}. Let us discuss now {\it arbitrary} cusped manifolds.

\medskip

{\bf General cusped manifolds}\qua  Let $M$ be an arbitrary
cusped manifold.  Given a quasi-geometric geodesic triangulation
$(T,w)$ of $M$, there exists a finite sequence $T \rightarrow \ldots
\rightarrow T'$ of positive $2\to 3$ moves such that $T'$ supports a
global branching $b'$ \cite[Theorem 3.4.9]{BP2}. For each move, we can
define the transit of cross-ratio moduli by (\ref{ideqmod}),
with $*=1$ everywhere.

\begin{defi} \label{W-T} {\rm We say that $M$ is {\it weakly-gentle}
    if there exists such a sequence $T \rightarrow \ldots \rightarrow
    T'$ of positive $2\to 3$ moves that lifts to a sequence of
    transits $(T,w) \to \ldots \rightarrow (T',w')$. In that case, we
    call $(T',b', (w')^{*_b})$ {\it an $\Ii$--triangulation} of $M$.
    For any flattening $f'$ for $(w')^{*_b}$ and integral charge $c'$,
    we say that $(T',b', (w')^{*_b},f',c')$ is a flat/charged
    $\Ii$--triangulation of $M$.}
\end{defi}

In this definition, the exponent $*_b$ means that the cross-ratio
moduli of a tetrahedron are turned to the inverse if its branching
orientation is negative. The authors do not know any example of non
weakly-gentle cusped manifolds. Roughly speaking, $M$ is not
weakly-gentle if, in order to give any quasi-geometric geodesic
triangulation of $M$ a global branching, we are forced to introduce
new interior vertices by performing some bubble moves. Dealing with
the {\it quantum} state sums, this gives rise to a technical
difficulty similar to the one that leads to the ``link fixing'' for
closed manifold $W$. This motivates the following definition.
\begin{defi} \label{N-W-T} {\rm If $M$ is a cusped manifold that is not
    weakly-gentle, we call a {\it marking} of $M$ the choice of an
    edge $l$ in the canonical Epstein--Penner cell decomposition of
    $M$. A {\it flat/charged $\Ii$--triangulation
    $(T',H',b',(w')^{*_b},f',c')$ of $(M,l)$} consists of: a
    triangulation $T'$ of the polyhedron $Q(M)$, obtained as in
    Definition \ref{W-T} from a quasi-geometric geodesic triangulation
    of $M$ via positive $2\to 3$ moves {\it and bubble moves}; a
    Hamiltonian subcomplex $H'$ of the $1$--skeleton of $T'$ with one
    edge $l_1$ isotopic to $l$, and the union of the other edges
    isotopic to a second copy $l_2$ of $l$ (we summarize this property
    by saying that $(T',H')$ is a {\it distinguished} triangulation of
    $(M,l)$); a flattening $f'$ for $(w')^{*_b}$; a collection $c'$
    of integral charges on the abstract tetrahedra of $T'$, such that
    the sum of the charges is equal to $4$ about $l_1$, $0$ about each
    edge of $l_2$, and $2$ about the other edges of $T'$.

If $M$ is weakly-gentle we set $l=\emptyset$,
so that the present definition incorporates Definition \ref{W-T}.
}
\end{defi}

From now on, an $\Ii$--triangulation of a cusped manifold $M$ always
mean one as in Definition \ref{N-W-T}, possibly forgetting the marking
and the flat/charge.

Let us say that a $\mz$--valued decoration of the
abstract edges of an $\Ii$--triangul\-ation is a {\it rough} integral
charge, flattening, or flat/charge if it is {\it locally} of that
form, on each tetrahedron.  Recall that a flat/charged
$\Ii$--tetrahedron $(\Delta,b,w,f,c)$ defines an enriched
$\Ii$--tetrahedron $(\Delta,b,w,a)$ as in Section \ref{LNbdilog} by
putting $a_j=f_j-*_bc_j$, $j=1,2,3$.

\begin{lem} \label{globalrough} Any rough flat/charge whose associated
  $N$th--branch map satisfies the
  relations (\ref{ideqmodprime}) comes from a pair $(f,c)$, where $f$
  is a {\it global} flattening and $c$ a global integral charge.
\end{lem} 

{\bf Proof}\qua For any rough flattening $f$ and integral
  charge $c$, any other is of the form $f'=f+b$ or $c'=c+d$, where $b$
  satisfies $b_0+b_1+b_2=0$, and similarly for $d$. So any rough
  flat/charge $a$ locally appears as $a=f-*_bc=f'-*_bc''$ for a
  suitable rough integral charge $c''$. In particular we can assume
  that $f$ is a global flattening. But imposing now the condition
  (\ref{ideqmodprime}) for all $N>1$ simultaneously, we realize that
  $c$ is necessarily a global integral charge.\endproof

{\bf Cohomological normalization of flattenings and charges}\qua
The reductions mod$(2)$ of both global flattenings and integral
charges induce cohomology classes in $H^1(W;\Z/2\Z)$ or
$H^1(M;\Z/2\Z)$.  Moreover, in the case of a cusped manifold $M$, both
the log-branches and integral charges induce classes in
$H^1(\partial M; \Z)$, where $\partial M$ denotes the toric boundary
components of the natural compactification of $M$. These classes are
defined as follows. For any mod($2$) (respectively integral) $1$--homology
class $a$ in $W$ or $M$ (respectively $b$ in $\partial M$), realize it by a
disjoint union of (respectively oriented and essential) closed paths
transverse to the triangulation and `without back-tracking', ie, such
that they never depart from a $2$--face of a tetrahedron
(respectively $1$--face of a triangle) by which they entered. Then the
mod($2$) sum of the flattenings or charges of the edges that we
encounter when following such paths in $W$ or $M$ define the value of
the corresponding class on $a$. Similarly, the {\it signed} sum of the
log-branches or charges of the edges whose ends are vertices that we
encounter when following such paths on $\partial M$ define the value
of the corresponding class on $b$; for each vertex $v$ the sign is
$*_b$ (respectively $-*_b$) if, with respect to $v$, the path goes in the
direction (respectively opposite to the one) given by the orientation of
$\partial M$. In Lemma 4.12 of \cite{BB2} we proved that the transits
of global integral charges for pairs $(W,L)$ preserve these
cohomology classes. This result extends immediately to cusped
manifolds $M$, to flattenings for both $(W,\rho)$ or $M$, and also to
log-branches for $\partial M$. So, from now on, we normalize the
global flattenings and integral charges by requiring that {\it all
these classes are trivial}. (Otherwise, the dilogarithmic invariants
would just depend on them).

 \begin{teo}\label{existencef/c}$\phantom{99}$ 
   
{\rm(1)}\qua Every pair $(W,\rho)$ or cusped manifold $M$ has
   $\Ii$--triangulations $\Tt_{\Ii}$, and every such triangulation
   admits global flattenings $f$.

{\rm(2)}\qua For every triple $(W,L,\rho)$ or pair $(M,l)$ there exist
distinguished $\Ii$--triang\-ulations, and every such triangulation
admits global integral charges $c$.
\end{teo}

The first claim in (1) is essentially proved in Subsection
\ref{ITdefsec}. For $(W,\rho)$ take any quasi-regular triangulation
with an idealizable cocycle representing $\rho$. For non-weakly-gentle
cusped manifolds, it is enough to use generic bubble transits in the
sequence of moves considered in Definition \ref{N-W-T} (see the proof
of Theorem \ref{connessione}). The first claim in (2) for triples
$(W,L,\rho)$ is a result of \cite{BB2}. The existence of global
integral charges for $(W,L)$ is shown in Chapter 2 of \cite{B}.  For
pairs $(M,l)$, see the proof of Theorem \ref{connessione}. The general
existence of global flattenings and integral charges are slight
adaptations of earlier results of W Neumann (first claim of Theorem 4.2
in \cite{N2}, and Theorem  2.4 $i)$ in \cite{N3} respectively).

Let $(\Tt_{\Ii},f)$ be a flattened $\Ii$--triangulation of $(W,\rho)$
or $M$. Let $(\Tt_{\Ii},f,c)$ be a flat/charged $\Ii$--triangulation of
$(W,L,\rho)$ or $(M,l)$. By using the
symmetrized classical dilogarithm $\Rr_1$, we can define the state sum
$$\textstyle
\Rr_1(\Tt_{\Ii},f) = \prod_{i} \Rr_1(\Delta^i,b^i,w^i,f^i)$$ which
reads as the exponential of $1/i\pi$ times a signed sum of uniformized
Rogers dilogarithms. By using the quantum matrix
dilogarithms $\Rr_N$, $N>1$, we can define the state sums 
$$\textstyle
\Rr_N(\Tt_{\Ii},f,c) = \sum_\alpha \prod_i
\Ll(\Delta^i,b^i,w^i,f^i,c^i)_\alpha$$
\noindent where the sum is over the states of $\Tt_{\Ii}$.  Let us
denote by $v$ the number of `interior' vertices of a given
triangulation ($v=0$ in the case of a weakly-gentle cusped manifold).
Put $X= (W,\rho)$ or $M$, and $\widetilde{X}= (W,L,\rho)$ or $(M,l)$
(recall that $l=\emptyset$ if $M$ is weakly-gentle). The main result
of this section is:

\begin{teo}\label{diloginv} $\phantom{99}$

  {\rm(1)}\qua The value of $\Rr_1(\Tt_{\Ii},f)$ does not depend, up to sign,
  on the choice of the flattened $\Ii$--triangulation $(\Tt_{\Ii},f)$.
  Hence, up to sign, $H_1(X) := \Rr_1(\Tt_{\Ii},f)$ is a well defined
  invariant called the {\rm classical dilogarithmic invariant} of $X$.

  {\rm(2)}\qua For every odd $N>1$, the value of $N^{-v}
  \Rr_N(\Tt_{\Ii},f,c)$ does not depend, up to sign and multiplication
  by $N$th roots of unity, on the choice of the flat/charged
  $\Ii$--triangulation $(\Tt_{\Ii},f,c)$.  Hence $H_N(\widetilde{X}) :=
  N^{-v}\Rr_N(\Tt_{\Ii},f,c)$ is a well defined invariant, up to the
  same ambiguity, called a {\rm quantum dilogarithmic invariant} of
  $\widetilde{X}$.
\end{teo}

We note that the normalization $N^{-v}$ of the quantum state sums
comes from the behaviour of $\Rr_N$ for bubble flat/charged
$\Ii$--transits.

The first step is to extend the proof of the transit invariance of the
state sums to the other transits, besides the $2\to 3$ one.

\begin{lem} \label{consmove} For any odd $N\geq 1$ and any flattened
  (respectively flat/charged) $\Ii$--transit between flattened
  (respectively flat/charged) $\Ii$--triangulations of $X$
  (respectively $\widetilde{X}$), the traces of the two patterns of
  associated matrix dilogarithms of rank $N$ are equal, possibly up to
  a sign and an $N$th root of unity phase factor.
\end{lem}

{\bf Proof}\qua Consider an $\Ii$--triangulation of $X$, and a
flattening transit that lifts the $\Ii$--transit of Figure
\ref{CQDidealt}. As usual, denote by $\Delta^j$ the tetrahedron
opposite to the $j$-th vertex for the ordering of the vertices induced
by the branching. Do a further $2 \rightarrow 3$ flattening transit
with $\Delta^0$ and $\Delta^2$. A mirror image of $\Delta^4$ appears,
which together with $\Delta^4$ forms the final configuration of a $0
\rightarrow 2$ flattening transit. The other two new flattened
$\Ii$--tetrahedra have the same decorations and gluings as $\Delta^1$
and $\Delta^3$. Hence Theorem \ref{Rtransit} implies that $0
\leftrightarrow 2$ transits correspond to identities of the form
(above for $\Delta^4$):
$$\Rr_1(\Delta,b,w,f) \ \Rr_1^{-1}(\bar{\Delta},b,w,f) = 0 .$$
Here $\bar{}$ denotes the opposite orientation; as remarked before Proposition
\ref{DdomI}, the mirror moduli and flattenings are the same. Finally,
we note that the final configuration of a bubble move is just obtained
by gluing two more faces in the final configuration of a $0
\rightarrow 2$ move, and the moduli and flattenings behave well with
respect to this gluing (that is the relations (\ref{ideq2mod}) and
(\ref{ideq2}) are satisfied). So we get the very same relation for a
bubble flattening transit. This proves the statement for $\Rr_1$.

We prove the statement for $\Rr_N$, $N>1$, and the $0
\leftrightarrow 2$ and bubble flat/charged $\Ii$--transits via the very
same arguments. The corresponding two term relations are:
$$\begin{array}{c} \Rr_N(\Delta,b,w,f,c)\ \Rr_N(\bar{\Delta},b,w,f,c)
\equiv_N {\rm Id}_{\mc^N} \otimes {\rm Id}_{\mc^N} \\ \\ {\rm
Trace}_i\bigl(\Rr_N(\Delta,b,w,f,c)\ \Rr_N(\bar{\Delta},b,w,f,c)\bigr)
\equiv_N N\ \cdot\ {\rm Id}_{\mc^N}
\end{array}$$
where $\bar{}$ is as above, and the trace is over one of the tensor
factors in the first relation. The trace appears for the bubble
transits because, as we said above, the final configuration of the
underlying moves are obtained from $0 \leftrightarrow 2$ ones by
gluing two more faces. Again, we use the fact that the charges behave
well under this gluing. \endproof

As the state sums are fully transit invariant, the main theorem
follows from the following triangulation connectedness theorem.

\begin{teo}\label{connessione} $\phantom{99}$

{\rm(1)}\qua Any two flattened $\Ii$--triangulations of $X$ can be connected via a
    finite number of transits of flattened $\Ii$--triangulations.

{\rm(2)}\qua Any two flat/charged $\Ii$--triangulations of $\widetilde{X}$ 
can be connected via a finite number of transits of flat/charged 
$\Ii$--triangulations.
\end{teo}
\begin{remarks}\label{improve}
{\rm (1)\qua A version of (1) was proved independently by Neumann in \cite{N4}. 

  (2)\qua For triples $(W,L,\rho)$, this theorem is a genuine refinement
  of the main result of \cite{BB2}. In that paper, we gave a complete
  proof of a weaker connectedness result, enough to show that the QHI
  are well defined, but not to get the invariance of the scissors
  congruence classes discussed in Section \ref{CQDSCISSORS}. The
  present strong version was mentioned in Section 4.5 of \cite{BB2}
  without proof. Along with the proof of Theorem \ref{connessione} we
  shall stress the main differences with respect to \cite{BB2}.}
\end{remarks}

In the course of the proof we shall use the following two
results. The first ensures the connectedness of the space of {\it
branched} (topological) ideal triangulations of an arbitrary compact
oriented $3$--manifold, and the second is an uniqueness (rigidity)
result for $\Ii$--triangulations of cusped manifolds with maximal
volume.

\begin{teo} \label{bconnessione} {\rm \cite{C}}\qua 
For any two branched triangulations $T$ and $T'$ of a same compact
  oriented $3$--manifold $Y$, there exists a finite sequence $T
  \rightarrow \ldots \rightarrow T'$ made of $2 \to 3$, $0 \to 2$,
  bubble moves or their inverses.
\end{teo}

\begin{teo}\label{unique} {\rm (See \cite{Dum} or \cite{F2})}\qua Let
  $(T,b,w)$ be an arbitrary $\Ii$--triangulation of a
  cusped $3$--manifold $M$. Then, among all $\Ii$--triangulations
  supported by $(T,b)$, this is the only one such that the algebraic
  sum of the volumes of its $\Ii$--tetrahedra equals Vol$(M)$.
\end{teo}

{\bf Proof of Theorem \ref{connessione}}\qua  In Proposition 4.27
of \cite{BB2} we proved for triples $(W,L,\rho)$ a version of Theorem
\ref{connessione} {\it up to branching changes}; we included
in the definition of a $\Dd$--triangulation that it was
{\it quasi-regular}, and we used {\it only} $\Ii$--triangulations
covered by such $\Dd$--triangulations, via the idealization. Moreover,
the flattenings were hidden (see Remark \ref{compareQHI}). We note
that for pairs $(W,\rho)$ the proof is easier because there is no
link.

The proof of Theorem \ref{connessione} up to branching changes
is enough for a weaker version of Theorem \ref{diloginv} $(1)$, where
the invariants are defined up to a sign and multiplication by $12$-th
roots of unity when $N=1$. We have only to apply Proposition
\ref{Rsym3}. On the other hand, to get the statement with the weaker
ambiguity we must show that even the branching changes can be realized
by means of a sequence of transits, and then exploit the global
compensations occuring when doing them (see the proof of Theorem
\ref{Rtransit}). For that, we use Theorem \ref{bconnessione}. 

The use of $\Ii$--triangulations dominated, via the idealization, by
quasi-regular $\Dd$--triangulations makes the proof simpler because
these triangulations support {\it generic} idealizable
$PSL(2,\mc)$--valued $1$--cocycles, and we have room for choosing paths
of transits between their idealizations. In \cite{BB2} we proved first that
quasi-regular $\Dd$--triangulations can be connected by (quasi-regular)
$\Dd$--transits, and we applied Lemma \ref{DdomI}. Here we want to do
it for $\Ii$--triangulations of a triple $(W,L,\rho)$ dominated by non
necessarily quasi-regular $\Dd$--triangulations, as it can actually
happen when $\rho$ is non trivial, and develop arguments that
eventually apply also to $\Ii$--triangulations of cusped manifolds,
that are not dominated by any $\Dd$--triangulation. The key point is to
show that we can get round of accidental stops to
the $\Ii$--transits. We do it as follows.

\begin{figure}[ht]
\begin{center}{\tiny
\psfraga <-2pt,0pt> {B}{$B$}
\psfraga <-2pt,0pt> {R}{$R$}
\psfraga <-2pt,0pt> {R'}{$R'$}
\psfraga <-2pt,0pt> {R''}{$R''$}
\psfraga <-2pt,0pt> {m}{$m$}
\includegraphics[width=10cm]{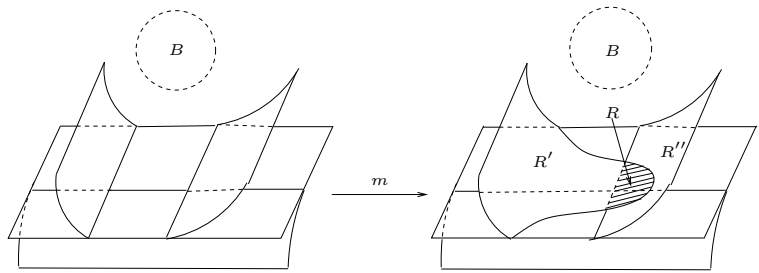}}
\caption{\label{GTbadmove1} A bad move} 
\end{center}
\end{figure}

Suppose we have a sequence of branched moves
between two $\Ii$--triangulations supported by $T$ and $T'$, such that
some of them do not allow an $\Ii$--transit (for instance, when $x=y$
in Figure \ref{CQDidealt}); call them {\it bad moves}. Take the first,
$m$. Dually (at the level of standard spines), we are in a situation
like in Figure  \ref{GTbadmove1}, where $B$ is a $3$--cell with the two
regions $R'$ and $R''$ lying on its boundary (immersed)
$2$--sphere. The shaded region $R$ is ``bad'', because one among the
cross-ratio moduli attached to its corners belongs to the forbidden
values $\{0,1\}$ (recall that the vertices are dual to
tetrahedra). Just before $m$, do a (positive) bubble $\Ii$--transit by
gluing a disk in the {\it interior} of $B$, which we call a {\it
capping disk}. For distinguished triangulations of triples
$(W,L,\rho)$, this bubble move takes place near the link $H$, and is
as explained before Proposition \ref{DdomI}. In particular, one,
$R_H$, of the two regions on which we glue the capping disk is dual to
an edge of $H$, and after the move, the capping disk as well as the
new region locally ``opposite'' to it and adjacent to $R_H$ are dual
to edges of $H$. (Recall that $H$ is Hamiltonian, so that $R_H$
exists).

\begin{figure}[ht]
\begin{center}{\tiny
\psfraga <-2pt,0pt> {R1}{$R1$}
\psfraga <-2pt,0pt> {R2}{$R2$}
\psfraga <-2pt,0pt> {R3}{$R3$}
\psfraga <-2pt,0pt> {R4}{$R4$}
\psfraga <-2pt,0pt> {R5}{$R5$}
\psfraga <-2pt,0pt> {m'}{$m'$}
 \includegraphics[width=10cm]{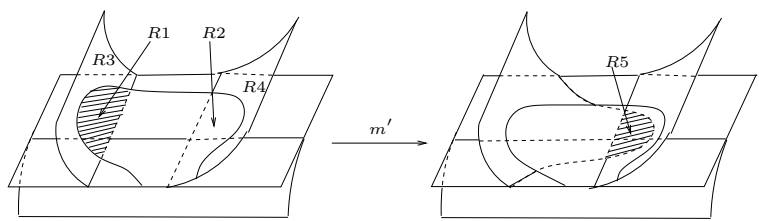}}
\caption{\label{GTbadmove2} The capping disk of a bubble move turns
  the bad move into a good one} 
\end{center}
\end{figure}

A key point is that there is one degree of freedom in choosing the
bubble $\Ii$--transit. So, for a generic choice of it we can slide via
successive $2 \to 3$ $\Ii$--transits a portion of the capping disk in
the position shown at the left of Figure  \ref{GTbadmove2}, and also
achieve the move $m'$ as an $\Ii$--transit (we have to keep track of
the region $R_H$, for otherwise the link $H$ would split). Of course,
now the cross-ratio moduli attached to the tetrahedra dual to the
vertices of that configuration, ie, at all the corners of the regions
$R_1,\ldots,R_5$, are altered compared to those already present in
Figure \ref{GTbadmove1}. By continuing the initial sequence between the
spines $P(T)$ and $P(T')$, the regions involved in the (dual) moves
following $m$ can intersect the capping disk. However, the arguments
of Proposition 4.23 in \cite{BB2} show that the capping disk is not an
obstruction for doing these subsequent moves, so that we eventually
reach the very same position as in $P(T')$ by sliding ``under'' the
capping disk (this is possible essentially because the moves are
purely local).  Hence, by applying the same procedure each time we
meet a bad move and using the genericity of bubble $\Ii$--transits, we
get a sequence of $\Ii$--transits from $T$ to a triangulation $T''$. We
can see that the standard spine $P(T'')$ dual to $T''$ is obtained
from $P(T)$ just by gluing successively some (non adjacent) $2$--disks,
each one corresponding to a capping disk.  Again by using the
arguments of Proposition 4.23 in \cite{BB2}, we can remove them one
after the other via further $\Ii$--transits, and eventually get a
sequence of $\Ii$--transits from $T$ to $T'$.

If the so obtained $\Ii$--triangulation supported by $T'$ is not the
one we started with, note that the very same arguments imply that we
can find sequences of $\Ii$--transits from $T'$ to a quasi-regular
triangulation $T_0$, thus supporting two distinct
$\Ii$--decorations. As there exist $\Dd$--triangulations supported by
$T_0$ that dominate each of these $\Ii$--decorations, we can apply the
proof of Proposition 4.27 in \cite{BB2} to see that they can be
connected via $\Ii$--transits.

Finally, the argument for the invariance with respect to the choice of
integral charge given in Proposition 4.27 of \cite{BB2} works
word-by-word also for the flattenings, because flattenings and
integral charges on the same triangulation are both affine spaces over
the {\it same} integral lattice. This follows directly from their
construction, see the references after Theorem \ref{existencef/c}.
This concludes the proof for the case of closed manifolds.

For ideal triangulations of cusped manifolds $M$, the same
arguments allow us to prove the following preliminary result:

\medskip

\noindent {\sl Let $\Tt_{\Ii}=(T,b,w)$ and $\Tt'_{\Ii}=(T',b',w')$ be
  two $\Ii$--triangulations of $M$ (quasi geometric if $M$ is gentle,
  as in Definition \ref{N-W-T} otherwise).
  Then there exists a finite sequence of $\Ii$--transits connecting
  $\Tt_{\Ii}$ to $\Tt''_{\Ii}= (T',b',w'')$, that is another
  $\Ii$--triangulation supported by $(T',b')$.}

\medskip

We stress that we possibly include also bubble transits, adding
interior vertices. The volume of $\Tt_{\Ii}$, which is the algebraic
sum of the volumes of its hyperbolic ideal tetrahedra, is not altered
by $\Ii$--transits and the transits used for defining
$\Ii$--triangulations of general (non gentle) cusped manifolds. So the
volume of $\Tt_{\Ii}$, $\Tt'_{\Ii}$ and $\Tt''_{\Ii}$ is that of the
hyperbolic manifold $M$, and Theorem \ref{unique} yields $w'=w''$.

We can include flattenings to the transits. However, the sum of the
charges about one edge of the initial triangle involved in a bubble
move goes from $2$ to $4$ (there was no link inside $M$), whereas in
the final configuration this sum is equal to $0$ for two of the new
edges, and is $2$ for all the other edges. Topologically, this means
choosing two of the new edges to form an unknotted and properly
embedded arc in the natural compactification of $M$ (with one boundary
torus at each end) that passes through the new ``interior'' vertex.
Hence, starting from $\Tt_{\Ii}$, when we use a charge transit we
obtain a rough integral charge which is not a global one as defined in
Section \ref{ITdefsec}.

If $M$ is weakly-gentle, as $\mathcal{T}_{\Ii}'$ does not contain
interior (non singular) vertices, any sequence of flat/charged
$\Ii$--transits connecting $\Tt_{\Ii}$ to $\mathcal{T}_{\Ii}'$
eventually gives a true global integral charge.

If $M$ is not weakly-gentle, we consider only distinguished
$\Ii$--triangulations $(\Tt_{\Ii},H)$ and $(\mathcal{T}_{\Ii}',H')$ of
$(M,l)$, for a marking $l$ of $M$. We do bubble charge transits as
usual, on $2$--simplices that contain an edge of $l_2$. Any sequence of
flattening $\Ii$--transits connecting $(\Tt_{\Ii},f)$ to
$(\mathcal{T}_{\Ii}',f')$ can be lifted to a sequence
$(\Tt_{\Ii},H,f,c) \rightarrow \ldots \rightarrow
(\mathcal{T}_{\Ii}',H',f'',c'')$ between distinguished flat/charged
$\Ii$--triangulations of $(M,l)$.

It remains to show the invariance with respect to the choice of
flattening and integral charge on a fixed $\Ii$--triangulation of $M$.
We do this as for closed manifolds.  Namely, as already said above,
the sets of flattenings and integral charges on a given ideal
triangulation are both affine spaces over the same integral lattice,
and the generators of this lattice depend only the local combinatorics
of the triangulation (each generator is associated to the abstract
star of an edge).\endproof

\begin{remark}\label{nodis}{\rm (No discrepancy for weakly-gentle
      cusped manifolds)\qua The flattenings of $\Ii$--triangulations of
    weakly-gentle cusped manifolds are just $-*_b$ times the integral
    charges. Thus, we can take $f-*_bc=-2*_bc$ for the flat/charges of
    the matrix dilogarithms in $\Rr_N(\Tt_{\Ii},f,c)$, $N >1$. Hence,
    the classical and quantum dilogarithmic invariants of
    weakly-gentle cusped manifolds are defined on the very same
    geometric supports, the flattened $\Ii$--triangulations.}
\end{remark} 

\begin{remark}
{\rm (Common  features of classical and quantum invariants)\qua 
We defined $\Rr_1(\Tt_{\Ii},f)$ by taking the {\it exponential of
$(1/i\pi)$ times} a sum of uniformized Rogers
dilogarithms. In this way, for every $N\geq 1$, we have:

(1)\qua $H_1(X)$ and $H_N(\widetilde{X})^N$, $N>1$, are well defined up to
    sign ambiguity;

(2)\qua $H_1(X)^* = H_1(X^*)$ and $H_N(\widetilde{X})^* =
H_N(\widetilde{X}^*)$, $N>1$, where $^*$ denotes the complex
conjugation, and $X^*$ (respectively $\widetilde{X}^*$) denotes $(-W,\rho^*)$
or $-M$ (respectively $(-W,L,\rho^*)$ or $(-M,l)$). For $N >1$ this is a consequence of
Proposition 4.29 in \cite{BB2}. For $N=1$ this follows easily from the
behaviour of the Rogers dilogarithms under complex conjugation of its
argument.}
\end{remark} 

\begin{remark}{\rm (On the phase ambiguity of $H_N(\tilde{X})$,
      $N >1$)\qua There is only a sign ambiguity for the flat/charged
      $\Ii$--transit corresponding to Figure \ref{CQDidealt}. So we could
      expect that the ambiguity up to $N$th roots of unity do vanish
      for {\it all} the transits in Corollary \ref{qtransit}, due to
      global compensations (see Theorem \ref{Rtransit} and Remark
      \ref{commut}). However, although this happens for a
      certain non trivial subset of transits, in general this is not
      the case.

      \noindent Indeed, by using Lemma \ref{idemQHI}, a statement of
      Proposition \ref{6jsym} with the full $N$th root of unity
      dependence follows from Proposition 6.4 in the Appendix of
      \cite{BB2}, for a specific choice of flattenings of the involved
      moduli. The symmetry defects appearing there are formally the
      exponentials of terms of the {\it very same} form as those
      coming from Lemma \ref{Rsym3}, so that we get a table as in
      Section \ref{cftrel}, where the defects are replaced with powers
      of $\zeta$ depending linearly on the charges. But we cannot
      deduce the global compensations, because, as
      $\zeta=\exp(2i\pi/N)$ and $N >1$, we cannot work with the
      charges mod($2$). At the time of this writing the authors do not
      see any way to renormalize $\Rr_N$ to avoid this discrepancy.
      Understanding the {\it geometric meaning} of the phase ambiguity
      is an important open question about the quantum dilogarithmic
      invariants.}
\end{remark}

\subsection{$H_1(X)$ and Cheeger--Chern--Simons classes}\label{CSV} 

For every pair $(W,\rho)$ set
$${\bf R}(W,\rho) := {\rm CS}(\rho) + i\ {\rm Vol}(\rho) \in
\mc/(\pi^2\Z)$$ where ${\rm CS}(\rho)$ and ${\rm Vol}(\rho)$ are
respectively the Chern-Simons invariant and the volume of $\rho$. As
we consider $\rho$ up to isomorphisms of flat bundles, it may be
identified with the conjugacy class of its holonomy representations.
We refer to Chapter 10-11 of \cite{D2} and the references therein for
details on these notions.

Meyerhoff extended in \cite{Mey} the definition of CS to cusped
manifolds $M$, so that we can consider again
$${\bf R}(M) := {\rm CS}(M) + i\ {\rm Vol}(M) \in \mc/(\pi^2\Z).$$
The following result holds (we use the notation of Section 
\ref{CQDSYMROGERS}):

\begin{teo}\label{RNChS} 
Let $X$ be either a pair $(W,\rho)$ or a cusped manifold $M$, and
$(\Tt,f)$ be a flattened $\Ii$--triangulation of $X$. Then 
$${\bf R}(X)=
\sum_{\Delta \subset T} *\ {\rm R}(\Delta,b,w,f) 
\quad
 {\rm mod}(\pi^2\mz).$$
Hence $H_1(X)$ is equal to $\exp\bigl((1/i\pi){\bf R}(X)\bigr)$, where
both invariants are defined up to a sign.
\end{teo}

This theorem is proved in \cite{N4}, using earlier deep results
of Dupont--Sah \cite{DS}, and Dupont \cite{D1}.

Our proof that
$$H_1(X)= \exp \biggl((1/i\pi) \sum_{\Delta \subset T} *\ {\rm
  R}(\Delta,b,w,f) \biggr)$$ is a well defined invariant of $X$ is
  independent of its identification with $\exp\bigl((1/i\pi){\bf
  R}(X)\bigr)$. This proof shows that $\textstyle \sum_{\Delta \subset
  T} *\ {\rm R}(\Delta,b,w,f)$ mod($\pi^2\mz$) itself does not depend
  on the choice of $(\Tt,f)$. It is based on direct geometric
  manipulations of decorated triangulations, and is structurally the
  same for the whole family of classical and quantum dilogarithmic
  invariants. On the other hand, in the classical case one can adopt
  the slightly different point of view of looking for {\it simplicial
  formulas} for the already known classical invariant ${\bf R}(X)$
  (thus obtaining, by the way, that the values of these formulas do
  not depend on the choice of the combinatorial support).

\section{Scissors congruence classes}\label{CQDSCISSORS}

We construct further invariants of $3$--manifolds called {\it scissors
  congruence class\-es}, that belong to suitably defined {\it
  (pre)-Bloch-like groups}. Later we discuss some problems about the
relations between these invariants and the dilogarithmic ones.

Fix one {\it base} oriented tetrahedron $\Delta$.  By varying the
respective decorations, we get the sets of $\Dd$--tetrahedra
$\{*_b(\Delta,b,z)\}$, $\Ii$--tetrahedra $\{*_b(\Delta,b,w)\}$,\break
 flattened $\Ii$--tetrahedra
$\{*_b(\Delta,b,w,f)\}$ and flat/charged $\Ii$--tetrahedra\break
$\{*_b(\Delta,b,w,f,c)\}$. Let us call
them generically $\Aa$--{\it tetrahedra}. We will specify $\Aa = \Dd,
\Ii, \Ii_f, \Ii_{fc}$ when necessary. We also denote by $\Aa$ the set
of all $\Aa$--tetrahedra and by $\Z[\Aa]$ the free $\Z$--module
generated by $\Aa$. We stipulate that the sign of $\Aa$--tetrahedra
is compatible with the algebraic sum in $\Z[\Aa]$, ie,
$-(-(\Delta,b,z))=(\Delta,b,z)$, and so on.

Any instance of $\Aa$--transit naturally induces a linear relation
between the involved $\Aa$--tetrahedra. Consider the relations
on $\Z[\Aa]$ generated by the five term identities corresponding to all
instances of $2 \leftrightarrow 3$ $\Aa$--transit. Denote by
$\Pp(\Aa)$ the resulting quotient of $\Z[\Aa]$. We call $\Pp(\Aa)$ the
$\Aa$--{\it (pre)-Bloch-like group}.

In the previous sections we described the behaviour of
$\Aa$--tetrahedra with respect to the tetrahedral symmetries. The
identification of $\Aa$--tetrahedra related by these symmetries gives
new relations in $\Z[\Aa]$, whence a quotient map $\Pp(\Aa)\to
\Pp'(\Aa)$.  Clearly there are forgetful maps $\Pp(\Ii_{fc})\to
\Pp(\Ii_f )\to \Pp(\Ii)$, and similar ones for the $\Pp'(\Aa)$,
making a commutative diagram with the quotient maps $\Pp(\Aa)\to
\Pp'(\Aa)$. If $\Tt$ is any $\Aa$--triangulation of any oriented
$3$--manifold $Y$, then the formal sum of the $\Aa$--tetrahedra of $\Tt$
determines an element $[\Tt]_\Aa \in \Pp(\Aa)$
(respectively $[\Tt]'_\Aa \in \Pp'(\Aa)$).

Some results of the previous section can be immediately
rephrased and somewhat illuminated in the present set up. Noting
that every $\Ii$--tetrahedron is the idealization of a
$\Dd$--tetrahedron, the fact that $\Dd$--transits dominate $\Ii$--transits
yields:
\begin{prop}\label{idealmor} 
The idealization induces a surjective homomorphism $\Ii_\Pp\colon$
$\Pp(\Dd)\to \Pp(\Ii).$
\end{prop}
Also, the results of Section \ref{CQDSYMROGERS} give the following.
\begin{prop}\label{blochfun} 
For any element $\textstyle \sum_i a_i*_{b_i}(\Delta,b_i,w_i,f_i)$ in
$\Pp(\Ii_{f} )$ (respectively $\Pp'(\Ii_{f})$), the formula
$$\sum_i a_i*_{b_i}{\rm R}(\Delta,b_i,w_i,f_i)$$ defines a function
${\rm R}\co  \Pp(\Ii_{f} )\to \C/\pi^2 \Z$ (respectively
${\rm R}'\co  \Pp'(\Ii_{f})\to \C/(\pi^2/6) \Z$).
\end{prop}
Obviously we can also define the function $H_1 = \exp({\rm R}/i\pi)$
on $\Pp(\Ii_f)$.
\begin{remark}\label{veroquot}{\rm We know from Section
    \ref{CQDSYMROGERS} that the function ${\rm R}'$ is well defined
    {\it only} mod$(\pi^2/6) \Z$, while the weaker ambiguity in the
    definition of ${\rm R}$ is due to remarkable global compensations
    occurring in the $2 \leftrightarrow 3$ flattened $\Ii$--transits.
    This means, in particular, that $\Pp'(\Ii_{f})$ is a {\it genuine}
    quotient of $\Pp(\Ii_{f})$.  Recall, on the contrary, that for the
    classical pre-Bloch group the tetrahedral symmetries are
    consequences of the five-term identity \cite{DS}.}
\end{remark}

With the usual notation, denote $X=(W,\rho)$ or $M$, and
$\widetilde{X}=(W,L,\rho)$ or $(M,l)$. Let $(\Tt,f)$ be a
flattened $\Ii$--triangulation for $X$, with the
associated element $[(\Tt,f)]_{\Ii_{f}} \in \Pp(\Ii_{f})$. Similarly, let
$(\Tt,f,c)$ be a flat/charged $\Ii$--triangulation for $\widetilde{X}$, with the associated element $[(\Tt,f,c)]_{\Ii_{fc}}
\in \Pp(\Ii_{fc})$. Here is the main result of this section:

\begin{teo}\label{scissorsclass} $\phantom{99}$

{\rm(1)}\qua The class $[(\Tt,f)]_{\Ii_f} \in \Pp(\Ii_f)$ does not depend on
the choice of $(\Tt,f)$. Hence $\cG_{\Ii_f} (X) = [(\Tt,f)]_{\Ii_f}$ is a well
defined invariant of $X$, called its $\Ii_f$--{\rm scissors congruence
  class}.

\noindent 
{\rm(2)}\qua The class $[(\Tt,f,c)]_{\Ii_{fc}} \in \Pp(\Ii_{fc})$ does not
depend on the choice of $(\Tt,f,c)$. Hence $\cG_{\Ii_{fc}} (X) =
[(\Tt,f,c)]_{\Ii_{fc}}$ is a well defined invariant of $X$, called its
$\Ii_{fc}$--{\rm scissors congruence class}.
\end{teo}
{\bf Proof}\qua By Theorem \ref{connessione}, any two
$\Aa$--triangulations of a same $3$--dimensional closed polyhedron can
be connected by means of {\it arbitrary} $\Aa$--transits, including
bubble and $0\to 2$ transits.  In particular, we can realize via
transits arbitrary global branching changes. (Without invoking this
fact we would get, in an easier way, the weaker result that the class
is well defined in $\Pp'(\Aa)$). Only the five term relations
associated to $2 \leftrightarrow 3$ $\Aa$--transits occur in the
definition of $\Pp(\Aa)$. Hence we have to prove that the relations
induced by all instances of the other transits follow from the five
term ones. This is done in the following two lemmas (compare with
Lemma \ref{consmove}). For simplicity we restrict the first to
$\Pp[\Dd]$.

Given a $\Dd$--tetrahedron $(\Delta,b,z)$, consider the tetrahedron
$\Delta !$ obtained by deforming an edge of $\Delta$ until it passes
through the opposite edge. We still denote $b$ and $z$ the branching
and cocycle induced on $\Delta !$. Note that we can identify
$\Delta$ and $\Delta !$ as \emph{bare} tetrahedra, since one is
obtained from the other by a cellular self-homeomorphism.

\begin{lem} \label{zerodeux} The following relation holds in
  $\Pp(\Dd)$: $(\Delta,b,z) = (\Delta ! ,b,z)$.
\end{lem}

{\bf Proof}\qua Let us prove a particular instance of this
relation. All the others come in exactly the same way, for there is no
restriction on the specific branching we choose in the arguments
below. Our arguments are based on a pictorial encoding with decorated
tetrahedra, but this is no loss of generality since the corresponding
algebraic relations in $\Pp(\Dd)$ may be thought as between abstract
elements.

\begin{figure}[ht]
\begin{center}
\begin{picture}(0,0)%
\includegraphics{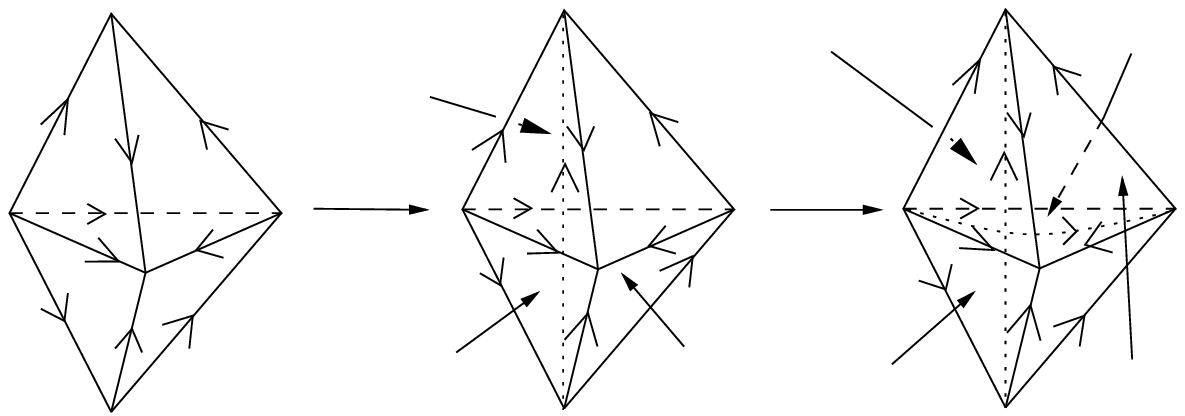}%
\end{picture}%
\setlength{\unitlength}{4144sp}%
\begingroup\makeatletter\ifx\SetFigFont\undefined%
\gdef\SetFigFont#1#2#3#4#5{%
  \reset@font\fontsize{#1}{#2pt}%
  \fontfamily{#3}\fontseries{#4}\fontshape{#5}%
  \selectfont}%
\fi\endgroup%
\begin{picture}(5627,1864)(6084,-4933)
\put(8216,-4757){\makebox(0,0)[lb]{\smash{\SetFigFont{6}{7.2}{\familydefault}{\mddefault}{\updefault}$(\Delta_3,\Dd_3)$}}}
\put(9281,-4703){\makebox(0,0)[lb]{\smash{\SetFigFont{6}{7.2}{\familydefault}{\mddefault}{\updefault}$(\Delta_4,\Dd_4)$}}}
\put(6084,-3401){\makebox(0,0)[lb]{\smash{\SetFigFont{5}{6.0}{\familydefault}{\mddefault}{\updefault}$(\Delta_1,\Dd_1)$}}}
\put(6176,-4757){\makebox(0,0)[lb]{\smash{\SetFigFont{5}{6.0}{\familydefault}{\mddefault}{\updefault}$(\Delta_2,\Dd_2)$}}}
\put(10208,-4813){\makebox(0,0)[lb]{\smash{\SetFigFont{6}{7.2}{\familydefault}{\mddefault}{\updefault}$(\Delta_6,\Dd_6)$}}}
\put(9914,-3201){\makebox(0,0)[lb]{\smash{\SetFigFont{5}{6.0}{\familydefault}{\mddefault}{\updefault}$(\Delta_5,\Dd_5)$}}}
\put(7973,-3410){\makebox(0,0)[lb]{\smash{\SetFigFont{6}{7.2}{\familydefault}{\mddefault}{\updefault}$(\Delta_5,\Dd_5)$}}}
\put(11245,-3183){\makebox(0,0)[lb]{\smash{\SetFigFont{5}{6.0}{\familydefault}{\mddefault}{\updefault}$(\Delta_8,\Dd_8)$}}}
\put(11216,-4796){\makebox(0,0)[lb]{\smash{\SetFigFont{6}{7.2}{\familydefault}{\mddefault}{\updefault}$(\Delta_7,\Dd_7)$}}}
\end{picture}
\end{center}
\caption{\label{fig1scissors} How to produce two term relations}
\end{figure}

\begin{figure}[ht]
\begin{center}
\begin{picture}(0,0)%
\includegraphics{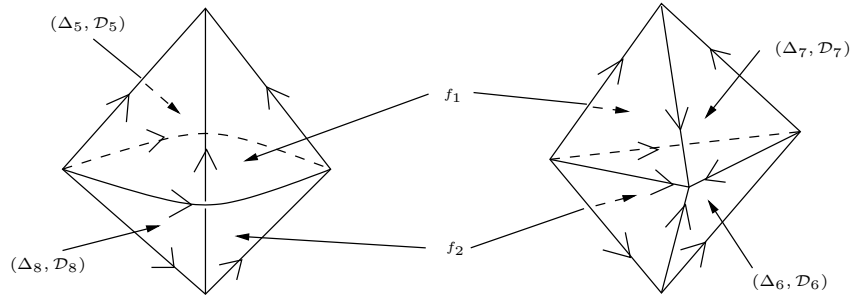}%
\end{picture}%
\setlength{\unitlength}{4144sp}%
\begingroup\makeatletter\ifx\SetFigFont\undefined%
\gdef\SetFigFont#1#2#3#4#5{%
  \reset@font\fontsize{#1}{#2pt}%
  \fontfamily{#3}\fontseries{#4}\fontshape{#5}%
  \selectfont}%
\fi\endgroup%
\begin{picture}(4749,1764)(5056,-4853)
\put(9496,-4786){\makebox(0,0)[lb]{\smash{\SetFigFont{6}{7.2}{\familydefault}{\mddefault}{\updefault}$(\Delta_6,\Dd_6)$}}}
\put(9631,-3391){\makebox(0,0)[lb]{\smash{\SetFigFont{6}{7.2}{\familydefault}{\mddefault}{\updefault}$(\Delta_7,\Dd_7)$}}}
\put(7651,-3661){\makebox(0,0)[lb]{\smash{\SetFigFont{5}{6.0}{\familydefault}{\mddefault}{\updefault}$f_1$}}}
\put(5056,-4681){\makebox(0,0)[lb]{\smash{\SetFigFont{6}{7.2}{\familydefault}{\mddefault}{\updefault}$(\Delta_8,\Dd_8)$}}}
\put(7666,-4576){\makebox(0,0)[lb]{\smash{\SetFigFont{5}{6.0}{\familydefault}{\mddefault}{\updefault}$f_2$}}}
\put(5296,-3241){\makebox(0,0)[lb]{\smash{\SetFigFont{6}{7.2}{\familydefault}{\mddefault}{\updefault}$(\Delta_5,\Dd_5)$}}}
\end{picture}
\end{center}
\caption{\label{fig2scissors} Decomposition of the final configuration in
Figure \ref{fig1scissors}}
\end{figure}

Consider the sequence of $2 \to 3$ $\Dd$--transits in
 Figure \ref{fig1scissors}, starting with an arbitrary cocycle
 transit. Denote by $\Dd_i=(b_i,z_i)$ the decoration of $\Delta_i$. We
 call $(\Delta_5,\Dd_5)$ and $(\Delta_8,\Dd_8)$ the two
 $\Dd$--tetrahedra glued along two faces in the final configuration
 (see Figure \ref{fig2scissors}); $(\Delta_8,\Dd_8)$ is glued to
 $(\Delta_6,\Dd_6)$ and $(\Delta_7,\Dd_7)$ along $f_1$ and $f_2$.

\noindent
The branchings and the cocycles of $(\Delta_6,\Dd_6)$ and
$(\Delta_7,\Dd_7)$ are respectively the same as for $(\Delta_1,\Dd_1)$
and $(\Delta_2,\Dd_2)$. Hence we may identify $(\Delta_6,\Dd_6)$ with
$(\Delta_1,\Dd_1)$ and $(\Delta_7,\Dd_7)$ with $(\Delta_2,\Dd_2)$. As
a consequence of the five term relations in $\Pp(\Dd)$, we deduce that
the composition of transits in Figure \ref{fig1scissors} translates into
$$\begin{array}{lll} (\Delta_1,\Dd_1) + (\Delta_2,\Dd_2) & = &
(\Delta_5,\Dd_5) - (\Delta_8,\Dd_8) + (\Delta_6,\Dd_6) +
(\Delta_7,\Dd_7) \\ & = & (\Delta_5,\Dd_5) - (\Delta_8,\Dd_8) +
(\Delta_1,\Dd_1) + (\Delta_2,\Dd_2)
\end{array}$$
where we note that $\Dd_8$ gives a negative orientation to
$\Delta_8$. This yields $(\Delta_5,\Dd_5) - (\Delta_8,\Dd_8) = 0$. As
the mirror image of $(\Delta_5,\Dd_5)$ is $(\Delta_8,\Dd_8)$, this proves the
lemma. \endproof

\begin{cor}\label{23full} 
The relations in $\mz[\Aa]$ corresponding to the $0 \leftrightarrow 2$
and bubble $\Aa$--transits are consequences of the relations
corresponding to the $2 \leftrightarrow 3$ $\Aa$--transits.
\end{cor}

{\bf Proof}\qua A $0 \leftrightarrow 2$ $\Dd$--transit leads to
mirror $\Dd$--tetrahedra $(\Delta_1,\Dd_1)$ and\break $(\Delta_2,\Dd_2)$. So the conclusion
for $0 \leftrightarrow 2$ $\Dd$--transits follows from Lemma
\ref{zerodeux}.

If we cut open the final configuration of a bubble
$\Dd$--transit along an interior face, we obtain the final
configuration of a $0 \to 2$ $\Dd$--transit. So, except for
what concerns the gluings of the two involved $\Dd$--tetrahedra, a
bubble $\Dd$--transit is abstractly given by the same data as a $0 \to
2$ $\Dd$--transit, and the two term relations in $\mz[\Dd]$ associated
to the $0 \leftrightarrow 2$ $\Dd$--transits and the bubble
$\Dd$--transits are the same.

By using Proposition \ref{idealmor} we deduce the result
for $\Pp[\Ii]$. The very same ``mirror argument'' of Lemma
\ref{zerodeux} lifts so as to eventually imply the other cases
$\Aa=\Ii_f$ and $\Aa=\Ii_{fc}$: the only difference is that there is
one degree of freedom for choosing the flattenings or the charges
during a flat/charged $2 \rightarrow 3$ transit. This causes no
problem because, given the first flattening or charge transit in Figure
\ref{fig1scissors}, we can always choose the second so as to produce,
again, mirror images in the final configuration. Also, we know from
the discussion before Proposition \ref{DdomI} that the bubble charge
transits depend on the choice of two interior edges in the final
configuration, where the sum of the charges is equal to $0$. As above,
if we cut it open along the interior face enclosed by these two edges,
we obtain the final configuration of a $0 \to 2$ charge transit. This
concludes the proof.  \endproof

It follows from Theorem \ref{scissorsclass} (1) and Proposition
\ref{blochfun} that

\begin{cor}\label{rogerofclass} We have $H_1(X) = H_1(\cG_{\Ii_f}(X))$.
\end{cor}

This means that the classical dilogarithmic invariants coincide with
the values of a function on $\Pp(\Ii_f)$ at the points corresponding
to the scissors congruence classes. As mentioned in Section \ref{CSV},
this function can be identified with the universal second
Cheeger--Chern--Simons class for flat $PSL(2,\mc)$--bundles.

A similar result, (by using $\Pp(\Ii_{fc})$ instead of $\Pp(\Ii_f)$) is
hopeless for the quantum dilogarithmic invariants $H_N$, $N>1$,
because the formal sums of tetrahedra that represent the points of the
{\it Abelian} group $\Pp(\Ii_{fc})$ do not encode any information
about $2$--face identifications, which, on the contrary, are essential
in the definition of the state sums. A way to overcome this problem
consists in defining an ``augmented'' scissors congruence
class, belonging to a Bloch-like group of further enriched
$\Ii_{fc}$--tetrahedra, where also the {\it states} are incorporated in
the augmented decorations; such a procedure is described in Section 5
of \cite{BB1}. However, this appears purely formal and risks hiding
more substantial questions. For instance, it makes sense to ask
whether the value of a quantum dilogaritmic state sum for a
$3$--manifold $\Ii_{fc}$--triangulation does actually only depend on the
corresponding $\Ii_{fc}$--scissors congruence class, though a positive answer
would be very surprising. More precisely:

\begin{ques}\label{quest1} {\rm Let $(W_j,L_j,\rho_j)$, $j=1,2$, be
    triples with the same $\Ii_{fc}$--scissors congruence classes:
    $\cG_{\Ii_{fc}}(W_1,L_1,\rho_1) = \cG_{\Ii_{fc}}(W_2,L_2,\rho_2)$.
    Do we have $$H_N(W_1,L_1,\rho_1) = H_N(W_2,L_2,\rho_2)$$ and
    similarly after replacing one $(W_j,L_j,\rho_j)$ or both with
    cusped manifolds ?}
\end{ques}

On the other hand, we proposed in Section 5 of \cite{BB2} a `Volume
Conjecture' for cusped manifolds $M$, saying that the dominant term of
the asymptotic expansion, when $N\to \infty$, of $H_N(M)^N$ grows
exponentially with $N^2$, and has a growth rate equal to $H_1(M)=
\exp\bigl((1/i\pi){\bf R}(M)\bigr)$ (see Subsection \ref{CSV} for
details on ${\bf R}$). In fact, the existence of the invariants
$H_N(M)$, $N \geq 1$, was only conjectural in that paper. Now this
conjecture is perfectly consistent for weakly-gentle cusped
manifolds. It is mostly motivated by the strong structural
coincidence between the classical and quantum invariants of such
manifolds (see Remark \ref{nodis}). 

With similar motivations, we also proposed a Volume Conjecture for
sequences $(W_n,L_n,\rho_n)$ of compact hyperbolic $3$--manifolds
converging geometrically to a cusped manifold $M$ (here, the $L_n$ are
the links isotopic to short simple closed geodesics that disappear in
the limit opening some cusps, and the $\rho_n$ are the holonomies of
the hyperbolic structures on the $W_n$). By Corollary
\ref{rogerofclass}, on the path to these conjectures there is the
seemingly weaker conjecture:
\begin{conj}\label{scissorvolconj} For every weakly-gentle cusped manifold $M$, 
  the dominant term of the asymptotic expansion of $H_N(M)$ when $N\to
  + \infty$ only depends on the scissors congruence class
  $\cG_{\Ii_{fc}}(M)$. Similarly, for the dominant term of
  $H_N(W_n,L_n,\rho_n)$ when both $N$ and $n$ tend to $+\infty$.
\end{conj}

Note that the asymptotic behaviour of $H_N(W,L,\rho)$ actually
depends, in general, on the link $L$ (see again Section 5 of \cite{BB2}),
whereas $H_1(W,\rho)$ only depends on $\cG_{\Ii_f}(W,\rho)$, which,
for any link $L$ in $W$, is the image of $\cG_{\Ii_{fc}}(W,L,\rho)$
via the natural forgetful map.

\small

\section{Appendix}\label{2GTAPP} 

In Section \ref{algprel} we recall some results on Heisenberg doubles.
The detailed constructions are given in \cite{B}, Chapter 3, but some
of them were already announced by Kashaev in \cite{K-}. This serves in
Section \ref{repR} to give a representation theoretic formulation of
the tensor $\widehat{\Ll}_N$ defined in Section \ref{CQDSYMQUANTUM}.

\subsection{Algebraic preliminaries}\label{algprel}

\noindent 
{\bf Heisenberg doubles.} Let $A = (1,\epsilon,m,\Delta,S)$ be a Hopf
algebra with unity over a ring $k$, where $1$, $\epsilon$, $m$,
$\Delta$ and $S$ are respectively the unit, the counit, the
multiplication, the comultiplication and the antipode of $A$. When $A$
is infinite dimensional we put $k=\mc[[h]]$, the ring of formal power
series over $\mc$ with indeterminate $h$, and we assume that $A =
U_h(\mathfrak{g})$, the quantum universal enveloping algebra (QUE) of
a complex finite dimensional Lie algebra $\mathfrak{g}$ (see
eg \cite[\S 6-8]{CP}). Every tensor product below is over $k$, and
when $A = U_h(\mathfrak{g})$ they are implicitely completed in
$h$--adic topology. When $A$ is finite dimensional the dual $k$--module
$A^* = {\rm Hom}(A,k)$ naturally inherits a {\it dual} Hopf algebra
structure $A^* = (\epsilon^*,1^*,\Delta^*,m^*,S^*)$ with
$$\begin{array}{l} \langle x,m(a \otimes b)\rangle = \langle m^*(x),a
\otimes b\rangle\quad ,\quad \langle x \otimes y,\Delta(a)\rangle =
\langle \Delta^*(x \otimes y),a \rangle\\ \\ \langle x,1\rangle =
1^*(x)\ ,\quad \langle \epsilon^*,a\rangle = \epsilon(a)\quad ,\quad
\langle x,S(a)\rangle = \langle S^*(x),a\rangle
\end{array}$$
\noindent where $a$, $b \in A$, $x$, $y \in A^*$, and $\langle \ ,\
\rangle \co  A^* \otimes A \rightarrow k$ is the canonical pairing. When
$A=U_h(\mathfrak{g})$ we have the Drinfeld's notion of QUE-\emph{dual}
Hopf algebra $A^* = U_h^*(\mathfrak{g})$ \cite[\S 6.3.D-8.3]{CP}. This
is a QUE-algebra isomorphic to $U_h(\mathfrak{g})$ as a
$\mc[[h]]$--module, with a dual Hopf algebra structure defined as
above.

For any $a \in A$, denote by ${\rm ev}_a \co  A^* \rightarrow
k$ the evaluation map: ${\rm ev}_a(x) = \langle x,a \rangle$. Let
$\pi_A\co  A \rightarrow {\rm End}_k(A^*)$ be the homomorphism defined by
$\pi_A(a) = (id \otimes {\rm ev}_a) \ m^*$, and $\pi_{A^*}\co  A^*
\rightarrow {\rm End}_k(A^*)$ be given by multiplication on the
left. The {\it Heisenberg double} $\mathcal{H}(A)$ of $A$ is the
subalgebra of ${\rm End}_k(A^*)$ (topologically) generated by the
image of $\pi_A$ and $\pi_{A^*}$ (this notion seems to have been
introduced in \cite{BS} and \cite{STS}). The image $R$ in $A \otimes
A^*\cong {\rm End}_k(A)$ of the identity morphism defines an
automorphism
$$(\pi_A \otimes \pi_{A^*})(R) = (\epsilon^* \otimes \Delta^*)\ (m^*
\otimes \epsilon^*)$$ of $(A^*)^{\otimes 2}$. We say that $R$ is the
{\it canonical element} of $\mathcal{H}(A)$. Given dual (topological)
basis $\{e_\alpha\}_{\alpha}$ and $\{e^{\beta}\}_{\beta}$ of $A$ and
$A^*$ for the pairing $\langle \ ,\ \rangle$, we write $\textstyle R =
\sum_{\gamma} e_\gamma \otimes e^\gamma$. Viewed as an element of
${\rm End}_k((A^*)^{\otimes 2})$, it has the following remarkable
(equivalent) properties:
\begin{eqnarray} \label{pentagone}
(1 \otimes e_{\alpha}) \ R = R \ \Delta(e_\alpha) \quad ,\quad R \
(e^{\alpha} \otimes 1) = m^*(e^\alpha) \ R \nonumber \\ R_{12} \
R_{13}\ R_{23} = R_{23}\ R_{12}\ \hspace{2cm}
\end{eqnarray}
\noindent where $R_{12} = R \otimes 1,\ R_{23} = 1 \otimes R$,
etc. The last identity is called the {\it pentagon} relation. The
first identity shows that for any linear representation $\rho$ of $A$
in the second tensor factor of the left hand side, $R^{-1}$ induces
embeddings of $\rho$ into the tensor product of two linear
representations of $A$ described via $\Delta(\rho)$, that is {\it
Clebsch-Gordan operators}. The pentagon relation means that $R$ (here
$R_{23}$ in the left hand side) induces matrices of change of basis
between the two possible ways of composing such embeddings, that is
{\it matrices of $6j$--symbols}. Finally, let us note that we can
reconstruct completely $\mathcal{H}(A)$ from the relations
(\ref{pentagone}). In fact, any solution of the pentagon equation
uniquely corresponds to a Hopf module over some Hopf algebra
\cite[Theorem 4.10]{BS}, \cite[Theorem 5.7]{Dav}.

\medskip

{\bf The case of $U_h(b(2,\mc))$}\qua Let us now specialize the
above considerations to the (positive) quantum Borel subalgebra
$A=U_h(b(2,\mc))$ of $U_h(sl(2,\mc))$. Recall that this is the QUE
Hopf algebra over $\mc[[h]]$ topologically generated by elements $H$
and $D$ such that $HD -DH=D$, with comultiplication and antipode given
by \cite[\S 6.4]{CP}, \cite[\S 17]{Kas}:
$$\begin{array}{l} \Delta(H) = H \otimes 1 + 1 \otimes H\quad , \quad
\Delta(D) = 1 \otimes D + D \otimes e^{hH}\\ S(H) = -H \quad , \quad
S(D) = -De^{-hH}\quad ,\quad \epsilon(H) = \epsilon(D) = 0\quad ,
\quad \epsilon(1) = 1.
\end{array}$$
(We denote by $1$ the identity of $U_h(b(2,\mc))$ and $\mc[[h]]$). For
technical reasons related to some of our choices below, let us
introduce $\mc((h))$, the field of fractions of $\mc[[h]]$, and
consider $U_h(sl(2,\mc))$ as a $\mc((h))$--module. The QUE-dual Hopf
algebra $U_h^*(b(2,\mc))$ is isomorphic as a topological Hopf algebra
over $\mc[[h]]$ to the {\it negative} quantum Borel subalgebra of
$U_h(sl(2,\mc))$, endowed with the opposite comultiplication
\cite[Proposition 8.3.2]{CP}. Hence $U_h^*(b(2,\mc))$ is topologically
generated over $\mc[[h]]$ by elements $\bar{H}$ and $\bar{D}$ such
that $\bar{H}\bar{D} - \bar{D}\bar{H} = -h\bar{D}$, with
comultiplication and antipode given by:
$$\begin{array}{ll} \Delta(\bar{H}) = \bar{H} \otimes 1 + 1 \otimes
\bar{H}\quad , \quad \Delta(\bar{D}) = 1 \otimes \bar{D} + \bar{D}
\otimes e^{-\bar{H}}\\ S(\bar{H}) = -\bar{H} \quad , \quad S(\bar{D})
= -\bar{D}e^{\bar{H}}\quad ,\quad \epsilon(\bar{H}) =
\epsilon(\bar{D}) = 0\quad , \quad \epsilon(1) = 1.
\end{array}$$
\noindent Clearly the map $\bar{H} \rightarrow -hH$, $\bar{D}
\rightarrow D$ is an isomorphism of topological algebras over
$\mc((h))$. It is shown in \cite{B}, Proposition 3.2.5, that the {\it Heisenberg double}
$\mathcal{H}_h(b(2,\mc))$ of $U_h(sl(2,\mc))$ is isomorphic to the
$\mc((h))$--algebra topologically generated over $\mc[[h]]$ by elements
$H$, $D$, $\bar{H}$ and $\bar{D}$ such that:
$$\begin{array}{ll} HD - DH = D\quad , \quad \bar{H} \bar{D} - \bar{D}
\bar{H} = -h\bar{D} \\ H\bar{H} - \bar{H} H = 1\quad , \quad \
D\bar{H} = \bar{H} D\\ H\bar{D} - \bar{D} H = -\bar{D}\quad ,\quad
D\bar{D} - \bar{D} D = (1 - q) \ e^{hH} \end{array}$$ where we put $q
= e^{-h}$. Moreover, we can write the canonical element of
$\mathcal{H}_h(b(2,\mc))$ as
\begin{equation}\label{can1}
R_h = e^{H \otimes \bar{H}}\ (D \otimes \bar{D};q)_{\infty}^{-1}\quad
.
\end{equation}
Here we denote by $(x;q)_{\infty}$ the {\it $q$--dilogarithm}, which is
the formal power series in $\mc(q)[[x]]$ given by
$$(x;q)_{\infty} = \prod_{n=0}^{\infty} (1-xq^n) = \sum_{n=0}^{\infty}
\frac{q^{\frac{n(n-1)}{2}}(-x)^n}{(q)_n}$$ where $(q)_n = (1-q)\ldots
(1-q^n)$. The proof of this result is instructive. A remarkable fact
is that the term $e^{H \otimes \bar{H}}$ in (\ref{can1}) is the
canonical element of the Heisenberg double $\mathcal{H}_h^0$ of the
Hopf subalgebra of $\mathcal{H}_h(b(2,\mc))$ topologically generated
over $\mc[[h]]$ by $H$. Note that in $\mathcal{H}_h^0$ we have
$H\bar{H} - \bar{H}H=1$. Then, it is easy to see that the pentagon
relation (\ref{pentagone}) for $e^{H \otimes \bar{H}}$ is a direct
consequence of the Baker-Campbell-Hausdorff formula for the complex
Lie algebra generated by $H$ and $\bar{H}$. Moreover, the pentagon
relation for $R_h$ splits into the product of the one for $e^{H
\otimes \bar{H}}$, together with the following {\it $q$--dilogarithm
equation}:
\begin{equation}\label{canz}
(U;q)_{\infty} \ \left(\frac{[U,V]}{(1-q)};q\right)_{\infty}\
(V;q)_{\infty} = (V;q)_{\infty}\ (U;q)_{\infty}
\end{equation}
where $U=1 \otimes D \otimes \bar{D}$ and $V=D \otimes \bar{D} \otimes
1$.

\medskip

{\bf Integral form of $\mathcal{H}_h(b(2,\mc))$ at a root of
unity}\qua We need to specialize the formal parameter $q=\exp(-h)$ to
specific (non zero) complex numbers. For that we must consider the
{\it integral form} of $\mathcal{H}_h(b(2,\mc))$. As defined in
\cite{B}, Definition 3.2.9, this is the $\mc[q,q^{-1}]$--algebra
$\mathcal{H}_q(b(2,\mc))$ generated by elements $E$,$E^{-1}$,
$\bar{E}$, $\bar{E}^{-1}$, $D$ and $\bar{D}$ such that
$$\begin{array}{ll}
EE^{-1} = E^{-1}E = 1\\
DE = qED\quad , \quad \bar{D} \bar{E} = q\bar{E} \bar{D} \\
E\bar{E} = q \bar{E}E\quad ,\quad D\bar{E} = \bar{E} D\\
E\bar{D} = q \bar{D}E\quad ,\quad D\bar{D} - \bar{D} D = (1 - q) \ E .
\end{array}$$
This algebra is obtained just by mimicking the commutation relations
between the elements $D$, $\bar{D}$, $E=\exp(hH)$ and
$\bar{E}=\exp(-\bar{H})$ of $\mathcal{H}_h(b(2,\mc))$, as is usual in
quantum group theory \cite[\S 9]{CP}. The subalgebra of
$\mathcal{H}_q(b(2,\mc))$ generated by $E$, $E^{-1}$ and $D$ is
isomorphic to the integral form $\mathcal{B}_q$ of $U_h(b(2,\mc))$,
which has comultiplication, counit and antipode given by:
$$\begin{array}{l} \Delta(E) = E \otimes E\quad ,\quad \Delta(D) = E
\otimes D + D \otimes 1\\ \epsilon(E)=1\quad ,\quad \epsilon(D)=0\quad
,\quad S(E) = E^{-1}\quad ,\quad S(D) = -E^{-1}D.
\end{array}$$
So we write $\mathcal{H}(\mathcal{B}_q)=\mathcal{H}_q(b(2,\mc))$. For
any non zero complex number $\epsilon$, we can now evaluate $q$ in
$\epsilon$, thus giving the specialization
$\mathcal{H}(\mathcal{B}_{\epsilon})$ of $\mathcal{H}_q(b(2,\mc))$.

A main problem is that the canonical element $R_h \in
\mathcal{H}_h(b(2,\mc))$ does not survive this procedure. First
because $(D \otimes \bar{D};q)_{\infty}$ is an infinite sum, which
moreover is ill defined when $q$ is a root of unity. Also, $\exp(H
\otimes \bar{H})$ cannot be written in terms of the generators of
$\mathcal{H}(\mathcal{B}_q)$. However, $R_h$ acts by conjugation as an
automorphism of $\mathcal{H}(\mathcal{B}_q)^{\otimes 2}$. We
will construct, when $q=\zeta^{-1}$ is a root of unity, a specific
element $R_{\zeta}$ of a suitable extension of
$\mathcal{H}(\mathcal{B}_{\zeta^{-1}})^{\otimes 2}$, that implements
the regular part of this action of $R_h$. We do this in two steps,
that we describe below. We fix a primitive $N$th root of unity
$\epsilon=\zeta^{-1}$ for a positive integer $N > 1$.

\medskip

\noindent 
{\bf The action of $\exp(H \otimes \bar{H})$ at a root of unity}\qua We
have seen that $\exp(H \otimes \bar{H})$ is the canonical element of
$\mathcal{H}_h^0$. Now, the specialization in $\zeta^{-1}$ of the
integral form of $\mathcal{H}_h^0$ is the subalgebra of
$\mathcal{H}(\mathcal{B}_{\zeta^{-1}})$ generated by $E$ and $\bar{E}$
such that $\bar{E}E=\zeta E\bar{E}$. This subalgebra can be endowed
with the very same structure of Hopf algebra as for
$\mathcal{B}_{\zeta}$, so we denote it by
$\mathcal{B}_{\zeta}^0$. Moreover, it is a central extension of a
Heisenberg double. Indeed, we have an
isomorphism (see \cite{B}, Section 3.2.3)
$$\mathcal{B}_{\zeta}^0/(E^N=\bar{E}^N=1) \cong \mathcal{H}(\mc[\mz/N\mz])$$
where $\mathcal{H}(\mc[\mz/N\mz])$ is the Heisenberg double of the
group algebra $\mc[\mz/N\mz]$ of $\mz/N\mz$, endowed with its usual
Hopf algebra structure (see eg \cite[\S 3]{Kas}). Hence the
canonical element $S_N$ of $\mathcal{H}(\mc[\mz/N\mz])$ is the natural
`periodic' analogue of $\exp(H \otimes \bar{H})$. We can define the
image of $S_N$ in the following extension of
$(\mathcal{B}_{\zeta}^0)^{\otimes 2}$. The algebra
$\mathcal{B}_{\zeta}^0$ has no zero divisors, since it is an Ore
extension of a polynomial ring \cite[\S 4]{Kas}. Then its center
$\mathcal{Z}(\mathcal{B}_{\zeta}^0)$ is an integral domain and we can
consider its quotient ring
$Q(\mathcal{Z}(\mathcal{B}_{\zeta}^0))$. Put
\begin{equation}\label{corps}
Q(\mathcal{B}_{\zeta}^0) := \mathcal{B}_{\zeta}^0
\otimes_{\mathcal{Z}(\mathcal{B}_{\zeta}^0)} Q(
\mathcal{Z}(\mathcal{B}_{\zeta}^0))
\end{equation}
\noindent 
and let $c_{E}$ and $c_{\bar{E}}$ be elements of
$\mathcal{Z}(\mathcal{B}_{\zeta}^0)$ such that $c_{E}^N = E^N$ and
$c_{\bar{E}}^N = \bar{E}^N$ (they exist because
$\mathcal{Z}(\mathcal{B}_{\zeta}^0)$ is integrally closed
\cite[Proposition 11.1.2]{CP}). The image of
$S_N$ in $Q(\mathcal{B}_{\zeta}^0)^{\otimes 2}$ can be written as (see
\cite{B}, Lemme 3.2.10)
\begin{equation}\label{Szeta}
S_{\zeta}= \frac{1}{N}\ \sum_{i,j =0}^{N-1} \zeta^{ij}\ (E')^{i} \otimes (\bar{E}')^j
\end{equation} 
\noindent 
where $E' = c_{E}^{-1}E$, $\bar{E} ' = c_{\bar{E}}^{-1}\bar{E} \in
Q(\mathcal{B}_{\zeta}^0)$. Note that the sum
$\textstyle N^{-1} \sum_{i} \zeta^{ij}(\bar{E}')^j$ is the normalized
inverse Fourier transform of $\bar{E}'$. We can verify that the action
by conjugation of $\exp(H \otimes \bar{H})$ on
$\mathcal{H}(\mathcal{B}_{\zeta^{-1}})^{\otimes 2}$ is the same as the
action by conjugation of $S_{\zeta}$.

\medskip

\noindent 
{\bf The action of $(D \otimes \bar{D};q)_{\infty}$ at a root of
unity}\qua This is best seen by considering the behaviour of the
$q$--dilogarithm $(x;q)_\infty$ when $q \rightarrow \zeta^{-1}$. Note
that it is a solution of the $q$--difference equation $(1-x)f(qx) -
f(x) = 0$. When $\vert q \vert < 1$, $(x;q)_{\infty}$ converges
normally on any compact domain of $\mc$ and defines an entire function
$E_q(x)$ that may be written in the form of the infinite product
$\textstyle E_q(x) = \prod_{n=0}^{\infty} (1-xq^n)$. The zeros of
$E_q(x)$ are simple and span the set $\{q^{-n} ,\ n \geq 0\}$. When
$\vert q \vert > 1$, the radius of convergence of $(x;q)_\infty$ is
$\vert q \vert$, and its sum defines in the open disc of convergence a
holomorphic function $e_q(x)$. Moreover, we can continue $e_q(x)$
meromorphically to the whole complex plane, so that $e_q(x) =
E_{q^{-1}}(xq^{-1})^{-1}$. The function $e_q(x)$ has no zeros and its
poles are simple and span the set $\{q^{n}\ ,\ n \geq 1\}$.

The remarkable fact is that $(x;q)_{\infty}$ has essential
singularities when $q$ tends to roots of unity, as we explain now. Let
$q=\exp(-\varepsilon/N^2)\zeta^{-1}$, where ${\rm Re}(\varepsilon) >
0$. Recall from Section \ref{matdefsec} and Section \ref{GTBASICMD}
the definition of the Euler dilogarithm ${\rm Li}_2$ and the function
$g$. When $\vert x \vert <1$, the $q$--dilogarithm $(x;q)_\infty$ has
for $\epsilon \rightarrow 0$ the following asymptotic behaviour (see
eg \cite{BR}):
$$(x;q)_{\infty}= g^{-1}(x) \ \ (1-x^N)^{1/2} \ 
\exp(-{\rm Li}_2(x^N)/\epsilon) \ \ \bigl(1 +
\mathcal{O}(\epsilon)\bigr).$$ Replacing $g^{-1}$ with its power
series expansion at $x=0$, we deduce that the `regular' part of $(D
\otimes \bar{D};q)_{\infty}$ when $q \rightarrow \zeta^{-1}$ lies in
the vector space $\widehat{\mathcal{H}}(\mathcal{B}_{\zeta^{-1}})$ of
formal power series in the generators of
$\mathcal{H}(\mathcal{B}_{\zeta^{-1}})$, with complex coefficients.
Similarly to (\ref{corps}), define
$$Q(\widehat{\mathcal{H}}(\mathcal{B}_{\zeta^{-1}})) := 
\widehat{\mathcal{H}}(\mathcal{B}_{\zeta^{-1}}) 
\otimes_{\mathcal{Z}(\mathcal{B}_{\zeta}^0)} Q(\mathcal{Z}
(\mathcal{B}_{\zeta}^0)).$$
Then  
\begin{equation}\label{cancyc}
R_{\zeta}= S_{\zeta}\ \ g(D \otimes \bar{D})
\end{equation}
acts by conjugation on
$Q(\widehat{\mathcal{H}}(\mathcal{B}_{\zeta^{-1}}))^{\otimes 2}$ as
the regular part of $R_h$. Denote by $r$ the central element $(1-(D^N
\otimes \bar{D}^N))^{1/N}$, viewed as the evaluation at $D^N \otimes
\bar{D}^N$ of the power series expansion of $(1-x)^{1/N}$ at $x=0$. We
can prove that $r^{(1-N)/2}R_{\zeta}$ verifies the pentagon relation
(\ref{pentagone}) (the element $r$ serves as a `determinant-like'
normalization, see Proposition \ref{unitarity}). As for $R_h$, this
relation splits into the product of the pentagon relation for
$S_{\zeta}$, and an identity obtained from (\ref{canz}) by replacing
each $q$--dilogarithm with the evaluation of $r^{(1-N)/2}g$ at certain
multiples of $D \otimes \bar{D}$ by central elements. We will not
describe this matter here (see Remark \ref{remref}). Rather, we
consider below the five term relations induced on the `cyclic' finite
dimensional representations of
$\mathcal{H}(\mathcal{B}_{\zeta^{-1}})$.

\begin{remark} \label{remref}{\rm The function $g$ (respectively the matrix
    $\Psi$ defined by (\ref{formmat2a})) is a `cyclic' analogue of
    $(x;q)_{\infty}$, because Lemma \ref{gdilocyc} below (respectively
    (\ref{factor2})) is a version for $q=\zeta$ of the $q$--difference
    equation $(1-x)f(qx) -f(x)=0$, that defines $(x;q)_{\infty}$ when
    $f(0)=1$. So the matrix identity (\ref{FadKas}) is a cyclic
    version of (\ref{canz}). An equivalent form of it was first
    announced in \cite{FK}. See also \cite{BR} for a description
    starting from the formal asymptotics of $(x;q)_{\infty}$.}
\end{remark}

\subsection{The $\Ll_N$ are representations of $R_{\zeta}$}\label{repR}

First we need to introduce some properties of the function $g$. As
before, let $N>1$ be any odd positive integer, and denote by log the
standard branch of the logarithm, which has the imaginary part in
$]-\pi,\pi]$. For any complex number $x \ne 0$ write
$x^{1/N}=\exp((1/N)\log(x))$. As remarked in the Introduction, the
function
$$g(x) = \prod_{j=1}^{N-1} (1-x\zeta^{-j})^{j/N}$$ is complex analytic
on $\mathfrak{D}_N$, the complex plane with cuts from the points
$x=\zeta^k$, $k=1,\ldots,N-1$, to infinity.

\begin{lem}\label{gdilocyc} 
For any $x \in \mc \setminus 
\{\zeta^j,\ j=1,\ldots, N-1\}$ and $k \in \mz$ we have
$$g(x\zeta^k) = g(x)\ \prod_{j=1}^{k} 
\ \frac{\ (1-x^N)^{1/N}}{1-x\zeta^j}.$$
\end{lem}

{\bf Proof}\qua We have
$$g(x\zeta^k) = \prod_{j=1}^{N-1}(1-x\zeta^{k-j})^{k/N}\ 
\prod_{j=k+1}^{N-1}(1-x\zeta^{k-j})^{(j-k)/N}\ 
\prod_{j=1}^{k}(1-x\zeta^{k-j})^{(j-k)/N}$$
$$\begin{array}{l} = \frac{(1-x^N)^{k/N}}{(1-x\zeta^k)^{k/N}} \
 \prod_{l=1}^{N-k-1}(1-x\zeta^{-l})^{l/N}\ \left(\
 \prod_{l=N-k+1}^{N}(1-x\zeta^{-l})^{l/N} \
 \prod_{l=0}^{k-1}(1-x\zeta^l)^{-1}\ \right)\\ \\ = g(x) \
 (1-x^N)^{k/N} \ \prod_{l=1}^{k}(1-x\zeta^{l})^{-1} = \ g(x) \
 \prod_{j=1}^{k} \frac{\ (1-x^N)^{1/N}}{1-x\zeta^j} .
\end{array}$$
\noindent 
Note that in the second equality we used the fact that $\textstyle
\sum_{j=1}^{N} {\rm log}(1-x\zeta^j) = {\rm log}(1-x^N)$ if $\vert x
\vert < 1$ and, by analytic continuation, if $x$ belongs to the cut
complex plane $\mathfrak{D}_N$. In fact if $x$ lies in the interior of
such a half ray, then $x^N \in (1,+\infty)$, so that the imaginary parts are
corrected by the same amount on both sides. Hence we find the desired
result.\endproof

Denote by $S(x \vert z)$ the rational function defined on the curve
$\{x^N + z^N =1\}$ by
$$S(x \vert z) = \sum_{k=1}^{N} \ \prod_{j=1}^k \frac{z}{1-x\zeta^j}.$$
\begin{lem}\label{factor1} 
When both sides are defined the following identities hold true:

{\rm(i)} \ $g(x)\ g(1/x) \equiv_N \frac{g(1)^2}{N}\ 
x^{\frac{1-N}{2}}\ \frac{\ 1-x^N}{1-x}$

{\rm(ii)} \ $x\ S(x\vert z\zeta) = (1-z)\ S(x \vert z)$ 

{\rm(iii)} \ $g(x)\ g(z/\zeta)\ S(x \vert z)\equiv_N x^{N-1}\ g(1)$

\noindent 
where $\equiv_N$ denotes the equality up to multiplication by $N$th
roots of unity.
\end{lem}

{\bf Proof}\qua (i)\qua Note that 
$$g(x)^N\ g(1/x)^N = x^{\frac{N(1-N)}{2}}\  
\prod_{j=1}^{N-1} (1-x\zeta^{-j})^j\ (x-\zeta^{-j})^j.$$
\noindent 
The product in the right hand side is a polynomial in $x$ of degree
$N(N-1)$, with zeros the $N$th roots of unity $\zeta^j$, $j\ne 0$,
each with multiplicity $N$. So there is a function $C(x)$, constant up
to roots of unity, such that
$$g(x)\ g(1/x)= C(x)\ x^{\frac{(1-N)}{2}}\ \frac{\ 1-x^N}{1-x} .$$
\noindent 
We find $C(x)\equiv_N g(1)^2/N$ by taking the limit $x\rightarrow 1$.

(ii)\qua This is a simple computation:
$$x\ S(x \vert z\zeta) = \sum_{k=1}^N \quad \prod_{j=1}^k 
\frac{z}{1-x\zeta^j}\ (x\zeta^k -1 +1) = S(x \vert z) - zS(x \vert z).$$

(iii)\qua Consider the rational function $\textstyle Q(x) =
\prod_{j=0}^{N-1} S(x\vert z\zeta^j)$. It does only depend on $x$,
because $z^N=1-x^N$ and $Q$ is invariant under the substitution $z
\mapsto z \zeta$. Remark that
$$\begin{array}{lll}
S(x\zeta\vert z) & = & \sum_{k=1}^N \ \prod_{j=1}^k 
\frac{z}{1-x\zeta^{j+1}} = \sum_{k=1}^N \ \frac{1-x\zeta}{z} \ 
\prod_{j=1}^{k+1} \frac{z}{1-x\zeta^j} = S(x \vert z)\ 
\frac{1-x\zeta}{z}.
\end{array}$$
\noindent So we have
$$Q(x\zeta )\ = \ Q(x) \ \ \frac{(1-x\zeta)^N}{z^N} \ =\ Q(x) \ \
\frac{(1-x\zeta)^N}{1-x^N}.$$ This shows that $Q$ has no zeros,
and that its poles are the roots of unity $\{\zeta^j\}$,
$j=1,\ldots,N-1$, where $\zeta^j$ has multiplicity $j$. Hence we find
$$Q(x) = \frac{\prod_{j=1}^{N-1}(1-\zeta^j)^j}{\prod_{j=1}^{N-1}
(1-x^{-1}\zeta^j)^j} = g(1)^N\ \frac{x^{\frac{N(N-1)}{2}}}{g(x)^N}$$
\noindent 
where the normalization constant is found by taking $x=1$, which gives
$S(1\vert 0)=1$ and $Q(1)=1$. Moreover, by applying (ii) directly to
the formula for $Q$ we get
$$Q(x) = S(x \vert z)^N\ x^{\frac{-N(N-1)}{2}} \ \prod_{j=1}^{N-1} 
(1-(z/\zeta)\ \zeta^{-j})^j.$$
$$S(x \vert z)^N = g(1)^N\ \frac{x^{N(N-1)}}{g(z/\zeta)^Ng(x)^N}\leqno{\rm Then}$$
which is just the $N$th power of (iii).\endproof

Consider now a complex linear representation $\rho$ of
$\mathcal{H}(\mathcal{B}_{\zeta^{-1}})$ with support $V_{\rho}$, such
that: each generator is mapped into ${\rm GL}(V_{\rho})$, and
$\rho(\Bar{D})\rho(D) = \zeta \ \rho(D)\rho(\bar{D})$
(ie $\rho(\Bar{D}D)= -\rho(E)$). We say that $\rho$ is {\it
cyclic}. One can check that $V_{\rho}$ is necessarily finite
dimensional, with ${\rm dim}_{\mc}(V_{\rho})=N$ if $\rho$ is
irreducible. The elements of the center $\mathcal{Z}$ of
$\mathcal{H}(\mathcal{B}_{\zeta^{-1}})$ (in particular the $N$th
powers of the generators) act as scalars on $V_\rho$, so we
have homomorphisms $\chi_\rho\co  \mathcal{Z} \rightarrow \mc$ called the
\emph{central characters}. We note that $\rho$ induces a pair of
cyclic irreducible representations of $\mathcal{B}_{\zeta^{-1}}$, by
considering its restriction to the algebras generated by $E$, $D$ and
$\bar{E}$, $\bar{D}$ respectively. 

\noindent Recall from Section \ref{CQDSYMQUANTUM} that the matrix valued map
$\hat{\Ll}_N\co  \widehat{\mc} \rightarrow {\rm M}_{N^2}(\mc/U_N)$ is
complex analytic, where $U_N$ is the
multiplicative group of $N$th roots of unity. Recall also that
$\equiv_N$ denotes the equality up to multiplication by $N$th roots
of unity, and that the integer $m$ is defined by $N=2m+1$.

\begin{teo} \label{repRzeta-} 
For any $(u;p,q) \in \widehat{\mc}$, there exists a cyclic irreducible
linear representation $\rho$ of
$\mathcal{H}(\mathcal{B}_{\zeta^{-1}})$ on $\mc^N$ such that $\chi_\rho(D^N)\chi_\rho(\bar{D}^N)= 1-u^{-1}$, and 
$$\widehat{\Ll}_N(u;p,q) \equiv_N (\rho \otimes 
\rho)\left( r^{\frac{1-N}{2}}R_\zeta\ (D' \otimes \bar{D}')^{-p(m+1)}\right)$$
where $r$ and $R_\zeta$ are defined in (\ref{cancyc}), $D'=c_{D}^{-1}D$ 
and $\bar{D}'=c_{\bar{D}}^{-1}\bar{D}$ for central elements 
$c_{D}$, $c_{\bar{D}} \in \mathcal{H}(\mathcal{B}_{\zeta^{-1}})$ 
with $c_{D}^N=D^N$ and $c_{\bar{D}}^N=\bar{D}^N$, 
and $\epsilon \in \{-1,+1\}$ is defined by 
$\log(u) + \log(1/(1-u)) + \log(1-u^{-1}) = \epsilon \pi i$.
\end{teo}

We note that the factor $(D' \otimes \bar{D}')^{-p(m+1)}$ only serves
to get the $N$th root $v_{-q}'$ of $1-u$ in $\widehat{\Ll}_N(u;p,q) =
\Ll_N(u_p',v_{-q}')$. Putting $p=q=0$ we recover the basic matrix
dilogarithms $\Ll_N(u)$ defined in (\ref{LNBdil}). The
ambiguity up to $N$th roots of unity depends on the arguments of
$u_p'$ and $v_{-q}'$, and is due to the use of Lemma \ref{factor1}
(i) and (iii).

\medskip

{\bf Proof}\qua  Let $Z$, $X$ and $Y = \zeta^{(m+1)}XZ$
be the $N \times N$--matrices with entries defined as: $Z_i^j =
\zeta^i$ $\delta (i-j)$, $X_i^j = \delta (i-j-1)$, and $Y_i^j =
\zeta^{(m+1) +j}$ $\delta (i-j-1)$ in the canonical basis of
$\mc^N$, where $\delta$ is, as usual, the Kronecker symbol
mod($N$). Consider a cyclic irreducible representation $\rho$ of $\mathcal{H}(\mathcal{B}_{\zeta^{-1}})$ given by
$$\begin{array}{ll}
\rho(E) = t_\rho^{-2}\ Z^{-1}\quad ,& \rho(D) = 
-\frac{x_\rho}{t_\rho}\ Y^{-1}\\
\rho(\bar{E}) = s_\rho^{-2}\ Y\quad ,& \rho(\bar{D}) = 
\frac{1}{s_\rho y_\rho}\ Z^{-1}Y = \frac{\zeta^{-(m+1)}}{s_\rho 
y_\rho}\ X
\end{array}$$
\noindent 
for non zero complex numbers $t_\rho$, $s_\rho$, $x_\rho$ and
$y_\rho$. Choose the central elements $c_E$ and $c_{\bar{E}}$ in
(\ref{Szeta}) such that $\rho(c_E)=t_\rho^{-2}{\rm Id}_{\mc^N}$ and
$\rho(c_{\bar{E}})=s_\rho^{-2}{\rm Id}_{\mc^N}$. We see immediately
that
\begin{equation}\label{Szetarep}
(\rho \otimes \rho)(S_{\zeta})= \frac{1}{N}\ \sum_{i,j =0}^{N-1} 
\zeta^{ij}\ Z^{-i} \otimes Y^j.
\end{equation} 
\noindent 
Consider $(\rho \otimes \rho)(r^{\frac{1-N}{2}}g(D
\otimes \bar{D})\ (D' \otimes \bar{D}')^{-p(m+1)})$. As $g(x)$ is complex analytic on the open unit disk $\vert x \vert <1$, by
considering its power series expansion $\textstyle
g(x)=\sum_{k=0}^{\infty} a_kx^k$ at $x=0$ we see that for any
unipotent matrix $M$ of order $N$ we have
$$\begin{array}{lll} \sum_{k=0}^{\infty} a_k (xM)^k & = &
\sum_{k=0}^{\infty} \ \sum_{s=0}^{N-1}\left(\frac{1}{N}\
\sum_{t=0}^{N-1} \zeta^{t(k-s)}\right)\ a_k (xM)^k \\ & = &
\frac{1}{N} \sum_{s,t=0}^{N-1} \ (\sum_{k=0}^{\infty} a_k
x^k\zeta^{tk})\ \zeta^{-st}M^s
\end{array}$$
whence
\begin{equation}\label{devgmat}
g(xM) = \frac{1}{N} \sum_{s,t=0}^{N-1} g(x\zeta^t)\ \zeta^{-st}M^s .
\end{equation}
This identity may be analytically continued to the cut complex plane
$\mathfrak{D}_N$ with respect to the variable $x$. Now, set
$x=\chi_\rho(D^N)\chi_\rho(\bar{D}^N)$. Given $a \in \mz$ choose
$x_\rho$ and $y_\rho$ so that $x_a'= -x_\rho/t_\rho s_\rho y_\rho$.
Moreover, let the central elements $c_D$ and $c_{\bar{D}}$ verify
$\rho(c_D)=-\zeta^{-(m+1)}x_\rho/t_\rho$ and
$\rho(c_{\bar{D}})=1/s_\rho y_\rho$. We have $(\rho \otimes
\rho)(r)=y_0'$, where $y=1-x$, and
$$\begin{array}{l}
(\rho \otimes \rho)(r^{\frac{1-N}{2}}g(D \otimes 
\bar{D})\ (D' \otimes \bar{D}')^{-b(m+1)}) \hspace{5cm} \\ \\
\hspace{1.5cm} = (y_0')^{\frac{1-N}{2}}\ g(x_a'(Y^{-1} \otimes
Z^{-1}Y)) \ (\zeta^{(m+1)}Y^{-1} \otimes Z^{-1}Y)^{-b(m+1)}\\ \\
 \hspace{1.5cm} = (y_{-b}')^{\frac{1-N}{2}}\ \frac{1}{N} \
 \sum_{k,t=0}^{N-1} g(x_a'\zeta^k)\ \zeta^{-kt}\ (Y^{-1} \otimes
 Z^{-1}Y)^{t-b(m+1)} \\
\hspace{1.5cm} = (y_{-b}')^{\frac{1-N}{2}}g(x_a')\ \frac{1}{N} \
\sum_{k,t=0}^{N-1} \prod_{j=1}^k
\left(\frac{y_{-b}'}{1-x_a'\zeta^j}\right)\ \zeta^{-kt}\ (Y^{-1}
\otimes Z^{-1}Y)^t.
\end{array}$$
The second equality follows from (\ref{devgmat}) and the third from
Lemma \ref{gdilocyc} together with a sum reordering $t \mapsto t
-b(m+1)$. Moreover, by Lemma \ref{factor1} (ii) and (iii) we have
$$\sum_{k=0}^{N-1} \prod_{j=1}^k \left(\frac{\
y_{-b}'}{1-x_a'\zeta^j}\right)\ \zeta^{-kt} \equiv_N
\frac{g(1)(x_a')^{N-1}}{g(y_{-b}'\zeta^{-1})g(x_a')}\
\prod_{j=1}^{N-t} \frac{1-y_{-b}'\zeta^{j-1}}{x_a'}.$$ On
another hand, the lemmas \ref{gdilocyc} and \ref{factor1} (i) give
respectively
$$g(y_{-b}'\zeta^{-1}) = g(y_{-b}') \ \prod_{j=1}^{N-1}
\frac{x_0'}{1-y_{-b}'\zeta^j} = g(y_{-b}') \ \frac{1-y_{-b}'}{x_0'}$$
$$g(y_{-b}')\ g((y_{-b}')^{-1}) \equiv_N \frac{g(1)^2}{N}\
(y_{-b}')^{\frac{N-1}{2}}\ \frac{x}{1-y_{-b}'}.$$ Since
$$\prod_{j=1}^{N-t} \frac{1-y_{-b}'\zeta^{j-1}}{x_a'} =
\prod_{s=t+1}^{N} \frac{1-y_{-b}'\zeta^{-s}}{x_a'} = \prod_{s=1}^{t}
\frac{x_a'}{1-y_{-b}'\zeta^{-s}}$$ we find that $(\rho \otimes
\rho)(r^{\frac{1-N}{2}}\ g(D \otimes \bar{D})\ (D' \otimes
\bar{D}')^{-b(m+1)})$ is equal, up to multiplication by $N$th roots of
unity, to
\begin{equation}\label{formmat}
\frac{g((y_{-b}')^{-1})}{g(1)}\ \sum_{t=0}^{N-1} \prod_{s=1}^{t}
\frac{x_a'}{1-y_{-b}'\zeta^{-s}} \ (Y^{-1} \otimes Z^{-1}Y)^t.
\end{equation} 
It is easily seen that $(Y^{-1} \otimes Z^{-1}Y)^t\ _{k,l}^{i,j} =
\zeta^{-kt-t(t+1)(m+1)}\ \delta (l-j-t)\ \delta (k+t-i)$. Then
the $_{k,l}^{i,j}$--component of this sum is
$$\begin{array}{l} \frac{g((y_{-b}')^{-1})}{g(1)}\ \sum_{t=0}^{N-1}
\zeta^{-kt-t(t+1)(m+1)}\ \delta (l-j-t)\ \delta (k+t-i)\
\prod_{s=1}^{t} \frac{x_a'}{1-y_{-b}'\zeta^{-s}}\hspace{2cm}\\ \\
\hspace{1.5cm} = \frac{g((y_{-b}')^{-1})}{g(1)}\ \zeta^{-k(i-k)}\
\delta (l+k-i-j)\ \prod_{s=1}^{i-k}
\frac{x_a'\zeta^{-s}}{1-y_{-b}'\zeta^{-s}}\\ \\
\hspace{1.5cm} = \frac{g((y_{-b}')^{-1})}{g(1)}\ \zeta^{-k(i-k)}\
\delta (l+k-i-j)\ \prod_{s=1}^{i-k}
\frac{-x_a'(y_{-b}')^{-1}}{1-(y_{-b}')^{-1}\zeta^{s}}\\ \\
\hspace{1.5cm} = \frac{g((y_{-b}')^{-1})}{g(1)}\ \zeta^{-k(i-k)}\
\omega((y_{-b}')^{-1},(-x/y)_{a+b+\epsilon}'\vert i-k)\ \delta
(l+k-i-j)
\end{array}$$
\noindent 
where $\epsilon$ is defined as in the statement, replacing $u$ by
$x$. Finally, denote by $\{ e_j \}_j$ the canonical basis of
$\mc^N$. It follows from (\ref{Szetarep}) that $(\rho \otimes
\rho)(S_{\zeta}) (e_i \otimes e_j) = e_i \otimes Y^i e_j$. As
$(Y^i)_l^j = \zeta^{(m+1)i^2 + ij} \delta (l-i-j)$ we get
$$\begin{array}{l} (\rho \otimes \rho)(r^{\frac{1-N}{2}}\ R_\zeta\ (D'
\otimes \bar{D}')^{-b(m+1)})_{k,l}^{i,j} \hspace{5cm} \\ \\
\hspace{1cm} \equiv_N \frac{g((y_{-b}')^{-1}))}{g(1)}\ \sum_{s,t
=0}^{N-1} \zeta^{(m+1)s^2+st}\ \delta (k-s) \ \delta (l-s-t) \\
\hspace{3cm} \times \ \zeta^{-s(i-s)}\ \omega((y_{-b}')^{-1},
(-x/y)_{a+b+\epsilon}'\vert i-s) \ \delta (s+t-i-j) \\ \\
\hspace{1cm} \equiv_N \frac{g((y_{-b}')^{-1})}{g(1)}\
\zeta^{(m+1)k^2+kj} \
\omega((y_{-b}')^{-1},(-x/y)_{a+b+\epsilon}'\vert i-k) \ \delta
(i+j-l).
\end{array}$$
Renaming the variables as $u=1/y=1/(1-x)$, $p=b$ and
$q=-a-b-\epsilon$, this is exactly the formula (\ref{newL}) for
$\widehat{\Ll}_N(u;p,q)=\Ll_N(u_p',v_{-q}')$. \endproof

\begin{remark} \label{rempar} {\rm The isomorphism classes of cyclic
    irreducible representations of $\mathcal{B}_\zeta$ are in one-one
    correspondence with the elements of a dense subset of a Borel
    subgroup $B$ of $PSL(2,\mc)$; a parametrization is given by $[\nu]
    \mapsto (x_\nu^N,y_\nu^N)$, where $\chi_\nu(E^N)=x_\nu^{2N}$ and
    $\chi_\nu(D^N)=x_\nu^{N}y_{\nu}^N$. In particular, the
    representation of $\mathcal{B}_{\zeta^{-1}}$ defined by $\rho(E)$
    and $\rho(D)$ in Theorem \ref{repRzeta-} is isomorphic to the {\it
      dual} of the representation given by $(x_\nu^N,y_\nu^N)$. In
    fact, Theorem \ref{repRzeta-} and the discussion at the begining
    of Subsection \ref{algprel} show that the $\widehat{\Ll}_N$ are
    matrices of $6j$--symbols for the cyclic irreducible
    representations of $\mathcal{B}_\zeta$: they describe the
    associativity of the tensor product of such representations (see
    \cite{BB2}, Appendix). The matrix Schaeffer's identity for
    $\widehat{\Ll}_N$ (Theorem \ref{basNtr}) is a version of the
    pentagon relation satisfied by the $6j$--symbols.}
\end{remark}

We conclude this section with a `unitarity-like' property of
$\Ll_N(u_p',v_{-q}')$:

\begin{prop}\label{unitarity} 
The matrix $\Ll_N(u_p',v_{-q}')$ is invertible and
has determinant $1$. Moreover, we have
$$\Ll_N^{-1}(u_p',v_{-q}') \equiv_N U\left(\
\Ll_N((u_p')^*,(v_{-q}')^*)^T\ \right)^*U^{-1}$$ where $T$ is the
transposition, $*$ the complex conjugation, and $U$ is the symmetric
$N^2\times N^2$--matrix given by $U_{k,l}^{i,j}=\delta(k+i)\
\delta(l+j)$. In components:
$$\Ll_N^{-1}(u_p',v_{-q}')_{k, l}^{i,j} = [u_p']\ \frac{g(1)}{g(u_p')}\ \zeta^{-i l
-(m+1)i^2}\ (w(u_p'\zeta^{-1},v_{-q}' \vert k - i))^{-1}\
\delta(k+l-j).$$
\end{prop}

{\bf Proof}\qua By Theorem \ref{repRzeta-}, the invertibility
of $\Ll_N(u_p',v_{-q}')$ follows from that of $R_{\zeta}$ and $r$,
together with the fact that $(\rho \otimes \rho)(D' \otimes \bar{D}')$
is invertible (the representation $\rho$ is cyclic). We
compute the determinant as follows. Recall that $S_\zeta$ is defined
in (\ref{Szeta}). Since the matrices $Z$ and $Y$ are unipotent of
order $N$ (they satisfy $Z^N=Y^N={\rm Id}_{\mc^N}$), the
eigenvalues of $(\rho \otimes \rho)(S_\zeta)$ span the set of $N$th
roots of unity $\zeta^j$, $j=0,\ldots,N-1$, each with multiplicity
$N$. So ${\rm det}\left((\rho \otimes \rho)(S_\zeta)\right)=1$.
Moreover, with the notation of the proof of Theorem \ref{repRzeta-}, we
see that ${\rm det}((\rho \otimes \rho) (r^{\frac{1-N}{2}}g(D \otimes
\bar{D})))$ is equal to
$$\begin{array}{l} y^{\frac{N(1-N)}{2}}\ {\rm det}\left(g(x_a'(Y^{-1}
\otimes Z^{-1}Y))\right) = y^{\frac{N(1-N)}{2}}\ \prod_{j=0}^{N-1}
g(x_a'\zeta^j)^N \hspace{3.5cm} \\
 \hspace{2.5cm} = y^{\frac{N(1-N)}{2}}\ g(x_a')^N\ \prod_{j=1}^{N-1}
 \prod_{i=1}^{j} \frac{y}{(1-x_a'\zeta^{i})^N} \\
 \hspace{2.5cm} = y^{\frac{N(1-N)}{2}}\ \prod_{j=1}^{N-1}
 (1-x_a'\zeta^{-j})^{Nj}\
 \frac{y^{\frac{N(N-1)}{2}}}{\prod_{j=1}^{N-1}
 (1-x_a'\zeta^{j})^{N(N-j)}} = 1
\end{array}$$
where we use Lemma \ref{gdilocyc} in the second equality. As ${\rm
det}\left((\rho \otimes \rho)(D' \otimes \bar{D}')\right) = {\rm
det}(Y^{-1} \otimes Z^{-1}Y)) = 1$, by Theorem \ref{repRzeta-} we
eventually find ${\rm det}(\Ll_N(u_p',v_{-q}')) =1$.

As for the last property, note that we have $g(x^*)^*=
(-x)^{\frac{N-1}{2}} \zeta^{\frac{(N-1)(2N-1)}{6}}g(1/x)$. Let us
write $[x] = N^{-1} (1-x^N)/(1-x)$; then, Lemma \ref{factor1} $i)$ gives
$$\frac{g(x^*)^*}{g(1)^*} = \frac{g(x^*)^*}{g(1)^*}\
\frac{g(1)^*}{g(1)}\
(-1)^{\frac{N-1}{2}}\zeta^{-\frac{(N-1)(2N-1)}{6}} = \frac{g(1/x)\
x^{\frac{N-1}{2}}}{g(1)} \equiv_N [x]\ \frac{g(1)}{g(x)}.$$ On
another hand, we have
$$\begin{array}{lll} (w((u')^*,(v')^* \vert n))^* & = &
\prod_{i=1}^{n} \frac{v'}{1-u'\zeta^{-i}} = \prod_{i=n+1}^{N}
\frac{1-u'\zeta^{-i}}{v'} = \prod_{j=0}^{N-n-1}
\frac{1-u'\zeta^{j}}{v'} \\ & = & (w(u'\zeta^{-1},v' \vert -n))^{-1}
\end{array}$$
where, to simplify the notation, we write $u'$ and $v'$ for given
$N$th roots of $u$ and $v$. So we get (recall that $N=2m+1$):
$$\begin{array}{l} \left( U\left(\ \Ll_N((u')^*,(v')^*)^T\
\right)^*U^{-1} \right)_{k, l}^{i,j} =
(\Ll_N((u')^*,(v')^*)^*)_{-i,-j}^{-k,-l} \hspace{6cm}\\
\hspace{2.5cm} \equiv_N [u']\ \frac{g(1)}{g(u')}\ \zeta^{-i l
-(m+1)i^2}\ (w(u'\zeta^{-1},v' \vert k - i))^{-1}\
\delta(k+l-j).
\end{array}$$
We show that this is the inverse of $\Ll_N(u',v')$ as follows. First we have
$$\begin{array}{l} [u']\ \sum_{s,t=0}^{N-1} \zeta^{-sj
-s^2(m+1)}\ \zeta^{sl +(m+1)s^2}\ \frac{w(u',v' \vert k -
s)}{w(u'\zeta^{-1},v' \vert i - s)}\ \delta(i+j-t)\ \delta(k+l-t) \\
\stackrel{(*)}{=} [u']\delta(i+j-k-l)\zeta^{-i(j-l)}w(u',v' \vert k -
i)\sum_{s=0}^{N-1} \zeta^{(i-s)(j-l)}\ \frac{w(u'\zeta^{k-i},v' \vert
i - s)}{w(u'\zeta^{-1},v' \vert i - s)}.
\end{array}$$
Next we observe that the above sum is obtained from 
$$f(x,y\vert z) = \sum_{s=0}^{N-1} \ \prod_{j=1}^{s} \
\frac{1-y\zeta^j}{1-x\zeta^j}\ z^s$$ by setting $x=u'\zeta^{k-i}$,
$y=u'\zeta^{-1}$ and $z=\zeta^{j-l}$. Straightforward manipulations
very similar to that of Lemma \ref{factor1} $ii)$ imply that
$$\begin{array}{l}
f(x,y \vert z\zeta) = \frac{1-z}{x-yz\zeta}\ f(x,y \vert z) \\ \\
f(x\zeta,y \vert z) = \frac{(1-x\zeta)(x-yz)}{z(x-y)}\ f(x,y \vert z).
\end{array}$$
By iterating these two identities we find
$$\begin{array}{lll} f(x\zeta^t,x\zeta^{-1} \vert \zeta^s) & = &
x^{-[s-1]_N} \ \prod_{a=1}^{s-1}
\frac{1-\zeta^{s-a}}{\zeta^t-\zeta^{s-a}}\ f(x\zeta^t,x\zeta^{-1}
\vert \zeta) \\ & = & x^{-[s-1]_N}\ \prod_{a=1}^{s-1}
\frac{1-\zeta^{s-a}}{\zeta^t-\zeta^{s-a}}\ \delta(t)\ f(x,x\zeta^{-1}
\vert \zeta) \\ & = & x^{-[s-1]_N}\ \delta(t)\ f(x,x\zeta^{-1} \vert
\zeta)
\end{array}$$
where $[a]_N$ denotes the rest of the Euclidean division of $a \in
\mathbb{N}$ by $N$. Also, 
$$\begin{array}{lll} f(x,x\zeta^{-1} \vert \zeta) & = & 1 +
  \frac{(1-x)\zeta}{1-x\zeta} + \frac{(1-x)\zeta^2}{1-x\zeta^2} + \ldots
  + \frac{(1-x)\zeta^{N-1}}{1-x\zeta^{N-1}} \nonumber \\ & = & (1-x) \
  \biggl(\frac{1}{1-x} + \frac{\zeta}{1-x\zeta} + \ldots +
  \frac{\zeta^{N-1}}{1-x\zeta^{N-1}} \biggr) \nonumber\\ & = & (1-x) \
  \frac{d}{dx}\left( - \log (1-x^N) \right) = N x^{N-1}\frac{1-x}{\
    1-x^N}\ .\nonumber
\end{array}$$
Hence $f(u'\zeta^{k-i},u'\zeta^{-1} \vert \zeta^{j-l}) = x^{-[j-l-1]_N} \
\delta(k-i) \ (u')^{N-1}[u']^{-1}$, which shows that $(*)$ above is
equal to $\delta(k-i)\ \delta(j-l)$. This concludes the proof of the
last claim.\endproof

\end{document}